\documentclass[11pt]{article}
\usepackage{amsthm,amsmath,amssymb}
\usepackage[left=1.00in,bottom=1.25in]{geometry}
\usepackage{graphicx}
\usepackage{gensymb}
\usepackage{float}
\usepackage{url}
\usepackage{titling}
\usepackage{hyperref}
\usepackage{nicematrix}
\usepackage{comment}

\usepackage{titlesec}
\usepackage[labelformat=simple]{subcaption}

\usepackage{adjustbox}
\usepackage{tikz}
\usepackage{wrapfig}
\usepackage{dsfont}
\usepackage{color}
\usepackage{xcolor}
\usepackage{bm}
\usepackage{setspace}
\usepackage{array,multirow}
\usepackage{enumerate}

\definecolor{myalcolor}{rgb}{0.18, 0.75, 0.}

\definecolor{myamcolor}{rgb}{0.06,0.4,0.06}
\titleformat{\subsection}[runin]
  {\normalfont\large\bfseries}{\thesubsection}{1em}{}
\titlelabel{\thetitle.\quad}

\titleformat{\section}[runin]
  {\normalfont\large\bfseries}{\thesection}{1em}{}

\newtheorem*{remark}{Remark}
\newtheorem{theorem}{Theorem}[section]
\newtheorem{lemma}[theorem]{Lemma}

\newtheorem{corollary}[theorem]{Corollary}
\numberwithin{equation}{section}

\title{Convergence analysis and a novel Lagrange multiplier partitioned method for fluid-poroelastic interaction}
\author{Amy de Castro
\thanks{Department of Mathematics, University of Utah, Salt Lake City, UT 84123
    ({\tt amy.de.castro@utah.edu}). Partially supported by the NSF under grant number DMS-2207971.}
\and
Hyesuk Lee
  \thanks{School of Mathematical and Statistical Sciences, Clemson
        University, Clemson, SC  29634-0975 ({\tt hklee@clemson.edu}).
    Partially supported by the NSF under grant number DMS-2207971.}
}
\date{}

\begin{document}

\maketitle
\vspace{-10mm}
\begin{abstract}
We propose a partitioned method for the monolithic formulation of the Stokes-Biot system that incorporates Lagrange multipliers enforcing the interface conditions. The monolithic system is discretized using finite elements, and we establish convergence of the resulting approximation. A Schur complement based algorithm is developed together with an efficient preconditioner, enabling the fluid and poroelastic structure subproblems to be decoupled and solved independently at each time step. The Lagrange multipliers approximate the interface fluxes and act as Neumann boundary conditions for the subproblems, yielding parallel solution of the Stokes and Biot equations. Numerical experiments demonstrate the effectiveness of the proposed algorithm and validate the theoretical error estimate.

\end{abstract}

	\section{Background and Research Goals}

Poroelastic materials, modeled by the Biot equations, represent the deformation of an elastic structural skeleton saturated by a fluid governed by Darcy's law. Poroelastic models appear in when describing groundwater flow or flow through fissured rocks \cite{Detournay_1993, Murad_2001} or biomedical actions such as arterial or vascular blood flow and drug transport \cite{Banks_2017, Bociu_2021, Calo_2008, Causin_2014}. Most commonly, the Biot equations may be represented using only displacement and pore pressure variables, known as a two-field model, or an additional variable representing the Darcy velocity can be included to arrive at the three-field model. Previous well-posedness results for coupled fluid-poroelastic interaction (FPSI) systems include linear fluid models (Stokes equations) coupled to three-field Biot equations \cite{Showalter_2005} and fully inertial two-field Biot equations \cite{Avalos_2024}, and nonlinear fluid models (Navier-Stokes equations) \cite{Ambartsumyan_2019, Bociu_2021,Cesmelioglu_2017}. 

We consider the coupling of the fully dynamic two-field Biot equations with the dynamic Stokes equations, using three LMs to enforce interface conditions between the subdomains. The introduction of LMs facilitates the development of the domain decomposition strategy presented in this paper. Well-posedness of our formulation was demonstrated in \cite{deCastro_2025_FPSIWP}. Additional well-posedness results for LM formulations of fluid-poroelastic models include \cite{Ambartsumyan_2018, Caucao_2022, Li_2022_Hydro}; see \cite{deCastro_2025_FPSIWP} for further discussion.

Turning to the numerical analysis point of view, monolithic approaches to solving fully dynamic FPSI problems often rewrite the Biot equations as first-order in time; examples include a discontinuous Galerkin scheme \cite{Wen_2020} and a stabilized scheme using the fluid pressure Laplacian technique \cite{Cesmelioglu_2020}.  
A monolithic approach using the quasi-static Biot equations is developed in \cite{Ambartsumyan_2018}. Well-posedness and error analysis are established for this formulation, which uses a Lagrange multiplier to enforce conservation of mass at the interface.

Partitioned approaches to solving FPSI systems include an optimization-based decoupling method where the normal component of the interface stress is used as a control for least squares minimization \cite{Cesmelioglu_2016}; see \cite{Kunwar_2020} for a similar approach with a second-order time discretization. An operator splitting method for a multilayered structure consisting of thin elastic and thick poroelastic layers is implemented in \cite{Bukac_2015OpSplit}, based on a kinematically coupled scheme developed for FSI problems with thin structures \cite{Guidoboni_2009}. This method allows for decoupling of the fluid and Biot problems, and unlike many domain decomposition methods, it does not require subiterations between the subproblems. 

One successful class of iterative scheme for poroelasticity is a fixed-stress splitting scheme, which solves the flow problem with a fixed mean stress and subsequently updates the mechanical subproblem \cite{Both_2017, Wheeler_2007}. Fixed-stress and fixed-strain iterative approaches are compared in \cite{Kim_2011}; the fixed-stress scheme is found to be unconditionally stable while the fixed-strain is only conditionally stable. The quasistatic, three-field Biot models are most often under consideration in these methods.

In \cite{Bukac_2015Nitsche}, a loosely coupled approach based on Nitsche's coupling is examined. Nitsche's method weakly enforces the interface conditions; three subproblems are solved sequentially corresponding to the structure, Darcy components, and Stokes variables. Stability is established with the addition of several stabilization operators in the weak formulation. However, the method does introduce splitting errors, which reduce the rate of convergence. For this reason, the authors consider using this splitting method as a preconditioner for a monolithic scheme. See \cite{Burman_2014NitscheRR} for a comparison of Nitsche splitting with Robin--Robin coupling methods for FSI. 

Other partitioned approaches include a method stabilized by a splitting technique based on a variational multiscale approach  \cite{Badia_2009_NSBiot}, and a formulation based on a non-conforming Crouzeix–Raviart (C-R) element discretization \cite{Wilfrid_2020}. 
Many methods use Lagrange multipliers (LMs) to implement interface conditions. Often, LMs are used to represent the structure velocity and pore pressure on the interface \cite{Li_2024, Li_2022}, but an LM for Stokes velocity is also used in \cite{Caucao_2022}.

In our approach, we consider the fully dynamic two-field Biot equations coupled with the dynamic Stokes equations, using LMs to enforce the interface conditions between the subdomains.
Rather than rewriting the Biot equations as a second-order system in time by introducing an additional variable representing the structural velocity, we work directly with the fully discrete system based on a first order in time formulation.  The normal flux components appearing in the interface conditions are also represented by the LMs, allowing the monolithic coupled system to be reduced to a Schur complement equation at the fully discrete level.
As noted above, many existing methods employ the three-field or quasistatic Biot formulations, and thus our particular formulation is unique. 
The monolithic system involving the LMs has been analyzed in our previous work \cite{deCastro_2025_FPSIWP}, where we established the well-posedness at both the semi-discrete and fully discrete levels and derived the stability estimate.

In this work, we develop a strongly coupled, non-iterative partitioned method centered on a Schur complement equation for the solution of this FPSI problem.  An efficient strategy for solving the resulting Schur complement equation is presented, along with the development of an effective preconditioner.
Our work is outlined as follows. We discuss the model equations in Section \ref{FPSI_FEM:sec:Continuous Model-3LMs} and present the weak form, leading to a well-posed semi-discrete formulation. In Section \ref{FPSI_FEM:sec:errorDiscrete}, we turn to the fully discrete model and present convergence results with respect to spatial discretization. Our partitioned Schur complement method is described in Section
\ref{FPSI_FEM:sec:matrices-partitioned}, along with a discussion of efficient computation of matrix--vector products with the Schur complement matrix, and we propose a preconditioner for the Schur complement equation. Lastly, we verify the performance of our algorithm in Section \ref{FPSI_FEM:sec:Numerical}.

	\section{Model Equations and Semi-Discrete Model }\label{FPSI_FEM:sec:Continuous Model-3LMs}
% This is the version you want to use
We begin by developing the weak form and semi-discrete monolithic formulation of the Stokes-Biot system. Consider a Lipschitz, polytopal domain $\Omega$ which is divided into two open regions: $\Omega_f \in \mathbb{R}^d$ containing the fluid, and $\Omega_p \in \mathbb{R}^d$ containing the poroelastic structure, for $d=2,3$. We assume these domains are non-overlapping and share an interface $\gamma$. 
As the problem is time dependent, we take $T>0$ to be a given final time.

The flow in the fluid domain is modeled by the transient Stokes equations, and the two-field Biot equations describe the poroelastic material in $\Omega_p$. The resulting unknown functions are the fluid velocity $\bm{u}(\bm{x},t)$, the fluid pressure $p_f(\bm{x},t)$, the structural displacement $\bm{\eta}(\bm{x},t)$, and the pore pressure $p_p(\bm{x},t)$. With given body forces $\bm{f}_f, \bm{f}_\eta,$ and source or sink $f_p$, the model problem reads:

\noindent \textit{Find  $\bm{u} \in \Omega_f \times (0,T] \mapsto \mathbb{R}^d, \hspace{1mm} p_f \in \Omega_f \times (0,T] \mapsto \mathbb{R}, \hspace{1mm} \bm{\eta} \in \Omega_p \times (0,T] \mapsto \mathbb{R}^d, p_p \in \Omega_p \times (0,T] \mapsto \mathbb{R} $ s.t.} 
\begin{align}
\rho_f  \frac{\partial \bm{u}}{\partial t} - 2 \nu_f \nabla \cdot D(\bm{u}) + \nabla p_f &= \bm{f}_f \hspace{5mm} \text{ in }\Omega_f \times (0,T] \label{FPSI_FEM:StokesMom} \\
\nabla \cdot \bm{u} &= 0 \hspace{8mm} \text{ in }\Omega_f \times (0,T]  \label{FPSI_FEM:StokesMass} \\
\rho_p \frac{\partial^2 \bm{\eta}}{\partial t^2} - 2\nu_p \nabla \cdot D(\bm{\eta}) - \lambda \nabla (\nabla \cdot \bm{\eta}) + \alpha \nabla p_p &= \bm{f}_\eta \hspace{6mm} \text{ in }\Omega_p \times (0,T] \label{FPSI_FEM:Structure1} \\
s_0 \frac{\partial p_p}{\partial t} + \alpha \nabla \cdot \frac{\partial \eta}{\partial t} - \nabla \cdot \kappa \nabla p_p &= f_p \hspace{7mm} \text{ in } \Omega_p \times (0,T]. \label{FPSI_FEM:Structure2}
\end{align}
 Above, $D(\cdot)$ is the deformation rate tensor, defined as $D(\bm{v}) := \frac{1}{2}(\nabla \bm{v} + (\nabla \bm{v})^T)$. Densities are denoted by $\rho_f, \rho_p$, fluid viscosity by $\nu_f$, Lam\'e parameters by $\nu_p, \lambda$, and the Biot-Willis constant by $\alpha$. The constrained specific storage coefficient is denoted by $s_0$, and $\kappa$ represents the hydraulic conductivity. Although in general $\kappa$ is a tensor, we simplify here by considering an isotropic porous material so that $\kappa$ becomes a scalar. Each parameter is assumed to be constant in time for our analysis.  

Initial conditions are provided for $\bm{u}, \bm{\eta}$, and $p_p$. With stress tensors $\bm{\sigma}_f := 2\nu_f D(\bm{u}) - p_f I$ and $\bm{\sigma}_p := 2\nu_p D(\bm{\eta}) + \lambda (\nabla \cdot \bm{\eta})I - \alpha p_p I$, boundary data is given as
\begin{align}
\begin{split}\label{FPSI_FEM:noninterBCs}
\bm{\sigma}_f  \bm{n}_f &= \bm{u}_N \hspace{5mm} \text{ on } \Gamma_N^f  \times (0,T], \hspace{8mm}
\bm{u} = \bm{0} \hspace{5mm} \text{ on } \Gamma_D^f \times (0,T] \\
\bm{\sigma}_p  \bm{n}_p &= \bm{\eta}_N \hspace{5mm} \text{ on } \widetilde{\Gamma}_N^p \times (0,T], \hspace{8mm}
\bm{\eta} = \bm{0} \hspace{5mm} \text{ on } \widetilde{\Gamma}_D^p \times (0,T] \\
(\kappa \nabla p_p) \cdot \bm{n}_p &= p_N \hspace{5mm} \text{ on } \Gamma_N^{p} \times (0,T], \hspace{8mm} p_p = 0 \hspace{5mm} \text{ on } \Gamma_D^p \times (0,T],
\end{split}
\end{align}
where $\Gamma^f = \Gamma^f_N \cup \Gamma^f_D \cup \gamma$ is the Lipschitz continuous boundary of $\Omega_f$. Likewise, the boundary $\Gamma^p$ of $\Omega_p$ may be written as $\Gamma^p = \Gamma^p_N \cup \Gamma^p_D \cup \gamma$ = $\widetilde{\Gamma}^p_N \cup \widetilde{\Gamma}^p_D \cup \gamma$. We employ two notations in order to allow different types of boundary conditions to be defined for the displacement and pore pressure along the same spatially coincident portion of $\Gamma^p$. In each domain $\Omega_k$, $k \in \{f,p\}$, we assume the measure of $\Gamma^k_N$ and $\Gamma^k_D$ are nonzero; this assumption is necessary for the well-posedness of the formulation \cite{deCastro_2025_FPSIWP}. Take the unit vectors $\bm{n}_k$ to be outward normal to the domains, and let the unit vector $\bm{\tau}_\gamma$ be tangential to the interface $\gamma$.

\noindent We complete the system by providing the following interface conditions representing mass conservation, balance of stresses, and the Beavers-Joseph-Saffman (BJS) condition, where $\beta$ is the resistance parameter in the tangential direction:
\begin{align}
\bm{u} \cdot \bm{n}_f &= -\left(\frac{\partial \bm{\eta}}{\partial t} - \kappa \nabla p_p\right) \cdot \bm{n}_p \hspace{35mm} \text{ on } \gamma \times (0,T] \label{FPSI_FEM:Inter1} \\
\bm{\sigma}_f \bm{n}_f &= -\bm{\sigma}_p \bm{n}_p \hspace{57mm} \text{ on } \gamma \times (0,T] \label{FPSI_FEM:Inter2} \\
\bm{n}_f\cdot \bm{\sigma}_f \bm{n}_f &= -p_p \hspace{62mm} \text{ on } \gamma \times (0,T] \label{FPSI_FEM:Inter3} \\
\bm{n}_f\cdot \bm{\sigma}_f \bm{\tau}_\gamma^\ell &= -\beta \left(\bm{u} - \frac{\partial \bm{\eta}}{\partial t}\right) \cdot \bm{\tau}_\gamma^\ell \hspace{5mm} \text{ for } 1 \leq \ell \leq d-1, \hspace{5mm} \text{ on } \gamma \times (0,T]. \label{FPSI_FEM:Inter4}
\end{align}

\subsection{Derivation of Weak Form}
Define the following continuous spaces:
\begin{align}\label{FPSI_FEM:ContSpaces}
\begin{split}
&U :=\{\bm{v} \in \bm{H}^1(\Omega_f): \bm{v} = \bm{0} \text{ on } \Gamma_D^f\}, \hspace{5mm}
X := \{\bm{\varphi} \in \bm{H}^1(\Omega_p): \bm{\varphi} = \bm{0} \text{ on } \widetilde{\Gamma}_D^p\} \\
&Q_f :=L^2(\Omega_f), \hspace{38mm} Q_p := \{q_p \in H^1(\Omega_p) : q_p = 0 \text{ on } \Gamma_D^p\}.
\end{split}
\end{align}

Boldface font is used to distinguish a vector-valued function, such as $\bm{u},\bm{\eta}$, from a scalar-valued function such as $p_p, p_f$. Likewise, function spaces are typeset in bold to indicate their correspondence to a vector-valued function.
Let $H^s(\Omega_i)$ be the Hilbert space of order $s$ defined on subdomain $\Omega_i$, $i \in \{ f,p\}$. The notations $(\cdot,\cdot) = (\cdot,\cdot)_{\Omega_i}$ and $(\cdot,\cdot)_{1,\Omega_i}$ represent the $L^2$ and $H^1$ inner products on $\Omega_i$, respectively. The subscript $\Omega_i$ may be dropped from the inner product or norm notation if it is clear from context. A duality product between $H^s$ and its dual space for $s>0$ is denoted by $\langle \cdot,\cdot \rangle_{\Omega_i}$.  We define the $\bm{H}^1$ norm for vector-valued functions $\bm{v} \in \bm{H}^1(\Omega_i)$ as $||\bm{v}||_{1,\Omega_i}^2 := ||\bm{v}||_{0,\Omega_i}^2 + ||D(\bm{v})||_{0,\Omega_i}^2$, with  $||w||_{1,\Omega_i}^2 := ||w||_{0,\Omega_i}^2 + ||\nabla w||_{0,\Omega_i}^2$ the corresponding norm for scalar-valued functions $w \in H^1(\Omega_i)$. Likewise along a portion of the boundary $\Gamma$, we take $(\cdot,\cdot)_\Gamma$ to be the $L^2$ inner product, and $\langle \cdot,\cdot \rangle_{\Gamma}$ to represent a dual product.

 For simplicity, we derive the weak formulation in $\mathbb{R}^2$. Three Lagrange multipliers (LMs) are introduced to handle the interface conditions in $\mathbb{R}^2$: $g_1 \in \Lambda_{g1} := H^{-1/2}(\gamma)$, $g_2 \in \Lambda_{g2} := L^2(\gamma)$, and $\lambda_p \in \Lambda_\lambda := H^{1/2}(\gamma)$, defined by
\begin{equation*}
g_1 := (\bm{\sigma}_f \bm{n}_f)\cdot \bm{n}_f \text{ on } \gamma \times (0,T], \quad g_2 := (\bm{\sigma}_f \bm{n}_f)\cdot \bm{\tau}_\gamma \text{ on } \gamma \times (0,T], \quad \text{ and }  \lambda_p := \kappa \nabla p_p \cdot \bm{n}_p   \text{ on } \gamma \times (0,T].
\end{equation*}
Now, multiplying by appropriate test functions and integrating by parts, we derive the weak form of \eqref{FPSI_FEM:StokesMom}-\eqref{FPSI_FEM:Structure2}, for given $\bm{f}_f \in (\bm{H}^1(\Omega_f))^*$, $\bm{f}_\eta \in (\bm{H}^1(\Omega_p))^*$, and $f_p \in (H^1(\Omega_p))^*$:

\noindent \textit{Find $\bm{u} \in U, p_f \in Q_f, \bm{\eta} \in X, p_p \in Q_p, g_1 \in \Lambda_{g1}, g_2 \in \Lambda_{g2},$ and $ \lambda_p \in \Lambda_\lambda$ s.t. for a.e. $t \in (0,T]$, }
\begin{align}
\begin{split} \label{FPSI_FEM:WF:subdomains}
&\rho_f \Big( \frac{\partial \bm{u}}{\partial t}, \bm{v} \Big)_{\Omega_f} + 2 \nu_f \left( D(\bm{u}), D(\bm{v}) \right)_{\Omega_f} - (p_f, \nabla \cdot \bm{v})_{\Omega_f} - \langle g_1 \bm{n}_f,\bm{v}\rangle_\gamma  - ( g_2 \bm{\tau}_\gamma,\bm{v})_\gamma \\
&\hspace{15mm}= \langle\bm{f}_f,\bm{v}\rangle_{\Omega_f} + \langle\bm{u}_N, \bm{v}\rangle_{\Gamma_N^f} \hspace{3mm} \forall \ \bm{v} \in U, \\
&( \nabla \cdot \bm{u}, q )_{\Omega_f} = 0 \hspace{3mm} \forall \ q \in Q_f, \\
&\rho_p \Big( \frac{\partial^2 \bm{\eta}}{\partial t^2}, \bm{\varphi} \Big)_{\Omega_p} +  2 \nu_p \left( D(\bm{\eta}), D(\bm{\varphi}) \right)_{\Omega_p} + \lambda ( \nabla \cdot \bm{\eta}, \nabla \cdot \bm{\varphi} )_{\Omega_p} - \alpha(p_p, \nabla \cdot \bm{\varphi})_{\Omega_p}- \langle g_1 \bm{n}_p, \bm{\varphi}\rangle_\gamma + ( g_2 \bm{\tau}_\gamma, \bm{\varphi})_\gamma \\
&\hspace{15mm}= \langle\bm{f}_\eta, \bm{\varphi}\rangle_{\Omega_p} + \langle\bm{\eta}_N, \bm{\varphi}\rangle_{\widetilde{\Gamma}_N^p} \hspace{3mm} \forall \  \bm{\varphi} \in X, \\
&s_0 \Big( \frac{\partial p_p}{\partial t}, w\Big)_{\Omega_p} + \alpha\left( \nabla \cdot \frac{\partial \bm{\eta}}{\partial t},w \right)_{\Omega_p} + \kappa (\nabla p_p, \nabla w)_{\Omega_p} - (\lambda_p, w)_{\gamma}= \langle f_p,w \rangle_{\Omega_p} + \langle p_N, w \rangle_{\Gamma^p_N} \hspace{3mm} \forall \ w \in Q_p.
\end{split}
\end{align}

\noindent The boundary integrals involving $g_1 \bm{n}_f$ and $ g_2 \bm{\tau}_\gamma$  in \eqref{FPSI_FEM:WF:subdomains} derive from \eqref{FPSI_FEM:Inter2}, which implies that $g_1\bm{n}_f + g_2 \bm{\tau}_\gamma = \bm{\sigma}_f \bm{n}_f = -\bm{\sigma}_p \bm{n}_p$. The three LMs allow us to rewrite the remaining interface conditions \eqref{FPSI_FEM:Inter1}, \eqref{FPSI_FEM:Inter3}, and \eqref{FPSI_FEM:Inter4} as  
\begin{align}\label{FPSI_FEM:interWith3LMS}
    \begin{split}
        \bm{u}\cdot \bm{n}_f + \frac{\partial \bm{\eta}}{\partial t} \cdot \bm{n}_p - \lambda_p = 0, \quad
        g_1 + p_p = 0, \quad
        g_2 + \beta \bm{u} \cdot \bm{\tau}_\gamma - \beta \frac{\partial \bm{\eta}}{\partial t} \cdot \bm{\tau}_\gamma= 0.
    \end{split}
\end{align}
\begin{remark}
The above conditions represent a restriction to the case $d=2$; however, the extension to $d=3$ would only require the definition of one more LM for the additional tangential direction (i.e., $g_2 = (\bm{\sigma}_f \bm{n}_f) \cdot \bm{\tau}_\gamma^1$ and $g_3 = (\bm{\sigma}_f \bm{n}_f) \cdot \bm{\tau}_\gamma^2$). In the analysis, this new LM could be grouped with $g_2$ without affecting the structure of the proofs; we continue with $d=2$ for simplicity. 
\end{remark}
\noindent We derive the weak form of \eqref{FPSI_FEM:interWith3LMS} by multiplying with test functions $s_1 \in \Lambda_{g1}, \mu \in \Lambda_\lambda$, and $s_2 \in \Lambda_{g2}$, respectively, and integrating.
%\begin{align}\label{FPSI_FEM:interWith3LMS-WF}
%    \begin{split}
%        \langle \bm{u}\cdot \bm{n}_f,s_1\rangle_\gamma + \left\langle \frac{\partial \bm{\eta}}{\partial t} \cdot \bm{n}_p ,s_1 \right\rangle_\gamma - \langle \lambda_p,s_1\rangle_\gamma &= 0 \hspace{5mm} \forall \hspace{1mm} s_1 \in \Lambda_{g1},\\
%        \langle g_1,\mu\rangle_\gamma + (p_p,\mu)_\gamma &= 0 \hspace{5mm} \forall \hspace{1mm} \mu \in \Lambda_\lambda ,\\
%        (g_2,s_2)_\gamma + \beta (\bm{u} \cdot \bm{\tau}_\gamma,s_2)_\gamma - \beta \left(\frac{\partial \bm{\eta}}{\partial t} \cdot \bm{\tau}_\gamma,s_2\right)_\gamma &= 0 \hspace{5mm} \forall \hspace{1mm} s_2 \in \Lambda_{g2}.
%    \end{split}
%\end{align}
\noindent The semi-discrete monolithic system is obtained through time discretization of \eqref{FPSI_FEM:WF:subdomains} and the integrated interface equations. We implement Backward Euler, adopting the following notation to signify time derivatives compactly: 
\begin{equation}\label{FPSI_FEM:eq:TimeDerivNotation}
    \bm{\dot{\eta}}^{n} := \frac{\bm{\eta}^{n} - \bm{\eta}^{n-1}}{\Delta t}.
\end{equation}
Thus, the second derivative in time may be written as $ \bm{\ddot{\eta}}^{\hspace{0.1mm} n} :=  \dfrac{\bm{\dot{\eta}}^n - \bm{\dot{\eta}}^{n-1}}{\Delta t} = \dfrac{\bm{\eta}^n - 2\bm{\eta}^{n-1} + \bm{\eta}^{n-2}}{\Delta t^2} $. The time-discretized weak form for the FPSI system becomes:

\noindent \textit{Find $\bm{u}^{n+1} \in U, p_f^{n+1} \in Q_f, \bm{\eta}^{n+1} \in X, p_p^{n+1} \in Q_p, g_1^{n+1} \in \Lambda_{g1}, g_2^{n+1} \in \Lambda_{g2},$ and $ \lambda_p^{n+1} \in \Lambda_\lambda$ s.t. }

\begin{align}\label{FPSI_FEM:WF:allTimeDisc_StabError}
    \begin{split}
       &\rho_f ( \bm{\dot{u}}^{n+1}, \bm{v} )_{\Omega_f} + 2  \nu_f \left( D(\bm{u}^{n+1}), D(\bm{v}) \right)_{\Omega_f} -  (p_f^{n+1}, \nabla \cdot \bm{v})_{\Omega_f} -  \langle g_1^{n+1} \bm{n}_f,\bm{v}\rangle_\gamma  -  ( g_2^{n+1} \bm{\tau}_\gamma,\bm{v})_\gamma\\
      &\hspace{15mm} =  \langle\bm{f}_f^{n+1},\bm{v}\rangle_{\Omega_f} +  \langle\bm{u}_N^{n+1}, \bm{v}\rangle_{\Gamma_N^f}  \hspace{3mm} \forall \ \bm{v} \in U \\
&( \nabla \cdot \bm{u}^{n+1}, q )_{\Omega_f} = 0 \hspace{3mm} \forall \ q \in Q_f \\
&\frac{\rho_p}{\Delta t} \left( \bm{\dot{\eta}}^{n+1} - \bm{\dot{\eta}}^n, \bm{\varphi} \right)_{\Omega_p} +  2  \nu_p \left( D\left(\bm{\eta}^{n+1}\right), D(\bm{\varphi}) \right)_{\Omega_p} + \lambda  \left( \nabla \cdot \bm{\eta}^{n+1}, \nabla \cdot \bm{\varphi} \right)_{\Omega_p} - \alpha (p_p^{n+1}, \nabla \cdot \bm{\varphi})_{\Omega_p}\\
&\hspace{15mm} - \langle g_1^{n+1}\bm{n}_p, \bm{\varphi}\rangle_\gamma +  ( g_2^{n+1}\bm{\tau}_\gamma, \bm{\varphi})_\gamma  =  \langle\bm{f}_\eta^{n+1}, \bm{\varphi}\rangle_{\Omega_p} +  \langle\bm{\eta}_N^{n+1}, \bm{\varphi}\rangle_{\widetilde{\Gamma}_N^p} \hspace{3mm} \forall \ \bm{\varphi} \in X \\
&s_0 ( \dot{p}_p^{n+1}, w)_{\Omega_p} + \alpha\left(\nabla \cdot \bm{\dot{\eta}}^{n+1},w\right)_{\Omega_p} + \kappa  (\nabla p_p^{n+1}, \nabla w)_{\Omega_p} -  ( \lambda_p^{n+1}, w)_{\gamma} \\
&\hspace{15mm}=  \langle f_p^{n+1},w \rangle_{\Omega_p} +  \langle p_N^{n+1}, w \rangle_{\Gamma_N^p}  \hspace{3mm} \forall \ w \in Q_p \\
  & \langle \bm{u}^{n+1} \cdot \bm{n}_f,s_1\rangle_\gamma + \left\langle \bm{\dot{\eta}}^{n+1} \cdot \bm{n}_p ,s_1 \right\rangle_\gamma   -  \langle \lambda_p^{n+1},s_1\rangle_\gamma = 0 \hspace{5mm} \forall \ s_1 \in \Lambda_{g1}\\
      & \langle g_1^{n+1},\mu\rangle_\gamma   + (  p_p^{n+1},\mu)_\gamma 
            + \overline{\epsilon}(\lambda_{p}^{n+1},\mu)_{1/2,\gamma}  = 0 \hspace{5mm} \forall \ \mu \in \Lambda_\lambda\\
      &  \frac{1}{\beta } ( g_2^{n+1},s_2)_\gamma  +  (\bm{u}^{n+1} \cdot \bm{\tau}_\gamma,s_2)_\gamma - \left(\bm{\dot{\eta}}^{n+1} \cdot \bm{\tau}_\gamma,s_2\right)_\gamma = 0 \hspace{5mm} \forall \ s_2 \in \Lambda_{g2}.
    \end{split}
\end{align}
\normalsize 

The system \eqref{FPSI_FEM:WF:allTimeDisc_StabError} has been shown to be well-posed through exploring its saddle point structure in our previous paper \cite{deCastro_2025_FPSIWP}, wherein the stability of the fully discrete system was also proved. 
Note that the system includes the 
 stabilization term $ \overline{\epsilon}(\lambda_{p}^{n+1},\mu)_{1/2,\gamma}$, which has been added to show the well-posedness of the semi-discrete form.
For the finite element discretization, assume that $\Omega_f, \Omega_p$ are convex polytopal domains. Let $h_1, h_2, \text{and } h_\gamma$ represent the mesh sizes of a quasi-uniform partition of $\Omega_f, \Omega_p, \text{and } \gamma$. The conforming discrete finite element spaces are denoted by adding the superscript $h$; i.e. $U^h \subset U$. For ease of notation, we do not distinguish between $h_1, h_2, h_\gamma$ in the notation as it is clearly inherited from the domain.

Numerical results suggest that the stabilization term $ \overline{\epsilon}(\lambda_{p}^{n+1},\mu)_{1/2,\gamma}$ is not needed in practice. All results shown in Section \ref{FPSI_FEM:sec:Numerical} are calculated for $\overline{\epsilon} = 0$. However, showing well-posedness of the continuous formulation with $\overline{\epsilon} = 0$ is an open question still. 
	\section{Spatial convergence analysis}\label{FPSI_FEM:sec:errorDiscrete}

%%%%%%%%%%%%%%%%%%%%%%%%%%%%%%%%%%%5
We first study the spatial convergence properties of our formulation, examining the error between the fully discrete solution as an approximation to the semi-discrete (in time) solution. Since the main goal of this work is to develop a domain decomposition scheme for the time discretized monolithic formulation, in the following analysis, we focus only on spatial discretization error in terms of mesh size $h$, and not on $\Delta t$.

Let $k_u$ be the degree of polynomial approximation for the space $U^h$, with polynomial degrees for other variables denoted similarly. Let $C$ denote a generic constant independent of mesh size $h$. For the fluid velocity, consider the Stokes-like projection operator $\mathcal{I}^{U^h} : U \rightarrow U^h$ given by the solution to the problem \cite{Ambartsumyan_2018, Li_2022_Hydro}:
\begin{align}\label{FPSI_FEM:eq:StokesProj_Eqs}
\begin{split}
    2\nu_f (D(\mathcal{I}^{U^h}(\bm{u})),D(\bm{v}_h))_{\Omega_f} - (\nabla \cdot \bm{v}_h, \widetilde{p}_{f,h})_{\Omega_f} &= 2\nu_f (D(\bm{u}),D(\bm{v}_h))_{\Omega_f} \ \forall \ \bm{v}_h \in U^h \\
    (\nabla \cdot \mathcal{I}^{U^h}(\bm{u}), q_h)_{\Omega_f} &= (\nabla \cdot \bm{u}, q_h)_{\Omega_f} \ \forall \ q_h \in Q_{f}^h.
    \end{split}
\end{align}
This operator satisfies the approximation property 
\begin{equation}\label{FPSI_FEM:eq:Approx_StokesProjection}
     || \bm{u} - \mathcal{I}^{U^h}(\bm{u})||_1 \leq C h_1^{k_u} ||\bm{u}||_{k_u+1}.
\end{equation}

\noindent For the displacement and pressures, we utilize $L^2$ projection operators \cite{Ciarlet_1978, Ern_2004}. Let $\mathcal{P}^{X^h}: X \rightarrow X^h$ satisfy  
      $(\mathcal{P}^{X^h}(\bm{\eta}),\bm{\varphi}_h) = (\bm{\eta},\bm{\varphi}_h) \quad \forall \  \bm{\varphi}_h \in X^h,$ 
with approximation properties
\begin{align}\label{FPSI_FEM:eq:ApproxProps_EtaPressures}
\begin{split}
        &|| \bm{\eta} - \mathcal{P}^{X^h}(\bm{\eta}) ||_0 \leq C h_2^{k_\eta+1} ||\bm{\eta}||_{k_\eta+1}, \quad \text{ and } \quad 
    || \bm{\eta} - \mathcal{P}^{X^h}(\bm{\eta})||_1 \leq C h_2^{k_\eta} ||\bm{\eta}||_{k_\eta+1}.\\
  %  &|| p_p - \mathcal{P}^{Q^{h}_p}(p_p) ||_0 \leq C h_2^{k_{pp}+1} ||p_p||_{k_{pp}+1}, \qquad  
 %   || p_p - \mathcal{P}^{Q^{h}_p}(p_p)||_1 \leq C h_2^{k_{pp}} ||p_p||_{k_{pp}+1} \\
  %    &|| p_f - \mathcal{P}^{Q_f^h}(p_f)||_0 \leq C h_1^{k_{pf}+1} ||p_f||_{k_{pf}+1}.
    \end{split}
\end{align}
Projectors $\mathcal{P}^{Q_f^h}, \mathcal{P}^{Q_p^h}$ are defined similarly. % Note, Layton uses the L2 projection, and Layton and Ern/Guermond give the best approximation property. Ciarlet gives all the stuff here. Can't find a reference stating that the L^2 projection actually is the element that achieves the infimum, so leaving it out.  
%% I think that reference is in the book by Brezzi, Fortin 2013 (Mixed Finite Elements) pg 131 Section 3.1
%  || \bm{\eta} - \mathcal{P}^{X^h}(\bm{\eta})||_0 = \underset{\bm{\varphi}_h \in X^h}{\text{inf}} || \bm{\eta} - \bm{\varphi}_h||_0,
   %\quad (\mathcal{P}^{Q_p^h}(p_p),w_h) = (p_p,w_h) \quad \forall \ w_h \in Q_p^h \\
     % & (\mathcal{P}^{Q_f^h}(p_f),q_h) = (p_f,q_h) \quad \forall \ q_h \in Q_f^h, 
For the Lagrange multiplier $g_2$, we use the following estimates, where $\mathcal{P}^{\Lambda_{g_2}^h}$ is the $L^2$ projection  \cite{Ciarlet_1978, Li_2022_Hydro}  satisfying
%note to self: L2 projection on boundary gives optimal error rate, interpolant gives suboptimal (Gunzburger 1992, remark pg 415)
\begin{align}\label{FPSI_FEM:eq:ApproxProps_G2}
&(\mathcal{P}^{\Lambda_{g_2}^h}(g_2),s_{h})_{0,\gamma} = (g_2,s_{h})_{0,\gamma} \ \forall \ s_{h} \in \Lambda_{g_2}^h, \quad ||g_2 - \mathcal{P}^{\Lambda_{g_2}^h}(g_2)||_{0,\gamma} \leq C h_\gamma^{k_{g2}+1}||g_2||_{k_{g2}+1,\gamma}. 
\end{align}
Assume $\Omega$ is a convex polyhedral domain. Approximation errors for $g_1$ and $\lambda_p$ use the following property of the subspaces $\Lambda_{g1}^h, \Lambda_{\lambda}^h$.
For $k^h \in \Lambda^h(\gamma)\subset H^{1/2}(\gamma)$, the following inverse inequality is satisfied \cite{Ciarlet_1978, Gunzburger_1992}:
$$||k^h||_{s,\gamma} \leq C h_\gamma^{t-s} ||k^h||_{t,\gamma} \quad \forall \ k^h \in \Lambda^h(\gamma),  \ -1/2 \leq t \leq s \leq 1/2.$$ 
%Let $\Lambda_{g1}^h \subset H^{1/2}(\gamma)$  {\color{red} do we need higher $(H^{1/2})$ regulairty for $g_1^h$ for error estimation? -- it's needed so that we can apply the properties from (2.14)-(2.17) in Max's paper} and $\Lambda_p^h \subset H^{1/2}(\gamma)$ where both discrete spaces satisfy the inverse inequality  {\color{red} \cite{Gunzburger_1992}} $$||k^h||_{s,\gamma} \leq C h_\gamma^{t-s} ||k^h||_{t,\gamma} \quad \forall \ k^h \in \Lambda^h(\gamma),  \ -1/2 \leq t \leq s \leq 1/2.$$ 
\begin{comment} % more general statement
Then $ \exists \ k \in \mathbb{Z}$ for which the following approximation properties hold for $r \in \{f,p\}$ \cite{Gunzburger_1992}: 
  
\begin{align}\label{FPSI_FEM:eq:ApproxProps_G1}
& \underset{g_{1,h} \in \Lambda_{g_1}^h}{\inf} ||g_1 - g_{1,h}||_{-1/2,\gamma} \leq Ch_\gamma^{m} \underset{\substack{w \in H^m(\Omega_r) \\ w|_\gamma = g_1}}{\inf}||w||_{m,\Omega_r} \ \forall \ g_1 \in H^{m}(\Omega_r) |_\gamma, 1 \leq m \leq k,
\end{align}

\begin{align}\label{FPSI_FEM:eq:H1/2approxProp}
 \underset{\lambda_{p,h} \in \Lambda_{p}^h}{\inf} ||\lambda_p - \lambda_{p,h}||_{1/2,\gamma} \leq Ch_\gamma^{m-1} \underset{\substack{w \in H^m(\Omega_r) \\ w|_\gamma = \lambda_p}}{\inf}||w||_{m,\Omega_r} \ \forall \ \lambda_p \in H^{m}(\Omega_r) |_\gamma, 1 \leq m \leq k.
\end{align}
\end{comment}
\noindent Consider $g_1, \lambda_p$ as the trace of $H^{k_{g1}+1}(\Omega_r),  H^{k_\lambda+1}(\Omega_r)$ functions, respectively, for $r \in \{f,p\}$:
\begin{align}
    &\underset{g_{1,h} \in \Lambda_{g_1}^h}{\inf} ||g_1 - g_{1,h}||_{-1/2,\gamma} \leq Ch_\gamma^{k_{g1}+1} \underset{\substack{w \in H^{k_{g1}+1}(\Omega_r) \\ w|_\gamma = g_1}}{\inf}||w||_{k_{g1}+1,\Omega_r} \ \forall \ g_1 \in H^{k_{g1}+1}(\Omega_r) |_\gamma, \label{FPSI_FEM:eq:ApproxProps_G1}\\  &\underset{\lambda_{p,h} \in \Lambda_{p}^h}{\inf} ||\lambda_p - \lambda_{p,h}||_{1/2,\gamma} \leq Ch_\gamma^{k_\lambda} \underset{\substack{w \in H^{k_\lambda+1}(\Omega_r) \\ w|_\gamma = \lambda_p}}{\inf}||w||_{k_\lambda+1,\Omega_r} \ \forall \ \lambda_p \in H^{k_\lambda+1}(\Omega_r) |_\gamma \label{FPSI_FEM:eq:H1/2approxProp}.
\end{align}

\noindent Lastly, we define an energy norm for the displacement
$$ ||\bm{\eta}^{n+1}||_E^2 := 2\nu_p ||D(\bm{\eta}^{n+1})||^2_{0,\Omega_p} + \lambda ||\nabla \cdot \bm{\eta}^{n+1}||_{0,\Omega_p}^2,$$
and list the following inequalities for all $\bm{u} \in U, \bm{\eta} \in X$, and all $p_p \in Q_p$ (see \cite{Bukac_2015OpSplit}, for example). Each constant only depends on the domain. 
\begin{align}\label{FPSI_FEM:eq:TraceIneqs} 
\begin{split}
&\text{Trace Inequalities: } \quad  ||\bm{u}||_{0,\gamma} \leq C_T ||\bm{u}||_{1,\Omega_f}, \quad ||\bm{\eta}||_{0,\gamma} \leq C_T ||\bm{\eta}||_{1,\Omega_p}, \quad ||p_p||_{0,\gamma} \leq C_T ||p_p||_{1,\Omega_p} \\
&\text{Korn's Inequalities: } \quad ||\nabla \bm{u}||_{0,\Omega_f} \leq C_K ||D(\bm{u})||_{0,\Omega_f}, \qquad ||\nabla \bm{\eta}||_{0,\Omega_p} \leq C_K ||D(\bm{\eta})||_{0,\Omega_p}. 
\end{split}
\end{align}
As functions in $U, X, Q_p$ satisfy homogeneous Dirichlet conditions on a component of the boundary, the Poincar\'{e} inequalities hold:
\begin{align}\label{FPSI_FEM:eq:PoincareIneqs}
\begin{split}
||\bm{u}||_{0,\Omega_f} \leq C_P ||\nabla \bm{u}||_{0,\Omega_f}, \qquad ||\bm{\eta}||_{0,\Omega_p} \leq C_P ||\nabla \bm{\eta}||_{0,\Omega_p}, \quad ||p_p||_{0,\Omega_p} \leq C_P ||\nabla p_p||_{0,\Omega_p}.
\end{split}
\end{align}
Combining Korn's and Poincar\'{e} inequalities, we see that there exists a constant $C_{KP}$ satisfying
\begin{equation}\label{FPSI_FEM:eq:1normbound}
||\bm{u}||_{1}^2 \leq C_{KP} ||D(\bm{u})||^2_{0}, \qquad ||\bm{\eta}||_{1}^2 \leq  C_{KP} ||D(\bm{\eta})||_0^2, \qquad ||p_p||_1^2 \leq C_{KP} ||\nabla p_p||_0^2.
\end{equation}
% technically $C_{KP} := C_P^2 C_K^2 + 1$ for u,\eta and $C_{P2} := (C_P^2 + 1)$ for p_p.

We restate \eqref{FPSI_FEM:WF:allTimeDisc_StabError} with test functions in the discrete subspaces, assuming homogeneous Neumann conditions for simplicity. Define $A = U \times Q_f \times X \times Q_p \times \Lambda_{g1} \times \Lambda_{g2} \times \Lambda_\lambda$ and its discrete equivalent $A^h$.  Then for all $\bm{a}_h := (\bm{v}_h,q_h, \bm{\varphi}_h,w_h,s_{1,h}, s_{2,h},\mu_{h}) \in A^h$, we look for $(\bm{u}_h^{n+1}, p_{f,h}^{n+1}, \bm{\eta}_h^{n+1}, p_{p,h}^{n+1}, g_{1,h}^{n+1}, g_{2,h}^{n+1}, \lambda_{p,h}^{n+1}) \in A^h$ such that 
\begin{align}\label{FPSI_FEM:WF:RestateForError_Discrete}
    \begin{split}
       &\rho_f ( \bm{\dot{u}}_h^{n+1}, \bm{v}_h )_{\Omega_f} + 2  \nu_f \left( D(\bm{u}_h^{n+1}), D(\bm{v}_h) \right)_{\Omega_f} -  (p_{f,h}^{n+1}, \nabla \cdot \bm{v}_h)_{\Omega_f} -  \langle g_{1,h}^{n+1} \bm{n}_f,\bm{v}_h\rangle_\gamma-  ( g_{2,h}^{n+1} \bm{\tau}_\gamma,\bm{v}_h)_\gamma\\
      &\hspace{10mm} =  \langle\bm{f}_f^{n+1},\bm{v}_h\rangle_{\Omega_f}  \\
&\frac{\rho_p}{\Delta t} \left( \bm{\dot{\eta}}_h^{n+1} - \bm{\dot{\eta}}_h^n, \bm{\varphi}_h \right)_{\Omega_p} +  2  \nu_p \left( D\left(\bm{\eta}_h^{n+1}\right), D(\bm{\varphi}_h) \right)_{\Omega_p} + \lambda  \left( \nabla \cdot \bm{\eta}_h^{n+1}, \nabla \cdot \bm{\varphi}_h \right)_{\Omega_p} \\
& \hspace{10mm}- \alpha (p_{p,h}^{n+1}, \nabla \cdot \bm{\varphi}_h)_{\Omega_p}- \langle g_{1,h}^{n+1}\bm{n}_p, \bm{\varphi}_h\rangle_\gamma +  ( g_{2,h}^{n+1}\bm{\tau}_\gamma, \bm{\varphi}_h)_\gamma =  \langle\bm{f}_\eta^{n+1}, \bm{\varphi}_h\rangle_{\Omega_p}  \\
&s_0 ( \dot{p}_{p,h}^{n+1}, w_h)_{\Omega_p} + \alpha\left(\nabla \cdot \bm{\dot{\eta}}_h^{n+1},w_h\right)_{\Omega_p} + \kappa  (\nabla p_{p,h}^{n+1}, \nabla w_h)_{\Omega_p} -  ( \lambda_{p,h}^{n+1}, w_h)_{\gamma} = \langle f_p^{n+1},w_h \rangle_{\Omega_p}  \\
   &  \frac{1}{\beta } ( g_{2,h}^{n+1},s_{2,h})_\gamma  +  (\bm{u}_h^{n+1} \cdot \bm{\tau}_\gamma,s_{2,h})_\gamma - \left(\bm{\dot{\eta}}_h^{n+1} \cdot \bm{\tau}_\gamma,s_{2,h}\right)_\gamma = 0\\
      &  \langle g_{1,h}^{n+1},\mu_h\rangle_\gamma  + (  p_{p,h}^{n+1},\mu_h)_\gamma + \overline{\epsilon}(\lambda_{p,h}^{n+1},\mu_h)_{1/2,\gamma} = 0 \\
      &( \nabla \cdot \bm{u}_h^{n+1}, q_h )_{\Omega_f} = 0 \\
        & \langle \bm{u}_h^{n+1} \cdot \bm{n}_f,s_{1,h}\rangle_\gamma + \left\langle \bm{\dot{\eta}}_h^{n+1} \cdot \bm{n}_p ,s_{1,h} \right\rangle_\gamma   -  \langle \lambda_{p,h}^{n+1},s_{1,h}\rangle_\gamma = 0. 
    \end{split}
\end{align}

%The constraint equations imply that the solution $(\bm{u}^{n+1},\bm{\dot{\eta}}^{n+1},p_p^{n+1},g_{2}^{n+1},\lambda_p^{n+1})$ is in the space 
%\begin{align*}
%    K &:= \{ (\bm{v}, \bm{\varphi}, w , s_{2},\lambda_{p}) \in M : \ b_{MZ}\left( (\bm{v},\bm{\varphi},w,s_{2},\lambda_{p}),(q,s_{1}) \right)= 0 \quad \forall \ (q,s_{1}) \in Z \}.
%\end{align*}

%\noindent Let $(\bm{u}^{n+1}, p_{f}^{n+1}, \bm{\eta}^{n+1}, p_{p}^{n+1}, g_{1}^{n+1}, g_{2}^{n+1}, \lambda_{p}^{n+1}) \in A$ be the semidiscrete solutions to \eqref{FPSI_FEM:WF:RestateForError_Discrete}.

\begin{theorem}\label{FPSI_FEM:thm:generalError}
 Let $(\bm{u}^{n+1}, p_{f}^{n+1}, \bm{\eta}^{n+1}, p_{p}^{n+1}, g_{1}^{n+1}, g_{2}^{n+1}, \lambda_{p}^{n+1}) \in A$ be the solution to the semi-discrete system \eqref{FPSI_FEM:WF:allTimeDisc_StabError} with $\bm{u}_N^{n+1}=\bm{\eta}_N^{n+1} =\bm{0}$.
Assume that this solution has sufficient regularity. Then with final time $T := N \Delta t$, the spatial error between the fully discrete solution and the semi-discrete solution at time step $M$ satisfies the following, for $0 < M \leq N$:
\begin{small}
\begin{align*}
 & ||\bm{\eta}_h^M - \bm{\eta}^M||_{1} +  ||\bm{u}_h^M - \bm{u}^M||_{0} +  ||\bm{\dot{\eta}}_h^M - \bm{\dot{\eta}}^M||_{0} +||p_{p,h}^M - p_p^M ||_{0} \\ 
    &+ \sqrt{\Delta t} \sum_{n=0}^{M-1} \Big[ ||\bm{u}_h^{n+1} - \bm{u}^{n+1}||_{1} + ||\bm{\dot{\eta}}_h^{n+1} - \bm{\dot{\eta}}^{n+1}||_{1} +  || p_{p,h}^{n+1} - p_{p}^{n+1}||_{1} +  || g_{2,h}^{n+1} - g_{2}^{n+1}||_{0,\gamma} \\
    &+  || \lambda_{p,h}^{n+1} - \lambda_{p}^{n+1}||_{1/2,\gamma} +  || p_{f,h}^{n+1} - p_{f}^{n+1}||_{0} +  || g_{1,h}^{n+1} - g_{1}^{n+1}||_{-1/2,\gamma} \Big] \\
      &\leq C_a \sqrt{2(7M+4)} \Big(   h_2^{k_\eta}||\bm{\eta}^M||_{k_\eta+1} +  h_1^{k_u}||\bm{u}^M||_{k_u+1}+   h_2^{k_\eta+1}  ||\bm{\dot{\eta}}^M||_{k_\eta+1}  + h_2^{k_{pp}+1}||p_p^M||_{k_{pp}+1} \\
     &+ ( \Delta t \ \overline{C})^{1/2} \sum_{n=0}^{M-1} \Big[h_1^{k_u}||\bm{\dot{u}}^{n+1}||_{k_u+1} + h_1^{k_u}||\bm{u}^{n+1}||_{k_u+1} +   h_2^{k_\eta} || \bm{\eta}^{n+1}||_{k_\eta+1}\\
     &+  h_2^{k_\eta}||\bm{\dot{\eta}}^{n+1}||_{k_\eta+1} + h_2^{k_{pp}+1}||p_p^{n+1}||_{k_{pp}+1} +  h_2^{k_{pp}} ||p_p^{n+1}||_{k_{pp}+1} +  h_\gamma^{k_{g2}+1} ||g_2^{n+1}||_{k_{g2}+1,\gamma} \\
     & + h_\gamma^{k_\lambda} \underset{\substack{w \in H^{k_\lambda+1}(\Omega_r) \\ w|_\gamma = \lambda_p }}{\inf} ||w||_{k_\lambda+1,\Omega_r} +  h_1^{k_{pf}+1} ||p_f^{n+1}||_{k_{pf}+1}  +  h_\gamma^{k_{g1}+1} \underset{\substack{w \in H^{k_{g1}+1}(\Omega_r) \\ w|_\gamma = g_1 }}{\inf} ||w||_{k_{g1}+1,\Omega_r} \Big]       \Big).
\end{align*}
\end{small}

\end{theorem}

\begin{proof} For clarity of presentation, we outline the major steps of the proof.\\
\textbf{Step 1: Define and separate approximation and truncation error terms}\\
Define $\bm{e}_u^{n+1} := \bm{u}_h^{n+1} - \bm{u}^{n+1}$ as the error between the fully discrete and semidiscrete solutions for the Stokes velocity, with the errors for other variables defined similarly. Then we have $\bm{\dot{u}}_h^{n+1} - \bm{\dot{u}}^{n+1} = \bm{\dot{e}}_u^{n+1}. $ Subtract the semidiscrete system \eqref{FPSI_FEM:WF:allTimeDisc_StabError}, with test functions in the discrete FE spaces, from the fully discrete weak form \eqref{FPSI_FEM:WF:RestateForError_Discrete}. Then for all $\bm{a}_h \in A^h$:
\begin{align}\label{FPSI_FEM:eq:FirstErrorEqnSubtracted}
    \begin{split}
       &\frac{\rho_f}{\Delta t} ( \bm{e}_u^{n+1} - \bm{e}_u^{n}, \bm{v}_h )_{\Omega_f} + 2  \nu_f \left( D(\bm{e}_u^{n+1}), D(\bm{v}_h) \right)_{\Omega_f} -  (e_{pf}^{n+1}, \nabla \cdot \bm{v}_h)_{\Omega_f} -  \langle e_{g1}^{n+1} \bm{n_f},\bm{v}_h\rangle_\gamma  \\
       & \hspace{10mm} -  ( e_{g2}^{n+1} \bm{\tau_\gamma},\bm{v}_h)_\gamma =0 \\
&\frac{\rho_p}{\Delta t} \left( \bm{\dot{e}}_\eta^{n+1} - \bm{\dot{e}}_\eta^n, \bm{\varphi}_h \right)_{\Omega_p} +  2  \nu_p \left( D\left(\bm{e}_\eta^{n+1}\right), D(\bm{\varphi}_h) \right)_{\Omega_p} + \lambda  \left( \nabla \cdot  \bm{e}_\eta^{n+1}, \nabla \cdot \bm{\varphi}_h \right)_{\Omega_p} \\
& \hspace{10mm}- \alpha (e_{pp}^{n+1}, \nabla \cdot \bm{\varphi}_h)_{\Omega_p}- \langle e_{g1}^{n+1}\bm{n_p}, \bm{\varphi}_h\rangle_\gamma +  ( e_{g2}^{n+1}\bm{\tau_\gamma}, \bm{\varphi}_h)_\gamma =  0  \\
&\frac{s_0}{\Delta t} ( e_{pp}^{n+1} - e_{pp}^n, w_h)_{\Omega_p} + \alpha\left(\nabla \cdot \bm{\dot{e}}_\eta^{n+1},w_h\right)_{\Omega_p} + \kappa  (\nabla e_{pp}^{n+1}, \nabla w_h)_{\Omega_p} -  ( e_{\lambda p}^{n+1}, w_h)_{\gamma} = 0  \\
   &  \frac{1}{\beta } ( e_{g2}^{n+1},s_{2,h})_\gamma  +  (\bm{e}_u^{n+1} \cdot \bm{\tau_\gamma},s_{2,h})_\gamma - \left(\bm{\dot{e}}_\eta^{n+1} \cdot \bm{\tau_\gamma},s_{2,h}\right)_\gamma = 0 \\
      &   \langle e_{g1}^{n+1},\mu_h\rangle_\gamma   + (  e_{pp}^{n+1},\mu_h)_\gamma + \overline{\epsilon}( e_{\lambda p}^{n+1},\mu_h)_{1/2,\gamma} = 0\\
       &( \nabla \cdot \bm{e}_u^{n+1}, q_h )_{\Omega_f} = 0 \\
        & \langle \bm{e}_u^{n+1} \cdot \bm{n_f},s_{1,h}\rangle_\gamma + \left\langle \bm{\dot{e}}_\eta^{n+1} \cdot \bm{n_p} ,s_{1,h} \right\rangle_\gamma   -  \langle e_{\lambda p}^{n+1},s_{1,h}\rangle_\gamma = 0.   
    \end{split}
\end{align}
Split the errors into a truncation error $\chi_{r,h}^{n+1}:=  r_h^{n+1} - \tilde{r}$ in the discrete subspace and an approximation error $\theta_{r}^{n+1} := r^{n+1} - \tilde{r}$ for a variable $r$, where $\tilde{r}$ represents the projection or interpolant of $r$ corresponding to the operators defined in \eqref{FPSI_FEM:eq:StokesProj_Eqs}- \eqref{FPSI_FEM:eq:H1/2approxProp}. For $r = g_{1}^{n+1}$ (respectively $\lambda_p^{n+1}$), $\tilde{r}$ is an arbitrary element $s_{1,h} \in \Lambda_{g1}^h$ (respectively, $\mu_h \in \Lambda_\lambda^h$).
\begin{comment}
Let $s_{1,h} \in \Lambda_{g1}^h$ and $\mu_h \in \Lambda_\lambda^h$ be arbitrary, and define, for \eqref{FPSI_FEM:eq:StokesProj_Eqs}
\begin{align*}
    \bm{e}_u^{n+1} &= (\bm{u}_h^{n+1} - \mathcal{I}^{U^h}(\bm{u}^{n+1}) ) - ( \bm{u}^{n+1} - \mathcal{I}^{U^h}(\bm{u}^{n+1}) ) := \bm{\chi}_{u,h}^{n+1} - \bm{\theta}_u^{n+1} \\
    \bm{e}_\eta^{n+1} &= (\bm{\eta}_h^{n+1} - \mathcal{P}^{X^h}(\bm{\eta}^{n+1} ) - ( \bm{\eta}^{n+1} - \mathcal{P}^{X^h}(\bm{\eta}^{n+1}) ) := \bm{\chi}_{\eta,h}^{n+1} - \bm{\theta}_\eta^{n+1} \\
      e_{pp}^{n+1} &= (p_{p,h}^{n+1} -  \mathcal{P}^{Q_p^h}(p_p^{n+1}) ) - ( p_p^{n+1} - \mathcal{P}^{Q_p^h}(p_p^{n+1}) ) := \chi_{pp,h}^{n+1} - \theta_{pp}^{n+1} \\
        e_{pf}^{n+1} &= (p_{f,h}^{n+1} -  \mathcal{P}^{Q_f^h}(p_f^{n+1}) ) - ( p_f^{n+1} - \mathcal{P}^{Q_f^h}(p_f^{n+1}) ) := \chi_{pf,h}^{n+1} - \theta_{pf}^{n+1} \\
        e_{g2}^{n+1} &= (g_{2,h}^{n+1} -  \mathcal{P}^{\Lambda_{g2}^h}(g_2^{n+1}) ) - ( g_2^{n+1} - \mathcal{P}^{\Lambda_{g2}^h}(g_{2}^{n+1}) ) := \chi_{g2,h}^{n+1} - \theta_{g2}^{n+1} \\
  %  e_{\lambda p}^{n+1} &= (\lambda_{p,h}^{n+1} -  \mathcal{P}^{\Lambda_\lambda^h}_{1/2}(\lambda_p^{n+1}) ) - ( \lambda_p^{n+1} - \mathcal{P}^{\Lambda_\lambda^h}_{1/2}(\lambda_p^{n+1}) ) := \chi_{\lambda p,h}^{n+1} - \theta_{\lambda p}^{n+1} \\
    e_{\lambda p}^{n+1} &= (\lambda_{p,h}^{n+1} - \mu_h) - (\lambda_p^{n+1} - \mu_h) := \chi_{\lambda p,h}^{n+1} - \theta_{\lambda p}^{n+1}\\
    e_{g1}^{n+1} &= (g_{1,h}^{n+1} - s_{1,h} ) - ( g_1^{n+1} -  s_{1,h}) := \chi_{g1,h}^{n+1} - \theta_{g1}^{n+1}.
\end{align*}
\end{comment}
\noindent Rewrite the system in terms of these two components, moving the approximation errors $\theta^{n+1}_r$ to the right hand side and leaving the truncation errors $\chi_{r,h}^{n+1}$ on the left. Thus for all $\bm{a}_h  \in A^h$,
\refstepcounter{equation}\label{FPSI_FEM:eq:splitChiTheta}
\begin{align}
  %  \begin{split}
       &\frac{\rho_f}{\Delta t} ( \bm{\chi}_{u,h}^{n+1} - \bm{\chi}_{u,h}^{n}, \bm{v}_h )_{\Omega_f} + 2  \nu_f \left( D(\bm{\chi}_{u,h}^{n+1}), D(\bm{v}_h) \right)_{\Omega_f} -  (\chi_{pf,h}^{n+1}, \nabla \cdot \bm{v}_h)_{\Omega_f} -  \langle \chi_{g1,h}^{n+1} \bm{n}_f,\bm{v}_h\rangle_\gamma\tag{\theequation (a)}\label{FPSI_FEM:eq:splitChiTheta_a}  \\
      &\hspace{10mm} -  ( \chi_{g2,h}^{n+1} \bm{\tau}_\gamma,\bm{v}_h)_\gamma  =  \frac{\rho_f}{\Delta t} ( \bm{\theta}_{u}^{n+1} - \bm{\theta}_{u}^{n}, \bm{v}_h )_{\Omega_f}  + 2  \nu_f ( D(\bm{\theta}_{u}^{n+1}), D(\bm{v}_h) )_{\Omega_f}  \notag \\
      &\hspace{10mm} -  (\theta_{pf}^{n+1}, \nabla \cdot \bm{v}_h)_{\Omega_f}-  \langle \theta_{g1}^{n+1} \bm{n}_f,\bm{v}_h\rangle_\gamma -  ( \theta_{g2}^{n+1} \bm{\tau}_\gamma,\bm{v}_h)_\gamma \notag \\
&\frac{\rho_p}{\Delta t} \left( \bm{\dot{\chi}}_{\eta,h}^{n+1} - \bm{\dot{\chi}}_{\eta,h}^n, \bm{\varphi}_h \right)_{\Omega_p} +  2  \nu_p ( D(\bm{\chi}_{\eta,h}^{n+1}), D(\bm{\varphi}_h) )_{\Omega_p} + \lambda  ( \nabla \cdot  \bm{\chi}_{\eta,h}^{n+1}, \nabla \cdot \bm{\varphi}_h )_{\Omega_p} \tag{\theequation (b)}\label{FPSI_FEM:eq:splitChiTheta_b}\\
& \hspace{10mm}- \alpha (\chi_{pp,h}^{n+1}, \nabla \cdot \bm{\varphi}_h)_{\Omega_p}- \langle \chi_{g1,h}^{n+1}\bm{n}_p, \bm{\varphi}_h\rangle_\gamma +  ( \chi_{g2,h}^{n+1}\bm{\tau}_\gamma, \bm{\varphi}_h)_\gamma =  2  \nu_p \left( D\left(\bm{\theta}_{\eta}^{n+1}\right), D(\bm{\varphi}_h) \right)_{\Omega_p}  \notag \\
& \hspace{10mm}+ \lambda  \left( \nabla \cdot  \bm{\theta}_{\eta}^{n+1}, \nabla \cdot \bm{\varphi}_h \right)_{\Omega_p} - \alpha (\theta_{pp}^{n+1}, \nabla \cdot \bm{\varphi}_h)_{\Omega_p} -\langle \theta_{g1}^{n+1}\bm{n}_p, \bm{\varphi}_h\rangle_\gamma +  ( \theta_{g2}^{n+1}\bm{\tau}_\gamma, \bm{\varphi}_h)_\gamma \notag \\
&\frac{s_0}{\Delta t} ( \chi_{pp,h}^{n+1} - \chi_{pp,h}^n, w_h)_{\Omega_p} + \alpha(\nabla \cdot \bm{\dot{\chi}}_{\eta,h}^{n+1},w_h)_{\Omega_p} + \kappa  (\nabla \chi_{pp,h}^{n+1}, \nabla w_h)_{\Omega_p} -  ( \chi_{\lambda p,h}^{n+1}, w_h)_{\gamma} \tag{\theequation (c)}\label{FPSI_FEM:eq:splitChiTheta_c} \\
&\hspace{10mm}=  \alpha(\nabla \cdot \bm{\dot{\theta}}_{\eta}^{n+1},w_h)_{\Omega_p} + \kappa  (\nabla \theta_{pp}^{n+1}, \nabla w_h)_{\Omega_p} -  ( \theta_{\lambda p}^{n+1}, w_h)_{\gamma}  \notag \\
  &  \frac{1}{\beta } ( \chi_{g2,h}^{n+1},s_{2,h})_\gamma  +  (\bm{\chi}_{u,h}^{n+1} \cdot \bm{\tau}_\gamma,s_{2,h})_\gamma - (\bm{\dot{\chi}}_{\eta,h}^{n+1} \cdot \bm{\tau}_\gamma,s_{2,h})_\gamma \tag{\theequation (d)}\label{FPSI_FEM:eq:splitChiTheta_d} \\
  & \hspace{10mm}=    (\bm{\theta}_{u}^{n+1} \cdot \bm{\tau}_\gamma,s_{2,h})_\gamma - (\bm{\dot{\theta}}_{\eta}^{n+1} \cdot \bm{\tau}_\gamma,s_{2,h})_\gamma  \notag \\
      &   \langle \chi_{g1,h}^{n+1},\mu_h\rangle_\gamma   + (  \chi_{pp,h}^{n+1},\mu_h)_\gamma + \overline{\epsilon}( \chi_{\lambda p,h}^{n+1},\mu_h)_{1/2,\gamma} \tag{\theequation (e)}\label{FPSI_FEM:eq:splitChiTheta_e} \\
      &\hspace{10mm} =  \langle \theta_{g1}^{n+1},\mu_h\rangle_\gamma + (  \theta_{pp}^{n+1},\mu_h)_\gamma + \overline{\epsilon} (\theta_{\lambda p}^{n+1},\mu_h)_{1/2,\gamma}  \notag \\
      & (\nabla \cdot \bm{\chi}_{u,h}^{n+1},q_h)_{\Omega_f} =  (\nabla \cdot \bm{\theta}_{u}^{n+1},q_h)_{\Omega_f}  \tag{\theequation (f)}\label{FPSI_FEM:eq:splitChiTheta_f}\\
    & \langle \bm{\chi}_{u,h}^{n+1} \cdot \bm{n}_f,s_{1,h}\rangle_\gamma + \langle \bm{\dot{\chi}}_{\eta,h}^{n+1} \cdot \bm{n}_p ,s_{1,h} \rangle_\gamma   -  \langle \chi_{\lambda p,h}^{n+1},s_{1,h}\rangle_\gamma \tag{\theequation (g)}\label{FPSI_FEM:eq:splitChiTheta_g}  \\
    &\hspace{10mm}= \langle \bm{\theta}_{u}^{n+1} \cdot \bm{n}_f,s_{1,h}\rangle_\gamma +\langle \bm{\dot{\theta}}_{\eta}^{n+1} \cdot \bm{n}_p ,s_{1,h} \rangle_\gamma   -  \langle \theta_{\lambda p}^{n+1},s_{1,h}\rangle_\gamma.  \notag   
 %   \end{split}
\end{align}

\noindent Note that the terms $(\bm{\dot{\theta}}_\eta^{n+1} - \bm{\dot{\theta}}_\eta^{n},\bm{\varphi}_h),$ $(\theta_{pp}^{n+1}-\theta_{pp}^n,w_h)$, and $(\theta_{g2}^{n+1},s_{2,h})_{0,\gamma}$ are zero by properties of the $L^2$ projection operators. Pick $(\bm{\chi}_{u,h}^{n+1}, \chi_{pf,h}^{n+1}, \bm{\dot{\chi}}_{\eta,h}^{n+1}, \chi_{pp,h}^{n+1}, \chi_{g1,h}^{n+1}, \chi_{g2,h}^{n+1},\chi_{\lambda p,h}^{n+1})$ $\in A^h$ as test functions in the above system. Simplifying the inner products of differencs using the identity $(a-b)a = \frac{1}{2}(a^2 - b^2 + (a-b)^2)$ and noting $ (\nabla \cdot \bm{\theta}_{u}^{n+1},\chi_{pf,h}^{n+1})_{\Omega_f} = 0$ by the properties of the Stokes projection operator \eqref{FPSI_FEM:eq:StokesProj_Eqs} results in
\begin{align*}
    \begin{split}
       &\frac{\rho_f}{2 \Delta t} \left( || \bm{\chi}_{u,h}^{n+1} ||_0^2 - ||\bm{\chi}_{u,h}^{n} ||_0^2 + || \bm{\chi}_{u,h}^{n+1} - \bm{\chi}_{u,h}^{n}||_0^2 \right) + 2  \nu_f || D(\bm{\chi}_{u,h}^{n+1}) ||_0^2 -  (\chi_{pf,h}^{n+1}, \nabla \cdot \bm{\chi}_{u,h}^{n+1})_{\Omega_f}  \\
       &\hspace{8mm}-  \langle \chi_{g1,h}^{n+1} \bm{n_f},\bm{\chi}_{u,h}^{n+1} \rangle_\gamma  -  ( \chi_{g2,h}^{n+1} \bm{\tau_\gamma},\bm{\chi}_{u,h}^{n+1})_\gamma  =  \frac{\rho_f}{\Delta t} (\bm{\theta}_u^{n+1}- \bm{\theta}_u^{n},\bm{\chi}_{u,h}^{n+1})_{\Omega_f} \\
       &\hspace{8mm}+ 2  \nu_f ( D(\bm{\theta}_{u}^{n+1}), D(\bm{\chi}_{u,h}^{n+1}) )_{\Omega_f} -  (\theta_{pf}^{n+1}, \nabla \cdot \bm{\chi}_{u,h}^{n+1})_{\Omega_f} -  \langle \theta_{g1}^{n+1} \bm{n_f},\bm{\chi}_{u,h}^{n+1} \rangle_\gamma  -  ( \theta_{g2}^{n+1} \bm{\tau_\gamma},\bm{\chi}_{u,h}^{n+1})_\gamma \\
&\frac{\rho_p}{2 \Delta t} \left( || \bm{\dot{\chi}}_{\eta,h}^{n+1} ||_0^2 - ||\bm{\dot{\chi}}_{\eta,h}^n||_0^2 + || \bm{\dot{\chi}}_{\eta,h}^{n+1} - \bm{\dot{\chi}}_{\eta,h}^{n} ||_0^2 \right) +  2  \nu_p ( D(\bm{\chi}_{\eta,h}^{n+1}),D(\bm{\dot{\chi}}_{\eta,h}^{n+1}))_{\Omega_p} \\
& \hspace{8mm}+ \lambda (\nabla \cdot  \bm{\chi}_{\eta,h}^{n+1},\nabla \cdot \bm{\dot{\chi}}_{\eta,h}^{n+1})_{\Omega_p}- \alpha (\chi_{pp,h}^{n+1}, \nabla \cdot \bm{\dot{\chi}}_{\eta,h}^{n+1})_{\Omega_p}- \langle \chi_{g1,h}^{n+1}\bm{n_p}, \bm{\dot{\chi}}_{\eta,h}^{n+1}\rangle_\gamma +  ( \chi_{g2,h}^{n+1}\bm{\tau_\gamma}, \bm{\dot{\chi}}_{\eta,h}^{n+1})_\gamma \\
&\hspace{8mm}=   2  \nu_p ( D(\bm{\theta}_{\eta}^{n+1}), D(\bm{\dot{\chi}}_{\eta,h}^{n+1}) )_{\Omega_p}  + \lambda  ( \nabla \cdot  \bm{\theta}_{\eta}^{n+1}, \nabla \cdot \bm{\dot{\chi}}_{\eta,h}^{n+1} )_{\Omega_p} - \alpha (\theta_{pp}^{n+1}, \nabla \cdot \bm{\dot{\chi}}_{\eta,h}^{n+1})_{\Omega_p} \\
&\hspace{8mm}- \langle \theta_{g1}^{n+1}\bm{n_p}, \bm{\dot{\chi}}_{\eta,h}^{n+1}\rangle_\gamma+  ( \theta_{g2}^{n+1}\bm{\tau_\gamma}, \bm{\dot{\chi}}_{\eta,h}^{n+1})_\gamma 
\end{split}
\end{align*}
\begin{align*}
\begin{split}
&\frac{s_0}{2 \Delta t} \left( ||\chi_{pp,h}^{n+1}||_0^2 - ||\chi_{pp,h}^n||_0^2 + ||\chi_{pp,h}^{n+1}-\chi_{pp,h}^n||_0^2 \right) + \alpha(\nabla \cdot \bm{\dot{\chi}}_{\eta,h}^{n+1},\chi_{pp,h}^{n+1})_{\Omega_p} + \kappa  ||\nabla \chi_{pp,h}^{n+1}||_0^2 \\
&\hspace{8mm}-  ( \chi_{\lambda p,h}^{n+1}, \chi_{pp,h}^{n+1})_{\gamma}= \alpha(\nabla \cdot \bm{\dot{\theta}}_{\eta}^{n+1},\chi_{pp,h}^{n+1})_{\Omega_p} + \kappa  (\nabla \theta_{pp}^{n+1}, \nabla \chi_{pp,h}^{n+1})_{\Omega_p} -  ( \theta_{\lambda p}^{n+1}, \chi_{pp,h}^{n+1})_{\gamma}  \\
  &  \frac{1}{\beta } || \chi_{g2,h}^{n+1}||_{0,\gamma}^2  +  (\bm{\chi}_{u,h}^{n+1} \cdot \bm{\tau_\gamma},\chi_{g2,h}^{n+1})_\gamma - (\bm{\dot{\chi}}_{\eta,h}^{n+1} \cdot \bm{\tau_\gamma},\chi_{g2,h}^{n+1})_\gamma =    (\bm{\theta}_{u}^{n+1} \cdot \bm{\tau_\gamma},\chi_{g2,h}^{n+1})_\gamma - (\bm{\dot{\theta}}_{\eta}^{n+1} \cdot \bm{\tau_\gamma},\chi_{g2,h}^{n+1})_\gamma \\
      &   \langle \chi_{g1,h}^{n+1},\chi_{\lambda p,h}^{n+1}\rangle_\gamma   + (  \chi_{pp,h}^{n+1},\chi_{\lambda p,h}^{n+1})_\gamma + \overline{\epsilon}|| \chi_{\lambda p,h}^{n+1}||^2_{1/2,\gamma} \\
      &\hspace{8mm}=   \langle \theta_{g1}^{n+1},\chi_{\lambda p,h}^{n+1}\rangle_\gamma+ (  \theta_{pp}^{n+1},\chi_{\lambda p,h}^{n+1})_\gamma + \overline{\epsilon} (\theta_{\lambda p}^{n+1},\chi_{\lambda p,h}^{n+1})_{1/2,\gamma} \\
        & (\nabla \cdot \bm{\chi}_{u,h}^{n+1},\chi_{pf,h}^{n+1})_{\Omega_f} =  (\nabla \cdot \bm{\theta}_{u}^{n+1},\chi_{pf,h}^{n+1})_{\Omega_f}\\
    & \langle \bm{\chi}_{u,h}^{n+1} \cdot \bm{n_f},\chi_{g1,h}^{n+1}\rangle_\gamma + \langle \bm{\dot{\chi}}_{\eta,h}^{n+1} \cdot \bm{n_p} ,\chi_{g1,h}^{n+1} \rangle_\gamma   -  \langle \chi_{\lambda p,h}^{n+1},\chi_{g1,h}^{n+1}\rangle_\gamma = \\
    &\hspace{8mm} \langle \bm{\theta}_{u}^{n+1} \cdot \bm{n_f},\chi_{g1,h}^{n+1}\rangle_\gamma + \langle \bm{\dot{\theta}}_{\eta}^{n+1} \cdot \bm{n_p} ,\chi_{g1,h}^{n+1} \rangle_\gamma   -  \langle \theta_{\lambda p}^{n+1},\chi_{g1,h}^{n+1}\rangle_\gamma.  
    \end{split}
\end{align*}

We now sum the equations, and after cancellation of many of the mixed terms involving truncation errors, we obtain:
\begin{align}\label{FPSI_FEM:eq:ErrorInnProducts}
    \begin{split}
       &\frac{\rho_f}{2 \Delta t} \left( || \bm{\chi}_{u,h}^{n+1} ||_0^2 - ||\bm{\chi}_{u,h}^{n} ||_0^2 + || \bm{\chi}_{u,h}^{n+1} - \bm{\chi}_{u,h}^{n}||_0^2 \right) + \frac{\rho_p}{2 \Delta t} \left( || \bm{\dot{\chi}}_{\eta,h}^{n+1} ||_0^2 - ||\bm{\dot{\chi}}_{\eta,h}^n||_0^2 + || \bm{\dot{\chi}}_{\eta,h}^{n+1} - \bm{\dot{\chi}}_{\eta,h}^{n} ||_0^2 \right)\\
       & \hspace{3mm} + \frac{s_0}{2 \Delta t} \left( ||\chi_{pp,h}^{n+1}||_0^2 - ||\chi_{pp,h}^n||_0^2 + ||\chi_{pp,h}^{n+1}-\chi_{pp,h}^n||_0^2 \right) + \frac{1}{\beta } || \chi_{g2,h}^{n+1}||_{0,\gamma}^2 + \overline{\epsilon}|| \chi_{\lambda p,h}^{n+1}||^2_{1/2,\gamma}\\
       & \hspace{3mm}+ 2  \nu_f || D(\bm{\chi}_{u,h}^{n+1}) ||_0^2  + \kappa  ||\nabla \chi_{pp,h}^{n+1}||_0^2     +  \frac{\nu_p}{\Delta t} \left(  ||  D(\bm{\chi}_{\eta,h}^{n+1})||_0^2 - || D(\bm{\chi}_{\eta,h}^{n})||_0^2 \right) + \Delta t \nu_p || D(\bm{\dot{\chi}}_{\eta,h}^{n+1})||_0^2 \\
       &\hspace{3mm}+   \frac{\lambda}{2\Delta t} \left( ||\nabla \cdot \bm{\chi}_{\eta,h}^{n+1}||_0^2 - ||\nabla \cdot \bm{\chi}_{\eta,h}^{n}||_0^2  \right) + \frac{\lambda \Delta t }{2} || \nabla \cdot \bm{\dot{\chi}}_{\eta,h}^{n+1} ||_0^2      \\     
     &  =    \frac{\rho_f}{\Delta t} (\bm{\theta}_u^{n+1}- \bm{\theta}_u^{n},\bm{\chi}_{u,h}^{n+1})_{\Omega_f} + 2  \nu_f ( D(\bm{\theta}_{u}^{n+1}), D(\bm{\chi}_{u,h}^{n+1}) )_{\Omega_f}   -  (\theta_{pf}^{n+1}, \nabla \cdot \bm{\chi}_{u,h}^{n+1})_{\Omega_f} \\
     &\hspace{3mm} -  \langle \theta_{g1}^{n+1} \bm{n}_f,\bm{\chi}_{u,h}^{n+1} \rangle_\gamma   -  ( \theta_{g2}^{n+1} \bm{\tau}_\gamma,\bm{\chi}_{u,h}^{n+1})_\gamma  + 2  \nu_p ( D(\bm{\theta}_{\eta}^{n+1}), D(\bm{\dot{\chi}}_{\eta,h}^{n+1}))_{\Omega_p}  \\
&\hspace{3mm} + \lambda  ( \nabla \cdot  \bm{\theta}_{\eta}^{n+1}, \nabla \cdot \bm{\dot{\chi}}_{\eta,h}^{n+1} )_{\Omega_p} - \alpha (\theta_{pp}^{n+1}, \nabla \cdot \bm{\dot{\chi}}_{\eta,h}^{n+1})_{\Omega_p}  - \langle \theta_{g1}^{n+1}\bm{n}_p, \bm{\dot{\chi}}_{\eta,h}^{n+1}\rangle_\gamma +  ( \theta_{g2}^{n+1}\bm{\tau}_\gamma, \bm{\dot{\chi}}_{\eta,h}^{n+1})_\gamma \\
& \hspace{3mm} + \alpha(\nabla \cdot \bm{\dot{\theta}}_{\eta}^{n+1},\chi_{pp,h}^{n+1})_{\Omega_p}  +\kappa  (\nabla \theta_{pp}^{n+1}, \nabla \chi_{pp,h}^{n+1})_{\Omega_p} -  ( \theta_{\lambda p}^{n+1}, \chi_{pp,h}^{n+1})_{\gamma} +   (\bm{\theta}_{u}^{n+1} \cdot \bm{\tau}_\gamma,\chi_{g2,h}^{n+1})_\gamma \\
&\hspace{3mm}- (\bm{\dot{\theta}}_{\eta}^{n+1} \cdot \bm{\tau}_\gamma,\chi_{g2,h}^{n+1})_\gamma   +   \langle \theta_{g1}^{n+1},\chi_{\lambda p,h}^{n+1}\rangle_\gamma +  (  \theta_{pp}^{n+1},\chi_{\lambda p,h}^{n+1})_\gamma + \overline{\epsilon} (\theta_{\lambda p}^{n+1},\chi_{\lambda p,h}^{n+1})_{1/2,\gamma} \\
&\hspace{3mm} + \langle \bm{\theta}_{u}^{n+1} \cdot \bm{n}_f,\chi_{g1,h}^{n+1}\rangle_\gamma  + \langle \bm{\dot{\theta}}_{\eta}^{n+1} \cdot \bm{n}_p ,\chi_{g1,h}^{n+1} \rangle_\gamma   -  \langle \theta_{\lambda p}^{n+1},\chi_{g1,h}^{n+1}\rangle_\gamma\\  
& := \sum_{j=1}^{21} R_j,
    \end{split}
\end{align}
where the inner products appearing on the left hand side were rewritten by:
\begin{align*}
      2  \nu_p ( D(\bm{\chi}_{\eta,h}^{n+1}),D(\bm{\dot{\chi}}_{\eta,h}^{n+1}))_{\Omega_p} &=   \frac{\nu_p}{\Delta t} \left(  ||  D(\bm{\chi}_{\eta,h}^{n+1})||_0^2 - || D(\bm{\chi}_{\eta,h}^{n})||_0^2 \right) + \Delta t \nu_p || D(\bm{\dot{\chi}}_{\eta,h}^{n+1})||_0^2\\
   \lambda (\nabla \cdot \bm{\chi}_{\eta,h}^{n+1},\nabla \cdot \bm{\dot{\chi}}_{\eta,h}^{n+1}) &= \frac{\lambda}{2\Delta t} \left( ||\nabla \cdot \bm{\chi}_{\eta,h}^{n+1}||_0^2 - ||\nabla \cdot \bm{\chi}_{\eta,h}^{n}||_0^2  \right) + \frac{\lambda \Delta t }{2} || \nabla \cdot \bm{\dot{\chi}}_{\eta,h}^{n+1} ||_0^2.
\end{align*}
\begin{comment}
\noindent Denote each inner product on the right hand side of \eqref{FPSI_FEM:eq:ErrorInnProducts} by $R_j$, for $j=1,\ldots,21$. With this notation, \eqref{FPSI_FEM:eq:ErrorInnProducts} is equivalent to 
\begin{align}\label{FPSI_FEM:Eq:ErrorLRSums}
      \begin{split}
       &\frac{\rho_f}{2 \Delta t} \left( || \bm{\chi}_{u,h}^{n+1} ||_0^2 - ||\bm{\chi}_{u,h}^{n} ||_0^2 + || \bm{\chi}_{u,h}^{n+1} - \bm{\chi}_{u,h}^{n}||_0^2 \right) + \frac{\rho_p}{2 \Delta t} \left( || \bm{\dot{\chi}}_{\eta,h}^{n+1} ||_0^2 - ||\bm{\dot{\chi}}_{\eta,h}^n||_0^2 + || \bm{\dot{\chi}}_{\eta,h}^{n+1} - \bm{\dot{\chi}}_{\eta,h}^{n} ||_0^2 \right)\\
       & \hspace{3mm} + \frac{s_0}{2 \Delta t} \left( ||\chi_{pp,h}^{n+1}||_0^2 - ||\chi_{pp,h}^n||_0^2 + ||\chi_{pp,h}^{n+1}-\chi_{pp,h}^n||_0^2 \right) + \frac{1}{\beta } || \chi_{g2,h}^{n+1}||_{0,\gamma}^2  + \overline{\epsilon}|| \chi_{\lambda p,h}^{n+1}||^2_{1/2,\gamma} \\
       &\hspace{3mm}+ 2  \nu_f || D(\bm{\chi}_{u,h}^{n+1}) ||_0^2  + \kappa  ||\nabla \chi_{pp,h}^{n+1}||_0^2 +   \frac{\nu_p}{\Delta t} \left(  ||  D(\bm{\chi}_{\eta,h}^{n+1})||_0^2 - || D(\bm{\chi}_{\eta,h}^{n})||_0^2 \right) \\
       &\hspace{3mm} + \frac{\lambda}{2\Delta t} \left( ||\nabla \cdot \bm{\chi}_{\eta,h}^{n+1}||_0^2 - ||\nabla \cdot \bm{\chi}_{\eta,h}^{n}||_0^2  \right)  + \Delta t \nu_p || D(\bm{\dot{\chi}}_{\eta,h}^{n+1})||_0^2 + \frac{\lambda \Delta t }{2} || \nabla \cdot \bm{\dot{\chi}}_{\eta,h}^{n+1} ||_0^2   =  \sum_{j=1}^{21} R_j.  
    \end{split}
\end{align}
\end{comment}
Each inner product $R_j$, $j=1,2,\ldots,21$ on the right hand side of \eqref{FPSI_FEM:eq:ErrorInnProducts} may be bounded using Cauchy-Schwartz, Young's inequality, and the inequalities \eqref{FPSI_FEM:eq:TraceIneqs}-\eqref{FPSI_FEM:eq:1normbound}, defining $C_1 := C_T \sqrt{C_{KP}}$. We also employ the inequality $||\cdot||_{1/2,\gamma} \leq ||\cdot||_{1/2,\Gamma^r} \leq || \cdot ||_{1,\Omega_r}$, for norms of $\bm{\chi}_{u,h}, \bm{\chi}_{\eta,h}$, $\bm{\theta}_u,$ and $\bm{\theta}_\eta$, and $r \in \{f,p\}$.
\begin{align*}
 R_{1} &:= \frac{\rho_f}{\Delta t} (\bm{\theta}_u^{n+1}- \bm{\theta}_u^{n},\bm{\chi}_{u,h}^{n+1})_{\Omega_f} \leq  \rho_f \left( \frac{\rho_f C_P^2 C_K^2}{\nu_f} || \bm{\dot{\theta}}_u^{n+1}||_0^2 + \frac{\nu_f }{4\rho_f C_P^2 C_K^2} || \bm{\chi}_{u,h}^{n+1}||_0^2 \right) \\
 &\hspace{15mm} \leq  \frac{\rho_f^2 C_P^2 C_K^2}{\nu_f} || \bm{\dot{\theta}}_u^{n+1}||_0^2 + \frac{\nu_f }{4 } || D(\bm{\chi}_{u,h}^{n+1})||_0^2 \\
 %\frac{\rho_f}{\Delta t} (\bm{\theta}_u^{n+1}- \bm{\theta}_u^{n},\bm{\chi}_{u,h}^{n+1})_{\Omega_f} \leq \frac{\rho_f}{\Delta t} \left( ||\bm{\theta}_u^{n+1}||_0 || \bm{\chi}_{u,h}^{n+1}||_0 + ||\bm{\theta}_u^n||_0 || \bm{\chi}_{u,h}^{n+1}||_0 \right)   \\
 %   &\hspace{20mm} \leq  \frac{\rho_f}{\Delta t} \Big( \frac{2C_P^2C_K^2 \rho_f}{\Delta t \nu_f} || \bm{\theta}_u^{n+1}||_0^2 + \frac{\Delta t \nu_f}{8C_P^2 C_K^2 \rho_f} || \bm{\chi}_{u,h}^{n+1}||_0^2 + \frac{2C_P^2C_K^2 \rho_f}{\Delta t \nu_f} || \bm{\theta}_u^{n}||_0^2 + \frac{\Delta t \nu_f}{8C_P^2 C_K^2 \rho_f} || \bm{\chi}_{u,h}^{n+1}||_0^2 \Big) \\ 
 %   &\hspace{15mm}\leq \frac{\rho_f}{\Delta t} \left(  \frac{2C_P^2 C_K^2 \rho_f}{\Delta t \nu_f}||\bm{\theta}_u^{n+1}||_0^2 + \frac{2C_P^2 C_K^2 \rho_f}{\Delta t \nu_f} ||\bm{\theta}_u^n||_0^2 \right) + \frac{\nu_f}{4}|| D(\bm{\chi}_{u,h}^{n+1})||_0^2 \\
    R_2 &:= 2  \nu_f ( D(\bm{\theta}_{u}^{n+1}), D(\bm{\chi}_{u,h}^{n+1}) )_{\Omega_f} \leq 2\nu_f \left( 4 || D(\bm{\theta}_{u}^{n+1})||_0^2 + \frac{1}{16} || D(\bm{\chi}_{u,h}^{n+1}) ||_0^2  \right)\\
    &\hspace{15mm} \leq  8\nu_f || \bm{\theta}_{u}^{n+1}||_1^2 + \frac{\nu_f}{8} || D(\bm{\chi}_{u,h}^{n+1}) ||_0^2   \\
     R_{3} &:= -  (\theta_{pf}^{n+1}, \nabla \cdot \bm{\chi}_{u,h}^{n+1})_{\Omega_f} \leq \frac{2dC_K^2}{\nu_f}||\theta_{pf}^{n+1}||_0^2 + \frac{\nu_f}{8d C_K^2} ||\nabla \cdot \bm{\chi}_{u,h}^{n+1}||_0^2 \\
     &\hspace{15mm} \leq \frac{2d C_K^2}{\nu_f}||\theta_{pf}^{n+1}||_0^2 + \frac{\nu_f}{8} ||D(\bm{\chi}_{u,h}^{n+1})||_0^2\\
    R_{4} &:= -  \langle \theta_{g1}^{n+1} \bm{n}_f,\bm{\chi}_{u,h}^{n+1} \rangle_\gamma \leq \frac{C_{KP}}{\nu_f} || \theta_{g1}^{n+1}||_{-1/2,\gamma}^2 + \frac{\nu_f}{4C_{KP}} ||\bm{\chi}_{u,h}^{n+1}||_{1/2,\gamma}^2 \\
    &\hspace{15mm}  \leq \frac{C_{KP}}{\nu_f} || \theta_{g1}^{n+1}||_{-1/2,\gamma}^2 + \frac{\nu_f}{4C_{KP}} ||\bm{\chi}_{u,h}^{n+1}||_{1}^2  \leq \frac{C_{KP}}{\nu_f} || \theta_{g1}^{n+1}||_{-1/2,\gamma}^2 + \frac{\nu_f}{4} ||D(\bm{\chi}_{u,h}^{n+1})||_{0}^2\\
      R_5 &:= -  ( \theta_{g2}^{n+1} \bm{\tau}_\gamma,\bm{\chi}_{u,h}^{n+1})_\gamma \leq \frac{C_1^2}{\nu_f} ||\theta_{g2}^{n+1}||_{0,\gamma}^2 + \frac{\nu_f}{4 C_1^2} ||\bm{\chi}_{u,h}^{n+1}||_{0,\gamma}^2 \leq \frac{C_1^2}{ \nu_f} || \theta_{g2}^{n+1}||_{0,\gamma}^2 + \frac{\nu_f}{4} || D(\bm{\chi}_{u,h}^{n+1})||_0^2     \\     
    R_6 &:= 2  \nu_p \left( D\left(\bm{\theta}_{\eta}^{n+1}\right), D(\bm{\dot{\chi}}_{\eta,h}^{n+1}) \right)_{\Omega_p} \leq 2\nu_p \left(  \frac{4}{\Delta t} || D\left(\bm{\theta}_{\eta}^{n+1}\right) ||_0^2 + \frac{\Delta t}{16} || D(\bm{\dot{\chi}}_{\eta,h}^{n+1})||_0^2 \right) \\
    & \hspace{15mm} \leq    \frac{8\nu_p}{\Delta t} || \bm{\theta}_{\eta}^{n+1} ||_1^2 + \frac{\Delta t \nu_p}{8} || D(\bm{\dot{\chi}}_{\eta,h}^{n+1})||_0^2 \\
    R_7 &:= \lambda  ( \nabla \cdot  \bm{\theta}_{\eta}^{n+1}, \nabla \cdot \bm{\dot{\chi}}_{\eta,h}^{n+1} )_{\Omega_p} \leq \lambda \left( \frac{2}{\Delta t} ||\nabla \cdot  \bm{\theta}_{\eta}^{n+1}||_0^2 + \frac{\Delta t}{8} || \nabla \cdot \bm{\dot{\chi}}_{\eta,h}^{n+1}||_0^2 \right) \\
    &\hspace{15mm} \leq \lambda \left( \frac{2dC_K^2}{\Delta t} ||\bm{\theta}_{\eta}^{n+1}||_1^2 + \frac{\Delta t}{8} || \nabla \cdot \bm{\dot{\chi}}_{\eta,h}^{n+1}||_0^2 \right)\\
    R_8 &:= - \alpha (\theta_{pp}^{n+1}, \nabla \cdot \bm{\dot{\chi}}_{\eta,h}^{n+1})_{\Omega_p} \leq \alpha \left( \frac{2\alpha}{\lambda \Delta t} || \theta_{pp}^{n+1}||_0^2 + \frac{\lambda \Delta t}{8\alpha} || \nabla \cdot \bm{\dot{\chi}}^{n+1}_{\eta,h} ||_0^2 \right)\\
    R_{9} &:= - \langle \theta_{g1}^{n+1}\bm{n}_p, \bm{\dot{\chi}}_{\eta,h}^{n+1}\rangle_\gamma  \leq \frac{C_{KP}}{\Delta t \nu_p} || \theta_{g1}^{n+1}||_{-1/2,\gamma}^2 + \frac{\Delta t \nu_p}{4} ||D(\bm{\dot{\chi}}_{\eta,h}^{n+1})||_{0}^2\\
      R_{10} &:=  ( \theta_{g2}^{n+1}\bm{\tau}_\gamma, \bm{\dot{\chi}}_{\eta,h}^{n+1})_\gamma \leq \frac{2C_1^2}{ \Delta t \nu_p}||\theta_{g2}^{n+1}||_{0,\gamma}^2 + \frac{\Delta t \nu_p}{8C_1^2} || \bm{\dot{\chi}}_{\eta,h}^{n+1}||_{0,\gamma}^2 \\
      &\hspace{15mm} \leq \frac{2C_1^2}{ \Delta t \nu_p}||\theta_{g2}^{n+1}||_{0,\gamma}^2 + \frac{\Delta t \nu_p}{8} || D(\bm{\dot{\chi}}_{\eta,h}^{n+1})||_0^2 \\
    R_{11} &:=  \alpha(\nabla \cdot \bm{\dot{\theta}}_{\eta}^{n+1},\chi_{pp,h}^{n+1})_{\Omega_p} \leq \alpha \left( \frac{2\alpha C_P^2}{\kappa}||\nabla \cdot \bm{\dot{\theta}}_{\eta}^{n+1}||_0^2 + \frac{\kappa}{8\alpha C_P^2}||\chi_{pp,h}^{n+1}||_0^2    \right) \\
    &\hspace{15mm} \leq \frac{2\alpha^2 C_P^2 d C_K^2}{\kappa}|| \bm{\dot{\theta}}_{\eta}^{n+1}||_1^2 + \frac{\kappa}{8} ||\nabla \chi_{pp,h}^{n+1}||_0^2   
                \end{align*}
    \begin{align*}
    R_{12} &:= \kappa  (\nabla \theta_{pp}^{n+1}, \nabla \chi_{pp,h}^{n+1})_{\Omega_p} \leq \kappa \left( 2 || \nabla \theta_{pp}^{n+1}||_0^2 + \frac{1}{8}|| \nabla \chi_{pp,h}^{n+1}||_0^2 \right) \leq \kappa \left( 2 || \theta_{pp}^{n+1}||_1^2 + \frac{1}{8}|| \nabla \chi_{pp,h}^{n+1}||_0^2 \right) \\
    R_{13} &:= -  ( \theta_{\lambda p}^{n+1}, \chi_{pp,h}^{n+1})_{\gamma} \leq \frac{C_1^2}{\kappa} ||\theta_{\lambda p}^{n+1}||_{0,\gamma}^2 + \frac{\kappa}{4C_1^2} || \chi_{pp,h}^{n+1}||_{0,\gamma}^2 \leq \frac{C_1^2}{ \kappa} ||\theta_{\lambda p}^{n+1}||_{1/2,\gamma}^2 + \frac{\kappa}{4} || \nabla \chi_{pp,h}^{n+1}||_0^2\\
    R_{14} &:=    (\bm{\theta}_{u}^{n+1}  ,\chi_{g2,h}^{n+1} \bm{\tau}_\gamma)_\gamma \leq \beta ||\bm{\theta}_u^{n+1}||_{0,\gamma}^2 + \frac{1}{4\beta} ||\chi_{g2,h}^{n+1}||_{0,\gamma}^2 \leq C_T^2 \beta  || \bm{\theta}_{u}^{n+1}||_1^2 + \frac{1}{4\beta } ||\chi_{g2,h}^{n+1}||_{0,\gamma}^2 \\
    R_{15} &:= - (\bm{\dot{\theta}}_{\eta}^{n+1} ,\chi_{g2,h}^{n+1}  \bm{\tau}_\gamma)_\gamma \leq C_T^2 \beta||\bm{\dot{\theta}}_{\eta}^{n+1}||_1^2 + \frac{1}{4\beta}||\chi_{g2,h}^{n+1}||_{0,\gamma}^2 \\
     R_{16} &:=  \langle \theta_{g1}^{n+1},\chi_{\lambda p,h}^{n+1}\rangle_\gamma \leq \frac{2}{\overline{\epsilon}} ||\theta_{g1}^{n+1}||_{-1/2,\gamma}^2 + \frac{\overline{\epsilon}}{8} ||\chi_{\lambda p,h}^{n+1}||_{1/2,\gamma}^2 \\
    R_{17} &:=  (  \theta_{pp}^{n+1},\chi_{\lambda p,h}^{n+1})_\gamma \leq \frac{1}{\overline{\epsilon}} ||  \theta_{pp}^{n+1}||_{0,\gamma}^2 + \frac{\overline{\epsilon}}{4} ||\chi_{\lambda p,h}^{n+1}||_{0,\gamma}^2 \leq  \frac{C_T^2}{\overline{\epsilon}} || \theta_{pp}^{n+1}||_1^2 + \frac{\overline{\epsilon}}{4} ||\chi_{\lambda p,h}^{n+1}||_{1/2,\gamma}^2\\
    R_{18} &:= \overline{\epsilon} (\theta_{\lambda p}^{n+1},\chi_{\lambda p,h}^{n+1})_{1/2,\gamma} \leq \overline{\epsilon} \left( 2 || \theta_{\lambda p}^{n+1}||_{1/2,\gamma}^2 + \frac{1}{8} ||\chi_{\lambda p,h}^{n+1}||_{1/2,\gamma}^2   \right)\\
    R_{19} &:=  \langle \bm{\theta}_{u}^{n+1} ,\chi_{g1,h}^{n+1}  \bm{n}_f\rangle_\gamma \leq \frac{1}{2\varepsilon_{19}} ||\bm{\theta}_{u}^{n+1}||_{1/2,\gamma}^2 + \frac{\varepsilon_{19}}{2} || \chi_{g1,h}^{n+1}||_{-1/2,\gamma}^2   \\
    &\hspace{15mm} \leq \frac{1}{2\varepsilon_{19}} ||\bm{\theta}_{u}^{n+1}||_{1}^2 + \frac{\varepsilon_{19}}{2} || \chi_{g1,h}^{n+1}||_{-1/2,\gamma}^2  \\
    R_{20} &:= \langle \bm{\dot{\theta}}_{\eta}^{n+1} ,\chi_{g1,h}^{n+1}  \bm{n}_p \rangle_\gamma  \leq  \frac{1}{2\varepsilon_{20}} || \bm{\dot{\theta}}_{\eta}^{n+1}||_{1}^2 + \frac{\varepsilon_{20}}{2} || \chi_{g1,h}^{n+1}||_{-1/2,\gamma}\\
    R_{21} &:= -  \langle \theta_{\lambda p}^{n+1},\chi_{g1,h}^{n+1}\rangle_\gamma \leq \frac{1}{2\varepsilon_{21}} ||\theta_{\lambda p}^{n+1}||_{1/2,\gamma}^2 + \frac{\varepsilon_{21}}{2} || \chi_{g1,h}^{n+1}||_{-1/2,\gamma}^2. 
\end{align*}

\noindent Combining these, taking $d=2$ since we work in the 2D case, yields 
\begin{align}\label{FPSI_FEM:eq:Err_RHSbounds}
    \begin{split}
        \sum_{j=1}^{21} &R_j \leq \left( 8\nu_f + C_T^2  \beta + \frac{1}{2\varepsilon_{19}}\right) || \bm{\theta}_{u}^{n+1}||_1^2 +  \frac{\rho_f^2 C_P^2 C_K^2}{ \nu_f} ||\bm{\dot{\theta}}_u^{n+1}||_0^2 +  \left( \frac{8\nu_p}{\Delta t}+ \frac{4\lambda C_K^2}{\Delta t} \right) || \bm{\theta}_{\eta}^{n+1} ||_1^2  \\
         &+ \left( \frac{4\alpha^2  C_P^2 C_K^2}{\kappa} + C_T^2 \beta+ \frac{1}{2\varepsilon_{20}} \right) || \bm{\dot{\theta}}_{\eta}^{n+1}||_{1}^2+ \frac{2\alpha^2}{\lambda \Delta t} ||\theta_{pp}^{n+1}||_0^2+ \left( 2\kappa  + \frac{C_T^2}{\overline{\epsilon}} \right) || \theta_{pp}^{n+1}||_1^2  \\
         &+ \left( \frac{2C_1^2}{ \Delta t \nu_p} + \frac{C_1^2}{ \nu_f} \right) || \theta_{g2}^{n+1}||_{0,\gamma}^2 + \left( \frac{C_1^2}{ \kappa} + 2 \overline{\epsilon} + \frac{1}{2\varepsilon_{21}} \right) ||\theta_{\lambda p}^{n+1}||_{1/2,\gamma}^2 + \frac{4C_K^2}{\nu_f}||\theta_{pf}^{n+1}||_0^2 \\
          & + \left( \frac{C_{KP}}{\nu_f} + \frac{C_{KP}}{\Delta t \nu_p}  + \frac{2}{\overline{\epsilon}} \right) ||\theta_{g1}^{n+1}||_{-1/2,\gamma}^2+  \nu_f || D(\bm{\chi}_{u,h}^{n+1})||_0^2 +  \frac{\Delta t \nu_p}{2}  ||D(\bm{\dot{\chi}}_{\eta,h}^{n+1})||_{0}^2 + \frac{\lambda  \Delta t}{4}  || \nabla \cdot \bm{\dot{\chi}}^{n+1}_{\eta,h} ||_0^2  \\
        & + \frac{\kappa}{2} || \nabla \chi_{pp,h}^{n+1}||_0^2  +  \frac{1}{2\beta }  ||\chi_{g2,h}^{n+1}||_{0,\gamma}^2 +  \frac{\overline{\epsilon}}{2} ||\chi_{\lambda p,h}^{n+1}||_{1/2,\gamma}^2   + \left( \frac{\varepsilon_{19}}{2}  + \frac{\varepsilon_{20}}{2} + \frac{\varepsilon_{21}}{2} \right) || \chi_{g1,h}^{n+1}||_{-1/2,\gamma}^2,
    \end{split}
\end{align}
where the variables $\varepsilon_{19},\varepsilon_{20}$, and $\varepsilon_{21}$ arising from Young's inequalities are yet to be defined.
\vspace{3mm}

\noindent \textbf{Step 2: Bound truncation error terms for dual variables using discrete inf-sup}\\
We pause to derive a bound for $||\chi_{pf,h}^{n+1}||_0^2 + || \chi_{g1,h}^{n+1}||_{-1/2,\gamma}^2$. Let $Z^h := Q_f^h \times \Lambda_{g1}^h$ and $M^h := U^h \times X^h \times Q_p^h \times \Lambda_{g2}^h \times \Lambda_\lambda^h$. The inf-sup condition proven in Theorem 3.1 of \cite{deCastro_2025_FPSIWP} states that for $(\chi_{pf,h}^{n+1}, \chi_{g1,h}^{n+1}) \in Z^h$, there exists a $\beta_2 > 0$ such that the following holds for $0 \neq \bm{m}_h := (\bm{v}_h,\bm{\varphi}_h,w_h,s_{2,h},\mu_h)$:
\begin{small}
\begin{align*}
   \beta_2 (||\chi_{pf,h}^{n+1}||_0^2 &+ || \chi_{g1,h}^{n+1}||_{-1/2,\gamma}^2)^{1/2} \leq  \underset{ \bm{m}_h\in M^h}{\text{sup}} \frac{-(\nabla \cdot \bm{v}_h, \chi_{pf,h}^{n+1}) - \langle \bm{v}_h \cdot \bm{n}_f,\chi_{g1,h}^{n+1}\rangle_{\gamma} - \langle \bm{\varphi}_h \cdot \bm{n}_p,\chi_{g1,h}^{n+1}\rangle_\gamma + \langle \mu_h,\chi_{g1,h}^{n+1}\rangle_\gamma }{||\bm{m}_h||_M}.
\end{align*}
\end{small}
To find an expression for the numerator, sum  \eqref{FPSI_FEM:eq:splitChiTheta_a}, \eqref{FPSI_FEM:eq:splitChiTheta_b}, and \eqref{FPSI_FEM:eq:splitChiTheta_e}:
\begin{align}\label{FPSI_FEM:eq:ErrorUseInfSup}
\begin{split}
   &\beta_2 (||\chi_{pf,h}^{n+1}||_0^2 + || \chi_{g1,h}^{n+1}||_{-1/2,\gamma}^2)^{1/2}\\
   & \leq  \underset{ \bm{m}_h\in M^h}{\text{sup}} \frac{1}{||\bm{m}_h||_M} \Big[  - \frac{\rho_f}{\Delta t} ( \bm{\chi}_{u,h}^{n+1} - \bm{\chi}_{u,h}^{n}, \bm{v}_h )_{\Omega_f} - 2  \nu_f ( D(\bm{\chi}_{u,h}^{n+1}), D(\bm{v}_h) )_{\Omega_f} +  ( \chi_{g2,h}^{n+1} \bm{\tau}_\gamma,\bm{v}_h)_\gamma  \\
   & + \frac{\rho_f}{\Delta t}(\bm{\theta}_u^{n+1}-\bm{\theta}_u^n,\bm{v}_h)_{\Omega_f} + 2  \nu_f \left( D(\bm{\theta}_{u}^{n+1}), D(\bm{v}_h) \right)_{\Omega_f} -  (\theta_{pf}^{n+1}, \nabla \cdot \bm{v}_h)_{\Omega_f} -  \langle \theta_{g1}^{n+1} \bm{n}_f,\bm{v}_h\rangle_\gamma \\
   &-  ( \theta_{g2}^{n+1} \bm{\tau}_\gamma,\bm{v}_h)_\gamma  -\frac{\rho_p}{\Delta t} ( \bm{\dot{\chi}}_{\eta,h}^{n+1} - \bm{\dot{\chi}}_{\eta,h}^n, \bm{\varphi}_h )_{\Omega_p} -  2  \nu_p ( D\left(\bm{\chi}_{\eta,h}^{n+1}\right), D(\bm{\varphi}_h) )_{\Omega_p} \\
   &- \lambda  ( \nabla \cdot  \bm{\chi}_{\eta,h}^{n+1}, \nabla \cdot \bm{\varphi}_h )_{\Omega_p} + \alpha (\chi_{pp,h}^{n+1}, \nabla \cdot \bm{\varphi}_h)_{\Omega_p} -  ( \chi_{g2,h}^{n+1}\bm{\tau}_\gamma, \bm{\varphi}_h)_\gamma  +  2  \nu_p \left( D\left(\bm{\theta}_{\eta}^{n+1}\right), D(\bm{\varphi}_h) \right)_{\Omega_p} \\
   & + \lambda  \left( \nabla \cdot  \bm{\theta}_{\eta}^{n+1}, \nabla \cdot \bm{\varphi}_h \right)_{\Omega_p} - \alpha (\theta_{pp}^{n+1}, \nabla \cdot \bm{\varphi}_h)_{\Omega_p} -\langle \theta_{g1}^{n+1}\bm{n}_p, \bm{\varphi}_h\rangle_\gamma +  ( \theta_{g2}^{n+1}\bm{\tau}_\gamma, \bm{\varphi}_h)_\gamma \\
     &   - (  \chi_{pp,h}^{n+1},\mu_h)_\gamma - \overline{\epsilon}( \chi_{\lambda p,h}^{n+1},\mu_h)_{1/2,\gamma} + \langle \theta_{g1}^{n+1},\mu_h\rangle_\gamma + (  \theta_{pp}^{n+1},\mu_h)_\gamma + \overline{\epsilon} (\theta_{\lambda p}^{n+1},\mu_h)_{1/2,\gamma} \Big]\\
     &:=  \underset{ \bm{m}_h\in M^h}{\text{sup}} \frac{1}{||\bm{m}_h||_M} \mathcal{J}_1.
     \end{split}
\end{align}

\noindent Using the inequalities $ ||\cdot||_{1/2,\gamma} \leq ||\cdot||_1$ and $||\cdot||_{0,\gamma} \leq ||\cdot||_{1/2,\gamma}$, $\mathcal{J}_1$ is bounded by
%||\cdot||_0 \leq ||\cdot||_1, ||D(\cdot)||_0 \leq ||\cdot||_1 obvious
\begin{align*}
    \mathcal{J}_1 &\leq ||\bm{v}_h||_1 \Big[   \frac{\rho_f}{\Delta t} || \bm{\chi}_{u,h}^{n+1} - \bm{\chi}_{u,h}^{n}||_0 + 2  \nu_f || D(\bm{\chi}_{u,h}^{n+1})||_0  +  C_T || \chi_{g2,h}^{n+1} ||_{0,\gamma} + \frac{\rho_f}{\Delta t} || \bm{\theta}_u^{n+1} - \bm{\theta}_u^{n+1}||_0 \\
    &+ 2  \nu_f || D(\bm{\theta}_{u}^{n+1})||_0+  \sqrt{2}C_K||\theta_{pf}^{n+1}||_0 +  || \theta_{g1}^{n+1}||_{-1/2,\gamma}  +  C_T|| \theta_{g2}^{n+1} ||_{0,\gamma} \Big]\\
    & + ||\bm{\varphi}_h||_1 \Big[ \frac{\rho_p}{\Delta t}||\bm{\dot{\chi}}_{\eta,h}^{n+1} - \bm{\dot{\chi}}_{\eta,h}^n||_0 +  2  \nu_p || D(\bm{\chi}_{\eta,h}^{n+1})||_0 + \lambda \sqrt{2}C_K || \nabla \cdot  \bm{\chi}_{\eta,h}^{n+1}||_0 \\
   &+ \alpha \sqrt{2}C_K||\chi_{pp,h}^{n+1}||_0 +  C_T|| \chi_{g2,h}^{n+1}||_{0,\gamma}  +  2  \nu_p || D(\bm{\theta}_{\eta}^{n+1})||_0 + \sqrt{2}C_K\lambda  ||\nabla \cdot  \bm{\theta}_{\eta}^{n+1}||_0 \\
   &+ \alpha \sqrt{2}C_K ||\theta_{pp}^{n+1}||_0 + ||\theta_{g1}^{n+1}||_{-1/2,\gamma}  +  C_T||\theta_{g2}^{n+1}||_{0,\gamma}    \Big]\\
   &+ ||\mu_h||_{1/2,\gamma} \Big[ || \chi_{pp,h}^{n+1}||_{0,\gamma} + \overline{\epsilon}|| \chi_{\lambda p,h}^{n+1}||_{1/2,\gamma} + ||\theta_{g1}^{n+1}||_{-1/2,\gamma} + || \theta_{pp}^{n+1}||_{0,\gamma} + \overline{\epsilon} || \theta_{\lambda p}^{n+1}||_{1/2,\gamma} \Big].
   \end{align*}

  \noindent For positive real numbers $\{a_i\}_{i=1}^n, \{b_i\}_{i=1}^n$, $a_1 b_1 + \ldots+ a_n b_n \leq (a_1 + \ldots + a_n)(b_1 + \ldots + b_n)$. This allows us to pull the test functions.
\begin{align*}
   &\mathcal{J}_1 \leq \Big( ||\bm{v}_h||_1 + ||\bm{\varphi}_h||_1 + ||\mu_h||_{1/2,\gamma} \Big) \Big(  \frac{\rho_f}{\Delta t} || \bm{\chi}_{u,h}^{n+1} - \bm{\chi}_{u,h}^{n}||_0 + 2  \nu_f || D(\bm{\chi}_{u,h}^{n+1})||_0   + \frac{\rho_p}{\Delta t}||\bm{\dot{\chi}}_{\eta,h}^{n+1} - \bm{\dot{\chi}}_{\eta,h}^n||_0 \\
   &+  2  \nu_p || D(\bm{\chi}_{\eta,h}^{n+1})||_0 + \lambda \sqrt{2}C_K || \nabla \cdot  \bm{\chi}_{\eta,h}^{n+1}||_0  + \left( \alpha \sqrt{2}C_K C_P +C_1 \right)||\nabla \chi_{pp,h}^{n+1}||_{0} +  2C_T || \chi_{g2,h}^{n+1} ||_{0,\gamma}\\
   &+ \overline{\epsilon}|| \chi_{\lambda p,h}^{n+1}||_{1/2,\gamma} + \frac{\rho_f}{\Delta t} ||\bm{\theta}_u^{n+1} - \bm{\theta}_u^{n}||_0  + 2  \nu_f || \bm{\theta}_{u}^{n+1}||_1  + \left( 2  \nu_p + d C_K^2 \lambda\right)  || \bm{\theta}_{\eta}^{n+1}||_1  + \alpha \sqrt{2}C_K ||\theta_{pp}^{n+1}||_0 \\
   &+ C_T || \theta_{pp}^{n+1}||_{1} + \overline{\epsilon} ||\theta_{\lambda p}^{n+1}||_{1/2,\gamma} +  2C_T|| \theta_{g2}^{n+1} ||_{0,\gamma} +  \sqrt{2}C_K||\theta_{pf}^{n+1}||_0 + 3 || \theta_{g1}^{n+1}||_{-1/2,\gamma}    \Big) \\
   &:= \Big( ||\bm{v}_h||_1 + ||\bm{\varphi}_h||_1 + ||\mu_h||_{1/2,\gamma} \Big) \mathcal{J}_2
  % &\leq \Big( ||\bm{v}_h||_1 + ||\bm{\varphi}_h||_1 + ||\mu_h||_{1/2,\gamma} + ||w_h||_1 + ||s_{2,h}||_{0,\gamma} \Big) \mathcal{J}_2\\
   %  &\leq \sqrt{5} \Big( ||\bm{v}_h||_1^2 + ||\bm{\varphi}_h||_1^2 + ||\mu_h||^2_{1/2,\gamma} + ||w_h||^2_1 + ||s_{2,h}||^2_{0,\gamma} \Big)^{1/2} \mathcal{J}_2
   \leq \sqrt{5} || \bm{m}_h ||_M \mathcal{J}_2.
\end{align*}

\noindent Combining with \eqref{FPSI_FEM:eq:ErrorUseInfSup} yields
\begin{align}\label{FPSI_FEM:eq:BoundZterms}
\begin{split}
     & ||\chi_{pf,h}^{n+1}||_0^2 + || \chi_{g1,h}^{n+1}||_{-1/2,\gamma}^2 \leq \frac{1}{\beta_2^2}  \left( \underset{ \bm{m}_h\in M^h}{\text{sup}} \frac{ \mathcal{J}_1 }{||\bm{m}_h||_M} \right)^2 \leq \frac{1}{\beta_2^2} \left(  \underset{ \bm{m}_h\in M^h}{\text{sup}} \frac{\sqrt{5} || \bm{m}_h ||_M \mathcal{J}_2}{||\bm{m}_h||_M}   \right)^2 = \frac{5}{\beta_2^2} \mathcal{J}_2^2 .
\end{split}
\end{align} 
Recalling that $||\bm{\theta}_u^{n+1} - \bm{\theta}_u^n||_0^2 = \Delta t^2 || \bm{\dot{\theta}}_u^{n+1}||_0^2$, the term $\mathcal{J}_2^2$ is bounded by
\begin{align}\label{FPSI_FEM:eq:boundJ2_squaredError}
    \begin{split}
         \mathcal{J}_2^2 &\leq 17\Big(  \frac{\rho_f^2}{\Delta t^2} || \bm{\chi}_{u,h}^{n+1} - \bm{\chi}_{u,h}^{n}||^2_0 + 4  \nu_f^2 || D(\bm{\chi}_{u,h}^{n+1})||^2_0 + \frac{\rho_p^2}{\Delta t^2}||\bm{\dot{\chi}}_{\eta,h}^{n+1} - \bm{\dot{\chi}}_{\eta,h}^n||^2_0 +  4  \nu_p^2 || D(\bm{\chi}_{\eta,h}^{n+1})||^2_0 \\
         &+ 2\lambda^2 C_K^2 || \nabla \cdot  \bm{\chi}_{\eta,h}^{n+1}||^2_0 + \left( \alpha \sqrt{2} C_K C_P + C_1 \right)^2 ||\nabla \chi_{pp,h}^{n+1}||^2_0   +  4C_T^2 || \chi_{g2,h}^{n+1} ||^2_{0,\gamma} + \overline{\epsilon}^2|| \chi_{\lambda p,h}^{n+1}||^2_{1/2,\gamma} \\
   & + \rho_f^2 || \bm{\dot{\theta}}_u^{n+1}||_0^2 + 4  \nu_f^2 || \bm{\theta}_{u}^{n+1}||^2_1 +  (2  \nu_p + 2 C_K^2 \lambda )^2 ||\bm{\theta}_{\eta}^{n+1}||_1^2 + 2\alpha^2  C_K^2 ||\theta_{pp}^{n+1}||^2_0  \\
   &  + C_T^2|| \theta_{pp}^{n+1}||^2_{1}  + \overline{\epsilon}^2 ||\theta_{\lambda p}^{n+1}||_{1/2,\gamma}^2 +  4 C_T^2|| \theta_{g2}^{n+1} ||^2_{0,\gamma} +  2 C_K^2||\theta_{pf}^{n+1}||^2_0 +  9|| \theta_{g1}^{n+1}||_{-1/2,\gamma}^2 \Big). 
    \end{split}
\end{align}

\noindent \textbf{Step 3: Combine bounds; move all truncation terms to LHS and simplify}\\
Now, add the expressions in \eqref{FPSI_FEM:eq:ErrorInnProducts} and \eqref{FPSI_FEM:eq:BoundZterms}, using the bounds for $\sum_{i=1}^{21} R_j$ and $\mathcal{J}_2$ given in \eqref{FPSI_FEM:eq:Err_RHSbounds} and  \eqref{FPSI_FEM:eq:boundJ2_squaredError}. We also multiply \eqref{FPSI_FEM:eq:BoundZterms} by a constant $\delta > 0$, to be defined later.

\begin{align*}%\label{FPSI_FEM:Eq:combinewithInfSup}
      \begin{split}
       &\frac{\rho_f}{2 \Delta t} \left( || \bm{\chi}_{u,h}^{n+1} ||_0^2 - ||\bm{\chi}_{u,h}^{n} ||_0^2 + || \bm{\chi}_{u,h}^{n+1} - \bm{\chi}_{u,h}^{n}||_0^2 \right) + \frac{\rho_p}{2 \Delta t} \left( || \bm{\dot{\chi}}_{\eta,h}^{n+1} ||_0^2 - ||\bm{\dot{\chi}}_{\eta,h}^n||_0^2 + || \bm{\dot{\chi}}_{\eta,h}^{n+1} - \bm{\dot{\chi}}_{\eta,h}^{n} ||_0^2 \right)\\
       & \hspace{3mm} + \frac{s_0}{2 \Delta t} \left( ||\chi_{pp,h}^{n+1}||_0^2 - ||\chi_{pp,h}^n||_0^2 + ||\chi_{pp,h}^{n+1}-\chi_{pp,h}^n||_0^2 \right) + \frac{1}{\beta } || \chi_{g2,h}^{n+1}||_{0,\gamma}^2  + \overline{\epsilon}|| \chi_{\lambda p,h}^{n+1}||^2_{1/2,\gamma} + 2  \nu_f || D(\bm{\chi}_{u,h}^{n+1}) ||_0^2 \\
       &\hspace{3mm} + \kappa  ||\nabla \chi_{pp,h}^{n+1}||_0^2 +   \frac{\nu_p}{\Delta t} \left(  ||  D(\bm{\chi}_{\eta,h}^{n+1})||_0^2 - || D(\bm{\chi}_{\eta,h}^{n})||_0^2 \right)  + \frac{\lambda}{2\Delta t} \left( ||\nabla \cdot \bm{\chi}_{\eta,h}^{n+1}||_0^2 - ||\nabla \cdot \bm{\chi}_{\eta,h}^{n}||_0^2  \right) \\
       &\hspace{3mm} + \Delta t \nu_p || D(\bm{\dot{\chi}}_{\eta,h}^{n+1})||_0^2 + \frac{\lambda \Delta t }{2} || \nabla \cdot \bm{\dot{\chi}}_{\eta,h}^{n+1} ||_0^2 + \delta ||\chi_{pf,h}^{n+1}||_0^2 + \delta || \chi_{g1,h}^{n+1}||_{-1/2,\gamma}^2    \\
       &\leq \sum_{j=1}^{21} R_j + \frac{5\delta}{\beta_2^2} \mathcal{J}^2_2\\
       &\leq   \left( 8\nu_f + C_T^2 \beta + \frac{1}{2\varepsilon_{19}} + \frac{85\delta}{\beta_2^2}4\nu_f^2 \right) || \bm{\theta}_{u}^{n+1}||_1^2 +  \left( \frac{85\delta \rho_f^2}{\beta_2^2 }  + \frac{\rho_f^2 C_P^2 C_K^2}{\nu_f} \right) ||\bm{\dot{\theta}}_u^{n+1}||_0^2   \\
         &+ \left( \frac{8\nu_p}{\Delta t} + \frac{4\lambda C_K^2}{\Delta t}+    \frac{85\delta}{\beta_2^2}(2  \nu_p + 2 C_K^2 \lambda )^2  \right) || \bm{\theta}_{\eta}^{n+1} ||_1^2 + \left( \frac{4\alpha^2 C_P^2  C_K^2 }{\kappa}+ C_T^2 \beta + \frac{1}{2\varepsilon_{20}} \right) ||\bm{\dot{\theta}}_{\eta}^{n+1}||_1^2 \\
         &+ \left( \frac{2 \alpha^2}{\lambda \Delta t} + \frac{85\delta}{\beta_2^2}2\alpha^2  C_K^2 \right) ||\theta_{pp}^{n+1}||_0^2 + \left( 2\kappa    + \frac{C_T^2}{\overline{\epsilon}}  + \frac{85\delta}{\beta_2^2}C_T^2 \right) ||\theta_{pp}^{n+1}||_1^2\\
         &+ \left( \frac{2C_1^2}{ \Delta t \nu_p} + \frac{C_1^2}{ \nu_f} + \frac{85\delta}{\beta_2^2}4C_T^2 \right) || \theta_{g2}^{n+1}||_{0,\gamma}^2 + \left( \frac{C_1^2}{ \kappa} +2\overline{\epsilon}  + \frac{1}{2\varepsilon_{21}} + \frac{85\delta}{\beta_2^2} \overline{\epsilon}^2 \right) ||\theta_{\lambda p}^{n+1}||_{1/2,\gamma}^2 \\
          & + \left( \frac{4C_K^2}{\nu_f}  +  \frac{85\delta}{\beta_2^2} 2 C_K^2 \right) ||\theta_{pf}^{n+1}||^2_0 + \left( \frac{C_{KP}}{\nu_f} + \frac{C_{KP}}{\Delta t \nu_p}  + \frac{2}{\overline{\epsilon}} + \frac{85\delta}{\beta_2^2} 9 \right) ||\theta_{g1}^{n+1}||_{-1/2,\gamma}^2 \\
        &+ \left(  \nu_f + \frac{85\delta}{\beta_2^2}4\nu_f^2 \right) || D(\bm{\chi}_{u,h}^{n+1})||_0^2 + \frac{85\delta \rho_f^2}{\beta_2^2 \Delta t^2}  || \bm{\chi}_{u,h}^{n+1} - \bm{\chi}_{u,h}^{n}||^2_0  + \frac{\Delta t \nu_p}{2} ||D(\bm{\dot{\chi}}_{\eta,h}^{n+1})||_{0}^2  \\
        & + \frac{\lambda  \Delta t}{4}  || \nabla \cdot \bm{\dot{\chi}}^{n+1}_{\eta,h} ||_0^2 + \frac{85\delta \rho_p^2}{\beta_2^2 \Delta t^2}  ||\bm{\dot{\chi}}_{\eta,h}^{n+1} - \bm{\dot{\chi}}_{\eta,h}^n||^2_0 + \frac{85\delta}{\beta_2^2} 4\nu_p^2  || D(\bm{\chi}_{\eta,h}^{n+1})||^2_0 +  \frac{85\delta}{\beta_2^2} 2\lambda^2 C_K^2 || \nabla \cdot  \bm{\chi}_{\eta,h}^{n+1}||^2_0 \\
        & + \left(  \frac{\kappa}{2} + \frac{85\delta}{\beta_2^2} (\alpha \sqrt{2} C_K C_P + C_1 )^2 \right) || \nabla \chi_{pp,h}^{n+1}||_0^2  + \left(  \frac{1}{2\beta} + \frac{85\delta}{\beta_2^2}4C_T^2 \right) ||\chi_{g2,h}^{n+1}||_{0,\gamma}^2\\
        &+ \left(  \frac{\overline{\epsilon}}{2}  + \frac{85\delta}{\beta_2^2} \overline{\epsilon}^2 \right) ||\chi_{\lambda p,h}^{n+1}||_{1/2,\gamma}^2   + \left( \frac{\varepsilon_{19}}{2}  + \frac{\varepsilon_{20}}{2} + \frac{\varepsilon_{21}}{2} \right) || \chi_{g1,h}^{n+1}||_{-1/2,\gamma}^2 . 
    \end{split}
\end{align*}

\noindent 
 In the next step, all the terms $\chi_{*,h}^{n+1}$ will be moved to the left side except 
the $|| D(\bm{\chi}_{\eta,h}^{n+1})||^2_0$ and  $|| \nabla \cdot  \bm{\chi}_{\eta,h}^{n+1}||^2_0$  terms.   
Let $K_i$, $i=1,2,\ldots,6$,  represent the constants that will appear on the left hand side of the inequality: 
\begin{align*}
K_1 &:= \frac{\rho_f}{2\Delta t} - \frac{85\delta \rho_f^2}{\beta_2^2 \Delta t^2} > 0, \quad K_2 := \frac{\rho_p}{2\Delta t} - \frac{85\delta \rho_p^2}{\beta_2^2 \Delta t^2} > 0, \quad K_3 := \frac{1}{2\beta} - \frac{340\delta C_T^2}{\beta_2^2} > 0,\\
K_4 &:= \frac{\overline{\epsilon}}{2} - \frac{85\delta \overline{\epsilon}^2}{\beta_2^2} > 0, \quad K_5 := \nu_f - \frac{340\delta \nu_f^2}{\beta_2^2} > 0, \quad K_6 := \frac{\kappa}{2} - \frac{85\delta}{\beta_2^2}(\alpha \sqrt{2} C_K C_P+C_1)^2 > 0\\
K_7 &:= \max \Big\{ \frac{170\nu_p}{\beta_2^2}, \frac{170 C_K^2 \lambda}{\beta_2^2} \Big\}.
\end{align*}

\noindent In order to ensure positivity of these constants, we set the parameter $\delta$ as 
\begin{align*}
\delta &:= \min\Big\{ \frac{\beta_2^2 \Delta t}{170 \rho_f}, \frac{\beta_2^2 \Delta t}{170 \rho_p}, \frac{\beta_2^2}{680 C_T^2 \beta}, \frac{\beta_2^2}{170 \overline{\epsilon}},\frac{\beta_2^2}{340 \nu_f},  \frac{\beta_2^2 \kappa}{170(\alpha \sqrt{2} C_K C_P + \sqrt{C_1})^2}, \frac{1}{2\Delta t K_7} \Big\},
\end{align*} and choose 
$\varepsilon_{19} = \delta / 2, \ \varepsilon_{20}= \varepsilon_{21} = \delta / 4. $
\begin{comment}
Also, to simplify constants, define
\begin{align*}
    \vartheta_1 &:=   8\nu_f + C_T^2 \beta + \frac{1}{\delta} + \frac{300\delta}{\beta_2^2}\nu_f^2 , \qquad \vartheta_2 :=  \frac{\rho_f^2 C_P^2 C_K^2}{ \nu_f}  \\
    \vartheta_3 &:= \frac{8\nu_p}{\Delta t} + \frac{2\lambda dC_K^2}{\Delta t}+    \frac{75\delta}{\beta_2^2}(2  \nu_p + d C_K\lambda )^2  , \qquad \vartheta_4 := \frac{2\alpha^2 C_P^2 d C_K^2 }{\kappa}+ C_T^2 \beta + \frac{2}{\delta}  \\
    \vartheta_5 &:= \frac{2 \alpha^2}{\lambda \Delta t} + \frac{75\delta}{\beta_2^2}\alpha^2 d C_K^2 , \qquad \vartheta_6 :=  2\kappa    + \frac{C_T^2}{\overline{\epsilon}}  + \frac{75\delta}{\beta_2^2}C_T^2 , \qquad \vartheta_7 := \frac{2C_1^2}{ \Delta t \nu_p} + \frac{C_1^2}{ \nu_f} + \frac{300\delta}{\beta_2^2}C_T^2 \\
    \vartheta_8 &:=  \frac{C_1^2}{ \kappa}   + 2\overline{\epsilon} + \frac{2}{\delta} , \qquad \vartheta_9 := \frac{2dC_K^2}{\nu_f}  +  \frac{75\delta}{\beta_2^2}  d C_K^2 , \qquad \vartheta_{10} :=   \frac{C_{KP}}{\nu_f} + \frac{C_{KP}}{\Delta t \nu_p}  + \frac{2}{\overline{\epsilon}} + \frac{675\delta}{\beta_2^2} 
\end{align*}
\end{comment}
By the definition of $K_7$, note that 
 $$\frac{340\delta}{\beta_2^2} \nu_p^2  || D(\bm{\chi}_{\eta,h}^{n+1})||^2_0 +  \frac{85\delta}{\beta_2^2} \lambda^2 d C_K^2 || \nabla \cdot  \bm{\chi}_{\eta,h}^{n+1}||^2_0 \leq \delta K_7 || \bm{\chi}_{\eta,h}^{n+1}||_E^2.$$
To simplify notation for the constants multiplying the approximation error terms, define $C_{\Delta t}$ as the maximum of each of the ten constants multiplying a norm of $|| \theta_*^{n+1}||$. The subscript $\Delta t$ here denotes the dependence of $C_{\Delta t}$ on the time step; in fact, it is inversely proportional to $\Delta t$. Note that for our spatial convergence analysis, $\Delta t$ is considered fixed. $C_{\Delta t}$ depends upon problem parameters, the discrete inf-sup parameter $\beta_2$, $\delta$, and constants from the use of trace, Korn, and Poincar\'{e} inequalities; it does not depend on mesh size.
%The reason the $\Delta t^{-1}$ factors appear is due to the presence of the $\Delta t$ multiplying the terms $||D(\bm{\dot{\chi}}_{\eta,h}^{n+1})||_0^2, || \nabla \cdot \bm{\dot{\chi}}_{\eta,h}^{n+1}||_0^2$ in the LHS of \eqref{FPSI_FEM:Eq:ErrorLRSums}. We had to introduce them in the Young's inequalities for $R_6$ through $R_{10}$ in order to tuck those terms onto the left hand side. 
Moving the truncation error terms to the left and simplifying yields
\begin{align}\label{FPSI_FEM:Eq:combinewithInfSup}
      \begin{split}
       &\frac{\rho_f}{2 \Delta t} \left( || \bm{\chi}_{u,h}^{n+1} ||_0^2 - ||\bm{\chi}_{u,h}^{n} ||_0^2 \right) + K_1 || \bm{\chi}_{u,h}^{n+1} - \bm{\chi}_{u,h}^{n}||_0^2  + \frac{\rho_p}{2 \Delta t} \left( || \bm{\dot{\chi}}_{\eta,h}^{n+1} ||_0^2 - ||\bm{\dot{\chi}}_{\eta,h}^n||_0^2 \right) + K_2|| \bm{\dot{\chi}}_{\eta,h}^{n+1} - \bm{\dot{\chi}}_{\eta,h}^{n} ||_0^2 \\
       & + \frac{s_0}{2 \Delta t} \left( ||\chi_{pp,h}^{n+1}||_0^2 - ||\chi_{pp,h}^n||_0^2 + ||\chi_{pp,h}^{n+1}-\chi_{pp,h}^n||_0^2 \right) + K_3|| \chi_{g2,h}^{n+1}||_{0,\gamma}^2  + K_4|| \chi_{\lambda p,h}^{n+1}||^2_{1/2,\gamma} + K_5 || D(\bm{\chi}_{u,h}^{n+1}) ||_0^2 \\
       & + K_6  ||\nabla \chi_{pp,h}^{n+1}||_0^2 +   \frac{\nu_p}{\Delta t} \left(  ||  D(\bm{\chi}_{\eta,h}^{n+1})||_0^2 - || D(\bm{\chi}_{\eta,h}^{n})||_0^2 \right)  + \frac{\lambda}{2\Delta t} \left( ||\nabla \cdot \bm{\chi}_{\eta,h}^{n+1}||_0^2 - ||\nabla \cdot \bm{\chi}_{\eta,h}^{n}||_0^2  \right) \\
       &+ \frac{\Delta t \nu_p}{2} || D(\bm{\dot{\chi}}_{\eta,h}^{n+1})||_0^2 + \frac{\lambda \Delta t }{4} || \nabla \cdot \bm{\dot{\chi}}_{\eta,h}^{n+1} ||_0^2 + \delta ||\chi_{pf,h}^{n+1}||_0^2 + \frac{\delta}{2} || \chi_{g1,h}^{n+1}||_{-1/2,\gamma}^2    \\
       &\leq C_{\Delta t} \Big(  || \bm{\theta}_{u}^{n+1}||_1^2 +  ||\bm{\dot{\theta}}_u^{n+1}||_0^2  +  || \bm{\theta}_{\eta}^{n+1} ||_1^2 +  ||\bm{\dot{\theta}}_{\eta}^{n+1}||_1^2 +  ||\theta_{pp}^{n+1}||_0^2 + ||\theta_{pp}^{n+1}||_1^2+  || \theta_{g2}^{n+1}||_{0,\gamma}^2 \\
         &+ ||\theta_{\lambda p}^{n+1}||_{1/2,\gamma}^2  +  ||\theta_{pf}^{n+1}||^2_0 +||\theta_{g1}^{n+1}||_{-1/2,\gamma}^2  \Big) + \frac{340\delta}{\beta_2^2} \nu_p^2  || D(\bm{\chi}_{\eta,h}^{n+1})||^2_0 +  \frac{85\delta}{\beta_2^2} \lambda^2 dC_K^2 || \nabla \cdot  \bm{\chi}_{\eta,h}^{n+1}||^2_0\\
        &\leq C_{\Delta t}\mathcal{J}_\theta^{n+1} + \delta K_7 ||\bm{\chi}_{\eta,h}^{n+1}||_E^2,
    \end{split}
\end{align}
where all terms involving approximation errors are identified as $\mathcal{J}^{n+1}_\theta$. 
Now, dropping the positive terms involving $ || \bm{\chi}_{u,h}^{n+1} - \bm{\chi}_{u,h}^n||_0^2$, $ || \bm{\dot{\chi}}_{\eta,h}^{n+1} - \bm{\dot{\chi}}_{\eta,h}^n||_0^2$, $||\nabla \cdot \bm{\dot{\chi}}_{\eta,h}^{n+1}||_0^2$, and $||\chi_{pp,h}^{n+1}-\chi_{pp,h}^n||_0^2$ from the left, we multiply by $2\Delta t$ and sum from $n=0$ to $n=M-1$, where $M=1,2,\ldots,N$. 
For simplicity, take $\bm{u}_h^0 := \mathcal{I}^{U^h}(\bm{u}^0)$, $\bm{\eta}_h^0 := \mathcal{P}^{X^h}(\bm{\eta}^0)$, $\bm{\dot{\eta}}_h^0 := \mathcal{P}^{X^h}(\bm{\dot{\eta}}^0)$, and $p_{p,h}^0 := \mathcal{P}^{Q_p^h}(p_p^0)$ so that the terms involving initial conditions become zero.
\begin{align*}%\label{FPSI_FEM:Eq:combinewithInfSup}
      \begin{split}
       &  || \bm{\chi}_{\eta,h}^{M}||_E^2 + \Delta t \sum_{n=0}^{M-1} \frac{1}{M\Delta t}\Big[ \rho_f  || \bm{\chi}_{u,h}^{M} ||_0^2   + \rho_p || \bm{\dot{\chi}}_{\eta,h}^{M} ||_0^2  + s_0 ||\chi_{pp,h}^{M}||_0^2 \Big]     \\
       &+ \Delta t \sum_{n=0}^{M-1} \Big[ 2K_3|| \chi_{g2,h}^{n+1}||_{0,\gamma}^2  + 2K_4|| \chi_{\lambda p,h}^{n+1}||^2_{1/2,\gamma}  + 2K_5  || D(\bm{\chi}_{u,h}^{n+1}) ||_0^2  + 2 K_6  ||\nabla \chi_{pp,h}^{n+1}||_0^2 \\
       & \hspace{3mm}  + \Delta t \nu_p || D(\bm{\dot{\chi}}_{\eta,h}^{n+1})||_0^2 + 2 \delta  ||\chi_{pf,h}^{n+1}||_0^2 + \delta || \chi_{g1,h}^{n+1}||_{-1/2,\gamma}^2    \Big] \\
       &\leq   \Delta t \sum_{n=0}^{M-1} 2\Big[ C_{\Delta t} \mathcal{J}_\theta^{n+1} \Big] +  \Delta t \sum_{n=0}^{M-2} \Big[ 2 \delta K_7  || \bm{\chi}_{\eta,h}^{n+1}||^2_E  \Big] + 2\Delta t \delta K_7 ||\bm{\chi}_{\eta,h}^{M}||_E^2. 
    \end{split}
\end{align*}
\textbf{Step 4: Apply discrete Gronwall's lemma}\\
To handle the sum of the terms $||\bm{\chi}_{\eta,h}^{n+1}||_E^2$ on the right hand side of the inequality, we apply the discrete Gronwall's lemma, restated in Lemma \ref{FPSI_FEM:lemma:DiscGronwall} for clarity (see, for example, \cite{Ambartsumyan_2018}). 
\begin{lemma}[Discrete Gronwall's lemma]\label{FPSI_FEM:lemma:DiscGronwall}
    Let $\Delta t > 0$, $B \geq 0$, and let $a^n, b^n, c^n, d^n$ for $n\geq 0$ be non-negative sequences such that  and $$ a^N + \Delta t \sum_{n=0}^{N-1} b^{n+1} \leq \Delta t \sum_{n=0}^{N-2} d^{n+1} a^{n+1} \ + \Delta t \sum_{n=0}^{N-1} c^{n+1} \ + B.$$
    Then \begin{equation*}a^N +  \Delta t \sum_{n=0}^{N-1} b^{n+1} \leq \text{\normalfont exp} \left( \Delta t \sum_{n=0}^{N-2} d^{n+1}\right)\left( \Delta t \sum_{n=0}^{N-1} c^{n+1} \ + B\right).\end{equation*}
\end{lemma}

Identifying $a^{n+1}$ with $|| \bm{\chi}_{\eta,h}^{n+1}||_E^2$, define $K_8 := 1 - 2\Delta t \delta K_7$, which is positive by the definition of $\delta$, and $C^* := \exp\left(\frac{2\delta K_7 \Delta t (M-1)}{K_8}\right).$ Take
    \begin{align*}
    b^{n+1} &:= \frac{1}{M\Delta t}\Big[ \rho_f  || \bm{\chi}_{u,h}^{M} ||_0^2   + \rho_p || \bm{\dot{\chi}}_{\eta,h}^{M} ||_0^2  + s_0 ||\chi_{pp,h}^{M}||_0^2 \Big]  + \Big[ 2K_3|| \chi_{g2,h}^{n+1}||_{0,\gamma}^2  + 2K_4|| \chi_{\lambda p,h}^{n+1}||^2_{1/2,\gamma}  \\
    &\hspace{3mm}+ 2K_5  || D(\bm{\chi}_{u,h}^{n+1}) ||_0^2  + 2 K_6  ||\nabla \chi_{pp,h}^{n+1}||_0^2  + \Delta t \nu_p || D(\bm{\dot{\chi}}_{\eta,h}^{n+1})||_0^2  + 2 \delta  ||\chi_{pf,h}^{n+1}||_0^2 + \delta || \chi_{g1,h}^{n+1}||_{-1/2,\gamma}^2    \Big]\\
    d^{n+1} &:= 2 \delta K_7 \\
    c^{n+1} &:= 2 C_{\Delta t} \mathcal{J}_\theta^{n+1}.
\end{align*}
Then applying Gronwall's lemma yields  $$K_8 a^M + \Delta t \sum_{n=0}^{M-1} b^{n+1} \leq \exp \left(\Delta t \sum_{n=0}^{M-2} \frac{1}{K_8} d^{n+1} \right)\left( \Delta t \sum_{n=0}^{M-1} c^{n+1} \right).$$ 
Bound the remaining the $||D(\cdot)||$ and $||\cdot||_E$ on the left side by $|| \cdot||_1$ to obtain:
\begin{align*}
\begin{split}
    &\frac{K_8}{C_\eta} || \bm{\chi}_{\eta,h}^{M}||_1^2 +  \rho_f  || \bm{\chi}_{u,h}^{M} ||_0^2   + \rho_p || \bm{\dot{\chi}}_{\eta,h}^{M} ||_0^2  + s_0 ||\chi_{pp,h}^{M}||_0^2 + \Delta t \sum_{n=0}^{M-1}  \Big[ 2K_3|| \chi_{g2,h}^{n+1}||_{0,\gamma}^2  + 2K_4|| \chi_{\lambda p,h}^{n+1}||^2_{1/2,\gamma} \\
    &\hspace{3mm} + \frac{2K_5}{C_{KP}}  || \bm{\chi}_{u,h}^{n+1} ||_1^2  + \frac{2 K_6}{C_{KP}}  || \chi_{pp,h}^{n+1}||_1^2   + \frac{\Delta t \nu_p}{C_{KP}} || \bm{\dot{\chi}}_{\eta,h}^{n+1}||_1^2  + 2 \delta  ||\chi_{pf,h}^{n+1}||_0^2 + \delta || \chi_{g1,h}^{n+1}||_{-1/2,\gamma}^2    \Big] \\
    &\leq C^* \Bigg( \Delta t \sum_{n=0}^{M-1}2 C_{\Delta t} \mathcal{J}_\theta^{n+1} \Bigg).
    \end{split}
\end{align*}

\noindent Take $K_9 := \min\big\{ \frac{K_8}{C_\eta}, \rho_f, \rho_p, s_0, 2K_3, 2K_4, \frac{2K_5}{C_{KP}}, \frac{2K_6}{C_{KP}}, \frac{\Delta t \nu_p}{C_{KP}}, \delta   \big\}$. Then
\begin{align}\label{FPSI_FEM:eq:ResultAfterGronwallError}
\begin{split}
    &|| \bm{\chi}_{\eta,h}^{M}||_1^2 +   || \bm{\chi}_{u,h}^{M} ||_0^2   +  || \bm{\dot{\chi}}_{\eta,h}^{M} ||_0^2  +  ||\chi_{pp,h}^{M}||_0^2 + \Delta t \sum_{n=0}^{M-1}  \Big[ || \chi_{g2,h}^{n+1}||_{0,\gamma}^2  + || \chi_{\lambda p,h}^{n+1}||^2_{1/2,\gamma}  +  || \bm{\chi}_{u,h}^{n+1} ||_1^2 \\
    &\hspace{3mm} +   || \chi_{pp,h}^{n+1}||_1^2   + || \bm{\dot{\chi}}_{\eta,h}^{n+1}||_1^2  +  ||\chi_{pf,h}^{n+1}||_0^2 +  || \chi_{g1,h}^{n+1}||_{-1/2,\gamma}^2    \Big] \leq \frac{C^*}{K_9} \Bigg( \Delta t \sum_{n=0}^{M-1}2 C_{\Delta t}\mathcal{J}_\theta^{n+1} \Bigg).
    \end{split}
\end{align}

%\begin{align}\label{FPSI_FEM:eq:Dividek9}
%\begin{split}
 %   & || \bm{\chi}_{\eta,h}||_{\ell^\infty(0,T;X)}^2 +   || \bm{\chi}_{u,h} ||_{\ell^\infty(0,T;\bm{L^2}(\Omega_f))}^2   + || \bm{\dot{\chi}}_{\eta,h} ||_{\ell^\infty(0,T;\bm{L^2}(\Omega_p))}^2  +  ||\chi_{pp,h}||_{\ell^\infty(0,T;L^2(\Omega_p))}^2 \\ %&+  || \bm{\chi}_{u,h} ||_{\ell^2(0,T;U)}^2   +  || \bm{\dot{\chi}}_{\eta,h}||_{\ell^2(0,T;X)}^2 
 %   +  ||\chi_{pp,h}||_{\ell^2(0,T;Q_p)}^2 +  || \chi_{g2,h}||_{\ell^2(0,T;\Lambda_{g2})}^2   \\
%     &+ || \chi_{\lambda p,h}||^2_{\ell^2(0,T;\Lambda_\lambda)} +   ||\chi_{pf,h}||_{\ell^2(0,T;Q_f)}^2 +  || \chi_{g1,h}||_{\ell^2(0,T;\Lambda_{g1})}^2    \\
%    &\leq \frac{2C^*}{K_9} \Bigg(  2T \vartheta_2||\bm{\dot{\theta}}_u||_{\ell^\infty(0,T;\bm{L^2}(\Omega_f))}^2 +  2T \vartheta_5 ||\theta_{pp}||_{\ell^\infty(0,T;L^2(\Omega_p))}^2 + 2\vartheta_3 || \bm{\theta}_{\eta}^{n+1} ||_{\ell^2(0,T;X)}^2  +2 \vartheta_1 || \bm{\theta}_{u}||_{\ell^2(0,T;U)}^2 \\
%    &+ 2 \vartheta_4 ||\bm{\dot{\theta}}_{\eta}||_{\ell^2(0,T;X)}^2  +  2\vartheta_6 || \theta_{pp}||_{\ell^2(0,T;Q_p)}^2  + 2\vartheta_7|| \theta_{g2}||_{\ell^2(0,T;\Lambda_{g2})}^2 + 2\vartheta_8 ||\theta_{\lambda p}||_{\ell^2(0,T;\Lambda_\lambda)}^2  \\
%    &+ 2\vartheta_9 ||\theta_{pf}||^2_{\ell^2(0,T;Q_f)} + 2\vartheta_{10} ||\theta_{g1}||_{\ell^2(0,T;\Lambda_{g1})}^2 \Bigg) \\
%          &:= \frac{2C^*}{K_9}\mathcal{RH}
 %         \end{split}
%\end{align}
\noindent \textbf{Step 5: Apply the approximation properties to the right hand side}\\
To provide an upper bound for $\mathcal{J}_\theta$, we take its infimum over $s_{1,h} \in \Lambda_{g1}^{h}$ and over $\mu_h \in \Lambda_\lambda^h$, applying the approximation properties in \eqref{FPSI_FEM:eq:Approx_StokesProjection}-\eqref{FPSI_FEM:eq:H1/2approxProp}. A constant $C$ without subscript represents a generic $C$ from the approximation properties which is not dependent on mesh size. 
\begin{align}\label{FPSI_FEM:eq:BoundJtheta}
\begin{split}
  &\underset{\mu_h \in\Lambda_{\lambda }^h}{\inf}     \underset{s_{1,h}\in\Lambda_{g1}^h}{\inf}  \mathcal{J}_\theta := || \bm{\theta}_{u}^{n+1}||_1^2 + ||\bm{\dot{\theta}}_u^{n+1}||_0^2  +  || \bm{\theta}_{\eta}^{n+1} ||_1^2 + ||\bm{\dot{\theta}}_\eta^{n+1}||_1^2 + ||\theta_{pp}^{n+1}||_0^2 + ||\theta_{pp}^{n+1}||_1^2 \\
         &\hspace{5mm} +  || \theta_{g2}^{n+1}||_{0,\gamma}^2 + \underset{\mu_h \in \Lambda_\lambda^h}{\inf}  ||\theta_{\lambda p}^{n+1}||_{1/2,\gamma}^2  + ||\theta_{pf}^{n+1}||^2_0 +  \underset{s_{1,h}\in\Lambda_{g1}^h}{\inf} ||\theta_{g1}^{n+1}||_{-1/2,\gamma}^2\\
         &\leq  Ch_1^{2k_u}|| \bm{u}^{n+1}||_{k_u+1}^2 + C h_1^{2k_u}||\bm{\dot{u}}^{n+1}||_{k_u+1}^2  +  C h_2^{2k_\eta}|| \bm{\eta}^{n+1} ||_{k_\eta+1}^2 +  C h_2^{2k_\eta} ||\bm{\dot{\eta}}^{n+1}||_{k_\eta+1}^2 \\
         &\hspace{5mm} + C h_2^{2(k_{pp}+1)}||p_p^{n+1}||_{k_{pp}+1}^2 +  C h_2^{2k_{pp}}||p_p^{n+1}||_{k_{pp}+1}^2+  C h_\gamma^{2(k_{g2}+1)}|| g_2^{n+1}||_{k_{g2}+1,\gamma}^2 \\
         &\hspace{5mm} + Ch_\gamma^{2k_{\lambda}} \underset{\substack{w \in H^{k_{\lambda}+1}(\Omega_r) \\ w|_\gamma = \lambda_p}}{\inf}||w||^2_{k_{\lambda+1},\Omega_r} + C h_1^{2(k_{pf}+1)}||p_f^{n+1}||^2_{k_{pf}+1} +  Ch_\gamma^{2(k_{g1}+1)} \underset{\substack{w \in H^{k_{g1}+1}(\Omega_r) \\ w|_\gamma = g_1}}{\inf}||w||^2_{k_{g1}+1,\Omega_r}.  %+Ch_\gamma^2 ||g_1^{n+1}||_{1/2,\gamma}^2.
         \end{split}
\end{align}

%Note that $||f||^2_{\ell^\infty(0,T,X)} \leq ||f||^2_{\ell^2(0,T,X)}$. 
\noindent To conclude, note that the triangle inequality allows us to write the norms of an error term in terms of its approximation and truncation errors; i.e., 
\begin{align*}
    || p_{p,h}^{n+1} - p_{p}^{n+1}||^2_{1} & \leq   \left(||\chi_{pp,h}^{n+1}||_{1} + ||\theta_{pp}^{n+1}||_{1} \right)^2 \leq 2\left( ||\chi_{pp,h}^{n+1}||^2_{1} + ||\theta_{pp}^{n+1}||^2_{1}  \right) .
    \end{align*}

\noindent For the final error bound, let $\mathcal{L}$ denote the left hand side of the following inequality:
\begin{align*}
    &\mathcal{L} := ||\bm{\eta}_h^M - \bm{\eta}^M||_{1}^2 +  ||\bm{u}_h^M - \bm{u}^M||_{0}^2 +  ||\bm{\dot{\eta}}_h^M - \bm{\dot{\eta}}^M||_{0}^2 +||p_{p,h}^M - p_p^M ||_{0}^2 \\ 
    &+ \Delta t \sum_{n=0}^{M-1} \Big[ ||\bm{u}_h^{n+1} - \bm{u}^{n+1}||_{1}^2 + ||\bm{\dot{\eta}}_h^{n+1} - \bm{\dot{\eta}}^{n+1}||_{1}^2 +  || p_{p,h}^{n+1} - p_{p}^{n+1}||^2_{1} \\
    &+  || g_{2,h}^{n+1} - g_{2}^{n+1}||^2_{0,\gamma} +  || \lambda_{p,h}^{n+1} - \lambda_{p}^{n+1}||^2_{1/2,\gamma} +  || p_{f,h}^{n+1} - p_{f}^{n+1}||^2_{0} +  || g_{1,h}^{n+1} - g_{1}^{n+1}||^2_{-1/2,\gamma} \Big] \\
    &\leq 2 \Big(   ||\bm{\chi}_{\eta,h}^M||_{1}^2 + ||\bm{\chi}_{u,h}^M||^2_{0} +  ||\bm{\dot{\chi}}_{\eta,h}^M||^2_{0} + ||\chi_{pp,h}^M||^2_{0}+ \Delta t \sum_{n=0}^{M-1} \Big[ ||\bm{\chi}_{u,h}^{n+1}||_{1}^2  + ||\bm{\dot{\chi}}_{\eta,h}^{n+1}||_{1}^2 +  || \chi_{pp,h}^{n+1}||_{1}^2 \\&
    +  || \chi_{g2,h}^{n+1}||_{0,\gamma}^2 +  || \chi_{\lambda p,h}^{n+1}||^2_{1/2,\gamma}+  || \chi_{pf,h}^{n+1}||^2_{0} +  || \chi_{g1,h}^{n+1}||^2_{-1/2,\gamma} \Big] + ||\bm{\theta}_{\eta}^M||_{1}^2 + ||\bm{\theta}_u^M||^2_{0} +   ||\bm{\dot{\theta}}_\eta^M||^2_{0}+  ||\theta_{pp}^M||^2_{0} \\
    &+ \Delta t \sum_{n=0}^{M-1} \Big[ ||\bm{\theta}_u^{n+1}||_{1}^2+ ||\bm{\dot{\theta}}_\eta^{n+1}||_{1}^2   + ||\theta_{pp}^{n+1}||^2_{1} 
   + ||\theta_{g2}^{n+1}||^2_{0,\gamma}  + ||\theta_{\lambda p}^{n+1}||^2_{1/2,\gamma}  + ||\theta_{pf}^{n+1}||^2_{0}  + ||\theta_{g1}^{n+1}||^2_{-1/2,\gamma}  \Big]    \Big).
    \end{align*}
Apply the bound for the truncation errors achieved in \eqref{FPSI_FEM:eq:ResultAfterGronwallError} and take the infimum of the entire inequality over $s_{1,h} \in \Lambda_{g1}^h$ and $\mu_h \in \Lambda_\lambda^h$, applying the approximation properties \eqref{FPSI_FEM:eq:Approx_StokesProjection}-\eqref{FPSI_FEM:eq:H1/2approxProp}:
\begin{align*}
   \mathcal{L} %&\leq 2\Bigg( \frac{2C^* \Delta t C_{\Delta t}}{K_9}  \sum_{n=0}^{M-1} \Big[ \mathcal{J}_\theta^{n+1} \Big]   + ||\bm{\theta}_{\eta}^M||_{1}^2 + ||\bm{\theta}_u^M||^2_{0} +   ||\bm{\dot{\theta}}_\eta^M||^2_{0} +  ||\theta_{pp}^M||^2_{0} + \Delta t \sum_{n=0}^{M-1} \Big[ ||\bm{\theta}_u^{n+1}||_{1}^2  \\
    %&+ ||\bm{\dot{\theta}}_\eta^{n+1}||_{1}^2  + ||\theta_{pp}^{n+1}||^2_{1}  + ||\theta_{g2}^{n+1}||^2_{0,\gamma}  + ||\theta_{\lambda p}^{n+1}||^2_{1/2,\gamma}  + ||\theta_{pf}^{n+1}||^2_{0}  + ||\theta_{g1}^{n+1}||^2_{-1/2,\gamma}  \Big]  \Bigg)\\
 &\leq 2\Big(  \frac{2C^* \Delta t C_{\Delta t}}{K_9}  \sum_{n=0}^{M-1} \Big[ \underset{\mu_h \in \Lambda_\lambda^h}{\inf} \underset{s_{1,h} \in \Lambda_{g1}^h}{\inf}\mathcal{J}_\theta^{n+1} \Big]  +  C h_2^{2k_\eta}||\bm{\eta}^M||_{k_\eta+1}^2 + C h_1^{2k_u}||\bm{u}^M||^2_{k_u+1} \\
    &+  C h_2^{2(k_\eta+1)}  ||\bm{\dot{\eta}}^M||^2_{k_\eta+1} +C h_2^{2(k_{pp}+1)}||p_p^M||^2_{k_{pp}+1} + \Delta t \sum_{n=0}^{M-1} \Big[ C h_1^{2k_u}||\bm{u}^{n+1}||_{k_u+1}^2   \\
     &+ C h_2^{2k_\eta}||\bm{\dot{\eta}}^{n+1}||_{k_\eta+1}^2+ C h_2^{2k_{pp}}||p_p^{n+1}||^2_{k_{pp}+1}+ C h_\gamma^{2(k_{g2}+1)}||g_2^{n+1}||^2_{k_{g2}+1,\gamma}  +  Ch_\gamma^{2k_{\lambda}} \underset{\substack{w \in H^{k_{\lambda}+1}(\Omega_r) \\ w|_\gamma = \lambda_p }}{\inf}||w||^2_{k_{\lambda}+1,\Omega_r} \\
     &+ C h_1^{2(k_{pf}+1)}||p_f^{n+1}||^2_{k_{pf}+1}  +  Ch_\gamma^{2(k_{g1}+1)} \underset{\substack{w \in H^{k_{g1}+1}(\Omega_r) \\ w|_\gamma = g_1}}{\inf}||w||^2_{k_{g1}+1,\Omega_r} \Big]       \Big).
\end{align*}
Let $\overline{C} := 1 + \dfrac{2C^* C_{\Delta t}}{K_9}$. Inserting the bound for $\mathcal{J}_\theta$ from \eqref{FPSI_FEM:eq:BoundJtheta} and simplifying yields

\begin{align}\label{FPSI_FEM:eq:AlmostFinalErrorBound}
\begin{split}
    \mathcal{L}  &\leq 2C \Big(   h_2^{2k_\eta}||\bm{\eta}^M||_{k_\eta+1}^2 +  h_1^{2k_u}||\bm{u}^M||^2_{k_u+1}+   h_2^{2(k_\eta+1)}  ||\bm{\dot{\eta}}^M||^2_{k_\eta+1}  + h_2^{2(k_{pp}+1)}||p_p^M||^2_{k_{pp}+1} \\
     &+ \Delta t \sum_{n=0}^{M-1} \Big[\frac{2C^* C_{\Delta t}}{K_9}h_1^{2k_u}||\bm{\dot{u}}^{n+1}||^2_{k_u+1} + h_1^{2k_u} \overline{C}||\bm{u}^{n+1}||_{k_u+1}^2 + \frac{2C^* C_{\Delta t}}{K_9} h_2^{2k_\eta} || \bm{\eta}^{n+1}||_{k_\eta+1}^2\\
     &+  h_2^{2k_\eta} \overline{C}||\bm{\dot{\eta}}^{n+1}||_{k_\eta+1}^2 + \frac{2C^* C_{\Delta t}}{K_9} h_2^{2(k_{pp}+1)}||p_p^{n+1}||_{k_{pp}+1}^2 +  h_2^{2k_{pp}} \overline{C}||p_p^{n+1}||^2_{k_{pp}+1} \\
     &+  h_\gamma^{2(k_{g2}+1)} \overline{C}||g_2^{n+1}||^2_{k_{g2}+1,\gamma}  +  h_\gamma^{2k_{\lambda}} \overline{C} \underset{\substack{w \in H^{k_{\lambda}+1}(\Omega_r) \\ w|_\gamma = \lambda_p}}{\inf}||w||^2_{k_{\lambda}+1,\Omega_r} +  h_1^{2(k_{pf}+1)} \overline{C}||p_f^{n+1}||^2_{k_{pf}+1}  \\
     &+ h_\gamma^{2(k_{g1}+1)} \overline{C} \underset{\substack{w \in H^{k_{g1}+1}(\Omega_r) \\ w|_\gamma = g_1}}{\inf}||w||^2_{k_{g1}+1,\Omega_r}\Big]       \Big).
     \end{split}
\end{align}
%The terms with $\frac{2C^* C_{\Delta t}}{K_9}$ can also be bounded above by $\overline{C}.$ 
Lastly, we bound $\mathcal{L}$ from below and the right hand side of the inequality from above. Recalling the definition of $\mathcal{L}$, 
\begin{align*}
    \mathcal{L} &\geq \frac{1}{7M+4} \Big( ||\bm{\eta}_h^M - \bm{\eta}^M||_{1} +  ||\bm{u}_h^M - \bm{u}^M||_{0} +  ||\bm{\dot{\eta}}_h^M - \bm{\dot{\eta}}^M||_{0} +||p_{p,h}^M - p_p^M ||_{0} \\ 
    &+ \sqrt{\Delta t} \sum_{n=0}^{M-1} \Big[ ||\bm{u}_h^{n+1} - \bm{u}^{n+1}||_{1} + ||\bm{\dot{\eta}}_h^{n+1} - \bm{\dot{\eta}}^{n+1}||_{1} +  || p_{p,h}^{n+1} - p_{p}^{n+1}||_{1} \\
    &+  || g_{2,h}^{n+1} - g_{2}^{n+1}||_{0,\gamma} +  || \lambda_{p,h}^{n+1} - \lambda_{p}^{n+1}||_{1/2,\gamma} +  || p_{f,h}^{n+1} - p_{f}^{n+1}||_{0} +  || g_{1,h}^{n+1} - g_{1}^{n+1}||_{-1/2,\gamma} \Big]    \Big)^2.
    \end{align*}
Let $C_a$ absorb constants C from the approximation errors. Noting $\dfrac{2C^* C_{\Delta t}}{K_9} \leq \overline{C}$ yields the bound in Theorem \ref{FPSI_FEM:thm:generalError}.
\end{proof}

\noindent We state this result specifically for the case tested in our numerical results in Section \ref{FPSI_FEM:sec:Numerical}. 
\begin{corollary}
        Assume sufficient regularity of the continuous solutions. With ($P_2,P_1$) used for $(\bm{u}_h, p_{f,h})$ and $P_2$ elements used for $(\bm{\eta}_h, p_{p,h})$ and each LM $g_{1,h}, g_{2,h},   \lambda_{p,h}$, the spatial error between the fully discrete and semi-discrete solutions at step $M$, for $0 < M \leq  N$ is given by 
\begin{align}
\begin{split}
      &||\bm{\eta}_h^M - \bm{\eta}^M||_{1} +  ||\bm{u}_h^M - \bm{u}^M||_{0} +  ||\bm{\dot{\eta}}_h^M - \bm{\dot{\eta}}^M||_{0} +||p_{p,h}^M - p_p^M ||_{0} \\ 
    &+ \sqrt{\Delta t} \sum_{n=0}^{M-1} \Big[ ||\bm{u}_h^{n+1} - \bm{u}^{n+1}||_{1} + ||\bm{\dot{\eta}}_h^{n+1} - \bm{\dot{\eta}}^{n+1}||_{1} +  || p_{p,h}^{n+1} - p_{p}^{n+1}||_{1} +  || g_{2,h}^{n+1} - g_{2}^{n+1}||_{0,\gamma} \\
    &+  || \lambda_{p,h}^{n+1} - \lambda_{p}^{n+1}||_{1/2,\gamma} +  || p_{f,h}^{n+1} - p_{f}^{n+1}||_{0} +  || g_{1,h}^{n+1} - g_{1}^{n+1}||_{-1/2,\gamma} \Big] \\
      &\leq C_a \sqrt{14M+8} \Big( h_2^{3} \left(  ||\bm{\dot{\eta}}^M||_{3} +||p_p^M||_{3} \right) + h_1^{2}||\bm{u}^M||_{3}+ h_2^{2}  ||\bm{\eta}^M||_{3}    \\
     &+  (\Delta t \ \overline{C})^{1/2} \sum_{n=0}^{M-1} \Big[h_1^{2} \left( ||\bm{\dot{u}}^{n+1}||_{3} +  ||\bm{u}^{n+1}||_{3} + ||p_f^{n+1}||_{2} \right) +  h_2^{2} \left( || \bm{\eta}^{n+1}||_{3}+ ||\bm{\dot{\eta}}^{n+1}||_{3} + ||p_p^{n+1}||_{3} \right) \\
     &+  h_\gamma^{3}\Big( ||g_2^{n+1}||_{H^{3}(\gamma)}   + \underset{\substack{w \in H^3(\Omega_r) \\ w|_\gamma = g_1 }}{\inf} ||w||_{3,\Omega_r} \Big) + h_2^3 ||p_p^{n+1}||_{3}  + h_\gamma^2  \underset{\substack{w \in H^{3}(\Omega_r) \\ w|_\gamma = \lambda_p}}{\inf} ||w||_{3,\Omega_r}   \Big]       \Big).
     \end{split}
\end{align}

\end{corollary}

	\section{Partitioned Method}\label{FPSI_FEM:sec:matrices-partitioned}

Let $\{v_j\}, \{q_j\}, \{\phi_j\}, \{w_j\}$, $\{s_{1,j}\}, \{s_{2,j}\}$ and $\{\mu_j\}$ be basis functions for the discrete solutions $\bm{u}_h, p_{f,h}, \bm{\eta}_h, p_{p,h}, g_{1,h}, g_{2,h}$, and $\lambda_{p,h}$ of \eqref{FPSI_FEM:WF:RestateForError_Discrete}. We may represent the discrete solutions as linear combinations of these basis functions by making use of $\bm{u},\bm{p}_f,\bm{\eta},\bm{p}_p, g_1, g_2$, and $\bm{\lambda}_p$ as the corresponding coefficient vectors; i.e., $\bm{u}_h(x,y,t) = \sum_j \bm{u}_j(t) v_j(x,y)$.  In the following, note that $\bm{u}, \bm{\eta}$ now represent the finite element coefficient vectors, instead of the unknown variables in the strong form. Substituting these linear combinations into \eqref{FPSI_FEM:WF:RestateForError_Discrete}, discretizing in time using Backward Euler, and picking appropriate basis functions as test functions yields the following linear system. 
\begin{align}\label{FPSI_FEM:MonolithicMatSystem}
\begin{split}
 M_f \bm{u}^{n+1} +  \Delta t K_f \bm{u}^{n+1} - \Delta t P_f \bm{p}_f^{n+1} - \Delta t G_{u,N}^T \bm{g}_1^{n+1} - \Delta t G_{u,\tau}^T \bm{g}_2^{n+1} &= \Delta t \bm{\overline{f}}_f^{n+1} + M_f \bm{u}^n \\
P_f^T \bm{u}^{n+1} &= \bm{0} \\
 M_\eta \bm{\eta}^{n+1} + \Delta t^2 (K_\eta + L_\eta) \bm{\eta}^{n+1} -\Delta t^2 P_p \bm{p}_p^{n+1} - \Delta t^2 G_{\eta,N}^T \bm{g}_1^{n+1} &+ \Delta t^2 G_{\eta,\tau}^T \bm{g}_2^{n+1} \\
 &= \Delta t^2 \bm{\overline{f}}_\eta^{n+1} + 2 M_\eta \bm{\eta}^n - M_\eta \bm{\eta}^{n-1}\\
M_p \bm{p}_p^{n+1} + \Delta t K_p \bm{p}_p^{n+1} + P_p^T \bm{\eta}^{n+1} - \Delta t G_p^T \bm{\lambda}_p^{n+1} &= \Delta t \bm{\overline{f}}_p^{n+1} + M_p \bm{p}_p^n + P_p^T \bm{\eta}^n \\
G_p \bm{p}_p^{n+1} + G_{1,\lambda} \bm{g}_1^{n+1}  &= \bm{0} \\
\Delta t G_{u,N} \bm{u}^{n+1}  + G_{\eta,N} \bm{\eta}^{n+1}  - \Delta t G_{1,\lambda}^T \bm{\lambda}_p^{n+1} &= G_{\eta,N} \bm{\eta}^n  \\
\Delta t M_{g2} \bm{g}_2^{n+1} +  \Delta t G_{u,\tau} \bm{u}^{n+1}  -  G_{\eta,\tau} \bm{\eta}^{n+1} &= - G_{\eta,\tau} \bm{\eta}^n .
\end{split}
\end{align}

\noindent A few notes and definitions are in order:
\begin{enumerate}
    \item $M_*$ and $K_*$ represent mass and stiffness matrices, $P_*$ are pressure matrices, $L_\eta$ contains divergence terms for the structural displacement, and both body forces and Neumann contributions to the right hand sides are contained in the vectors $\bm{\overline{f}}_f, \bm{\overline{f}}_\eta,$ and $\bm{\overline{f}}_p$. Parameters other than $\Delta t$ have been included implicitly in the definition of the corresponding matrices.
    \item Interactions between subdomain and interface bases are captured in the $G$ matrices. In two dimensions, the normal and tangential vectors $\bm{n_f}, \bm{n_p}, \bm{\tau_\gamma}$ each have two components, giving us the matrix definitions: 
\begin{align*}
        G_{u,N} = \begin{bmatrix}
        n_{f,1} (G_{u,N}^1)  & n_{f,2} (G_{u,N}^2)
        \end{bmatrix}, \hspace{5mm} &\text{ with } \Big(G_{u,N}^r\Big)_{i,j} = \langle v_j, s_{1,i} \rangle_\gamma \text{  for } r=1,2 \\
        G_{\eta,N} = \begin{bmatrix}
        n_{p,1} (G_{\eta,N}^1) & n_{p,2} (G_{\eta,N}^2)
        \end{bmatrix}, \hspace{5mm} &\text{ with } \Big(G_{\eta,N}^r\Big)_{i,j} = \langle \varphi_j, s_{1,i} \rangle_\gamma \text{  for } r=1,2 \\
         G_{u,\tau} = \begin{bmatrix}
        \tau_{\gamma,1} (G_{u,\tau}^1)  & \tau_{\gamma,2} (G_{u,\tau}^2)
        \end{bmatrix}, \hspace{5mm} &\text{ with } \Big(G_{u,\tau}^r\Big)_{i,j} = ( v_j, s_{2,i} )_\gamma \text{  for } r=1,2 \\
        G_{\eta,\tau} = \begin{bmatrix}
        \tau_{\gamma,1} (G_{\eta,\tau}^1) & \tau_{\gamma,2} (G_{\eta,\tau}^2)
        \end{bmatrix}, \hspace{5mm} &\text{ with } \Big(G_{\eta,\tau}^r\Big)_{i,j} =( \varphi_j, s_{2,i} )_\gamma \text{  for } r=1,2 \\
        \left(G_p\right)_{i,j} = ( \mu_i, w_j )_\gamma, \hspace{5mm}& \text{ and } \left(G_{1,\lambda} \right)_{i,j} = \langle  \mu_i, s_{1,j}  \rangle_\gamma.
    \end{align*}

    \end{enumerate}
  
    \begin{remark} As discussed in Section \ref{FPSI_FEM:sec:Continuous Model-3LMs}, an additional stabilization term was needed for the LM $\lambda_p$ to prove well-posedness in the continuous case and the same formulation was used for the discrete case. However, as will be observed in Section \ref{FPSI_FEM:sec:Numerical}, this stabilization is numerically unneeded.

    \end{remark}    
     
   We treat $\bm{g}_1$ and $\bm{g}_2$ together as they represent the normal and tangential components of the interface stress, defining the column vector $\bm{g} := [\bm{g}_1^T, \ \bm{g}_2^T]^T$. This grouping suggests the block matrices:
\begin{align}\label{FPSI_FEM:eq:GMgBlockDef}
\begin{split}
    G_{u}:= \begin{bmatrix}
        G_{u,N} \\ G_{u,\tau}
    \end{bmatrix}, \quad 
    G_\eta := \begin{bmatrix}
        -G_{\eta,N} \\ G_{\eta,\tau}
    \end{bmatrix}, \quad G_\lambda = \begin{bmatrix}
        G_{1,\lambda}^T \\ 0_{N_{g2} \times N_{\lambda}}
    \end{bmatrix}, \quad M_g = \begin{bmatrix}
        0_{N_{g1} \times N_{g1}} & 0_{N_{g1} \times N_{g2}} \\ 0_{N_{g2} \times N_{g1}} &   M_{g2}
    \end{bmatrix},
    \end{split}
\end{align}

\noindent where $N_*$ is the total number of degrees of freedom for a variable $*$. For vector-valued variables $\bm{u}_h$ and  $\bm{\eta}_h$, we use $\frac{1}{2} N_u$ and  $\frac{1}{2} N_\eta$ to represent the number of DoFs on which these variables are defined in one dimension. Thus in $2$-D the size of the full coefficient vectors $\bm{u}, \bm{\eta}$ are $N_u, N_\eta$ respectively. Also, let  $N_\gamma = N_{g1} + N_{g2}$ be the number of DoFs for $\bm{g}$. Next, define the following matrices and vectors: 
\begin{align}\label{FPSI_FEM:eq:Wdefns}
\begin{split}
   & W_f := M_f + \Delta t K_f, \hspace{7mm} W_\eta := M_\eta + \Delta t^2 (K_\eta + L_\eta), \hspace{7mm} W_p := M_p + \Delta t K_p \\
    &\bm{w_1}^{n+1} := \Delta t \bm{\overline{f}}_f^{n+1} + M_f \bm{u}^n, \hspace{6mm} 
    \bm{w_2}^{n+1} := \Delta t^2 \bm{\overline{f}}_\eta^{n+1} + 2 M_\eta \bm{\eta}^n - M_\eta \bm{\eta}^{n-1} \\
    &\bm{w_3}^{n+1} := \Delta t \bm{\overline{f}}_p^{n+1} + M_p \bm{p}_p^n + P_p^T \bm{\eta}^n, \quad
\bm{w_4}^{n+1} :=  -G_{\eta} \bm{\eta}^n. %\hspace{5mm} 
 %   \bm{w_5}^{n+1} := - G_{\eta,\tau} \bm{\eta}^n
 \end{split}
\end{align}

\noindent With this notation, we may rewrite the system \eqref{FPSI_FEM:MonolithicMatSystem} in simplified form as 
\begin{align}\label{FPSI_FEM:MonolithicMatSystemSimplified}
\begin{split}
 W_f \bm{u}^{n+1} - \Delta t P_f \bm{p}_f^{n+1} - \Delta t G_u^T \bm{g}^{n+1} &= \bm{w_1}^{n+1} \\
 W_\eta \bm{\eta}^{n+1} -\Delta t^2 P_p \bm{p}_p^{n+1} + \Delta t^2 G_\eta^T \bm{g}^{n+1} &= \bm{w_2}^{n+1}\\
W_p \bm{p}_p^{n+1}  + P_p^T \bm{\eta}^{n+1} - \Delta t G_p^T \bm{\lambda}_p^{n+1} &= \bm{w_3}^{n+1} \\
P_f^T \bm{u}^{n+1} &= \bm{0} \\
 \Delta t G_u \bm{u}^{n+1} - G_\eta \bm{\eta}^{n+1} + \Delta t M_{g} \bm{g}^{n+1} -\Delta t G_\lambda \bm{\lambda}_p &= \bm{w_4}^{n+1} \\
 G_p \bm{p}_p^{n+1} + G_{\lambda}^T \bm{g}^{n+1}  &= \bm{0}.
\end{split}
\end{align}
 \begin{comment}
 The dimensions of each matrix are as follows:
\begin{align*}
\begin{array}{rrrr}
    W_f: N_u \times N_u & W_\eta: N_\eta \times N_\eta &W_p: N_{pp} \times N_{pp} &  P_f: N_u \times N_{pf} \\
     P_p: N_\eta \times N_{pp} & G_u: N_\gamma \times N_u  &  G_\eta: N_\gamma \times N_\eta & G_p: N_\lambda \times N_{pp}\\
      G_\lambda: N_\gamma \times N_\lambda & G_{1,\lambda}: N_\lambda \times N_{g1} &   M_g:  N_{\gamma} \times N_\gamma &  M_{g2}: N_{g2} \times N_{g2}\\
    G_{u,N}: N_{g1} \times N_u &  G_{\eta,N}: N_{g1} \times N_\eta& G_{u,\tau}: N_{g2} \times N_u &   G_{\eta,\tau}: N_{g2} \times N_\eta.
\end{array}
\end{align*}
\end{comment}
Our aim is to derive a Schur complement equation expressing the variables functioning as LMs in our saddle point formulation (i.e., $\bm{p}_f,\bm{g},\bm{\lambda}_p$) implicitly in terms of the remaining variables. Grouping $\bm{y} := [\bm{p}_f^T, \bm{g}^T, \bm{\lambda}_p^T]^T$ gives rise to the block matrices and vector:
\begin{align*}
    &A_f := \begin{bmatrix}
         \Delta t P_f^T \\ \Delta t G_u \\ 0_{N_\lambda \times N_u} 
    \end{bmatrix} \hspace{5mm} A_1 := \begin{bmatrix}
        P_f^T \\ \Delta t G_u \\ 0_{N_\lambda \times N_u}
    \end{bmatrix} \hspace{5mm} A_2 := \begin{bmatrix}
         0_{N_{pf} \times N_\eta } \\  - G_{\eta} \\ 0_{N_\lambda \times N_\eta} 
    \end{bmatrix} \hspace{5mm} A_3 := \begin{bmatrix}
         0_{N_{pf} \times N_{pp}} \\ 0_{N_{\gamma} \times N_{pp}} \\  G_p
    \end{bmatrix} \\
    & A_4 := \begin{bmatrix}
       0_{N_{pf} \times N_{pf}} & 0_{N_{pf} \times N_\gamma} & 0_{N_{pf} \times N_\lambda} \\
       0_{N_\gamma \times  N_{pf}} & \Delta t M_g & -\Delta t G_\lambda \\
         0_{N_\lambda \times N_{pf}} &  G_\lambda^T & 0_{N_\lambda \times N_\lambda} 
    \end{bmatrix} \hspace{5mm} \bm{a}^{n+1} :=\begin{bmatrix}
        \bm{0}_{N_{pf}} \\ \bm{w_4}^{n+1}  \\ \bm{0}_{N_\lambda}
    \end{bmatrix} .
\end{align*}

Although $A_1$ and $A_f$ only differ by a factor of $\Delta t$, we have intentionally chosen not to multiply the incompressibility constraint by $\Delta t$ for the sake of weighting this term more heavily in the preconditioner to follow.
The matrix system \eqref{FPSI_FEM:MonolithicMatSystemSimplified} becomes 
\begin{align}\label{FPSI_FEM:CondensedMatSystem}
    \begin{bmatrix}
        W_f & 0 & 0 & -A_f^T \\
        0 & W_\eta & -\Delta t^2 P_p & -\Delta t^2 A_2^T \\
        0 & P_p^T & W_p & -\Delta t A_3^T \\
        A_1 & A_2 & A_3 & A_4
    \end{bmatrix} \begin{bmatrix}
        \bm{u}^{n+1} \\ \bm{\eta}^{n+1} \\ \bm{p}_p^{n+1} \\ \bm{y}^{n+1}
    \end{bmatrix} = \begin{bmatrix}
        \bm{w_1}^{n+1} \\ \bm{w_2}^{n+1} \\
        \bm{w_3}^{n+1} \\ 
        \bm{a}^{n+1}
    \end{bmatrix}.
\end{align}

\noindent  Defining the matrix $T_p := W_p + \Delta t^2 P_p^T W_\eta^{-1} P_p$, row reduction allows $\bm{y}^{n+1}$ to be expressed in terms of $\bm{u}^{n+1}, \bm{\eta}^{n+1}, \bm{p}_p^{n+1}$, yielding the Schur complement equation 
$S \bm{y}^{n+1} = \bm{b}^{n+1}$, where \begin{align}\label{FPSI_FEM:SchurCompEq}
\begin{split}
S &=   A_4 + A_1 W_f^{-1} A_f^T +\Delta t^2 A_2 W_\eta^{-1} A_2^T - (A_3 + \Delta t^2 A_2 W_\eta^{-1} P_p) T_p^{-1} (-\Delta t A_3^T + \Delta t^2 P_p^T W_\eta^{-1} A_2^T),
\\ \bm{b}^{n+1} &=  \bm{a}^{n+1} - A_1 W_f^{-1} \bm{w_1}^{n+1} - A_2 W_\eta^{-1} \bm{w_2}^{n+1} - (A_3 + \Delta t^2 A_2 W_\eta^{-1} P_p) T_p^{-1} (\bm{w_3}^{n+1} - P_p^T W_\eta^{-1} \bm{w_2}^{n+1}).
\end{split}
\end{align}

Thus, at each time step $t^{n+1}$, we propose the following steps to we solve the matrix system \eqref{FPSI_FEM:MonolithicMatSystem}, represented in simplified form by \eqref{FPSI_FEM:CondensedMatSystem}:
\begin{enumerate}
\item \textit{Update right hand sides:} Update $\bm{w_i}^{n+1}$, $i=1,\ldots,4$, as defined in \eqref{FPSI_FEM:eq:Wdefns}.
    \item \textit{Update Lagrange multipliers:} Solve $S \bm{y}^{n+1} = \bm{b}^{n+1}$ for $\bm{y}^{n+1}$, with $S$ and $\bm{b}^{n+1}$ defined in \eqref{FPSI_FEM:SchurCompEq}. As $\bm{y}^{n+1} = \begin{bmatrix}
    \bm{p}_f^{n+1} \\ \bm{g}^{n+1} \\ \bm{\lambda}_p^{n+1}
    \end{bmatrix}$, the fluid pressure or interface LMs may be extracted as needed.
    \item \textit{Update fluid velocity:} Solve the fluid system for $\bm{u}^{n+1}$:
    \begin{align*} 
    W_f \bm{u}^{n+1} = \bm{w_1}^{n+1} + A_f^T \bm{y}^{n+1} = \bm{w_1}^{n+1} + \Delta t P_f \bm{p}_f^{n+1} + \Delta t G_u^T \bm{g}^{n+1}.
    \end{align*}
    \item \textit{Update displacement and pore pressure:} Solve the system \begin{align*}
        \begin{bmatrix}
            W_\eta & -\Delta t^2 P_p \\ P_p^T & W_p
        \end{bmatrix} \begin{bmatrix}
            \bm{\eta}^{n+1} \\ \bm{p}_p^{n+1}
        \end{bmatrix} = \begin{bmatrix}
            \bm{w_2}^{n+1} +\Delta t^2 A_2^T \bm{y}^{n+1} \\ 
            \bm{w_3}^{n+1} + \Delta t A_3^T \bm{y}^{n+1}
        \end{bmatrix}.
    \end{align*}
   This step may be completed by solving the block system above, or by sequentially solving \begin{align*} T_p \bm{p}_p^{n+1} &= \bm{w_3}^{n+1} +\Delta t A_3^T \bm{y}^{n+1} - P_p^T W_\eta^{-1} (\bm{w_2}^{n+1} +\Delta t^2 A_2^T \bm{y}^{n+1}) \\
    &= \bm{w_3}^{n+1} - P_p^T W_\eta^{-1} \bm{w_2}^{n+1} + \Delta t^2 P_p^T W_\eta^{-1} G_\eta^T \bm{g}^{n+1} + \Delta t G_p^T \bm{\lambda}_p^{n+1} \end{align*} for the pore pressure $\bm{p}_p^{n+1}$, and then solving for the displacement by $$W_\eta \bm{\eta}^{n+1} = \bm{w_2}^{n+1} + \Delta t^2 P_p \bm{p}_p^{n+1} +\Delta t^2 A_2^T \bm{y}^{n+1} = \bm{w_2}^{n+1} + \Delta t^2 P_p \bm{p}_p^{n+1} - \Delta t^2 G_\eta^T \bm{g}^{n+1}.$$
\end{enumerate}
Note that steps 3 and 4 could be done in parallel; once the Schur complement equation is solved, the subdomains $\Omega_f, \Omega_p$ are effectively decoupled. 
The Schur complement equation in Step 2 is problematic at first glance due to the number of matrix inversions required; we dedicate the rest of this section to discussing its solution. 

\subsection{Efficient Computation of the Schur Complement Equation}\label{FPSI_FEM:subsec:SchurSolving}
The use of a direct method to solve the Schur complement equation \eqref{FPSI_FEM:SchurCompEq} is undesirable as it requires explicit construction of the matrix $S$. We also want to avoid the direct formation of matrix inverses of $W_f, W_\eta, T_p$. 
These three matrices frequently appear in linear solves at each time step, both in the Schur complement equation as well as in the subdomain updates for $\bm{u}^{n+1}, \bm{\eta}^{n+1},\bm{p}_p^{n+1}$. As they are invariant in time, we suggest computing their Cholesky factorizations at the first time step. For $s \in \{f,\eta\}$, we have $P_{C,s}^T W_s P_{C,s} = R_s^T R_s$, for upper triangular $R_s$ and permutation matrix $P_{C,s}$ which improves the sparsity structure of the factors $R_s$. 
Then the system $W_\eta \bm{z} = \bm{b}$, for example, is solved as:
\begin{equation}
    R_{\eta}^T \bm{\bm{x_1}} = P_{C,\eta}^T \bm{b} \quad \rightarrow \quad  R_\eta \bm{x_2} = \bm{\bm{x_1}} \quad \rightarrow \quad \bm{z} = P_{C,\eta} \bm{x_2}.
\end{equation}
For the matrix $T_p = W_p + \Delta t^2 P_p^T W_\eta^{-1} P_p$, we can utilize the Cholesky factorization of $W_\eta$. With $Q = R_\eta^{-T} P_{C,\eta}^T P_p$, we see that $T_p = W_p + \Delta t^2 Q^T Q$, and we may compute the Cholesky factorization of the dense $T_p$. Computation of the inverse of upper triangular matrix $R_\eta$ is needed, which can be done by substitution. However, since these three Cholesky factorizations allow us to solve both the Schur complement equation and the subdomain updates efficiently, this one-time cost is worth it. Future work will focus on the implementation of reduced order models to reduce the matrix sizes in these computations. Also, we note that the Cholesky factorizations provide a faster alternative to LU, and in our numerical experiments, factoring $T_p$ explicitly has proven to give large speedups over implementations where the larger block $\begin{bmatrix}
    W_\eta &-\Delta t^2 P_p \\ P_p^T & W_p
\end{bmatrix}$ is factored and relevant information extracted from there.
With these three Cholesky factorizations, the right hand side of the Schur complement equation \eqref{FPSI_FEM:SchurCompEq} can be constructed through linear solves with the Cholesky factors of $W_\eta, W_f, $ and $T_p$.

In the rest of this section, we propose an iterative method and a preconditioner to solve the Schur complement equation. The use of an iterative solver such as the biconjugate gradient stabilized method (BiCGStab) with preconditioner only requires the computation of matrix-vector products instead of the explicit construction of $S$.
 BiCGStab is an appealing option as it is intended for nonsymmetric systems, has better convergence properties than BiCG, and does not require multiplication with the transpose of the matrix. 

\subsubsection{Matrix-Vector Products with $S$}\label{FPSI_FEM:sec:MatVecs}
In order to compute the action of $S$ on a vector, expand $S$ as defined in \eqref{FPSI_FEM:SchurCompEq}. This yields $S := \begin{bmatrix}
        S_1 & S_2 & S_3
    \end{bmatrix}$, where
    \begin{align*}
    S_1 &= \begin{bmatrix}
         \Delta t P_f^T W_f^{-1} P_f \\  \Delta t^2 G_{u} W_f^{-1} P_f \\ 0
    \end{bmatrix},\\
    S_2 &= \begin{bmatrix}
         \Delta t P_f^T W_f^{-1} G_u^T \\ \Delta t M_g + \Delta t^2 G_{u} W_f^{-1} G_u^T + \Delta t^2 G_{\eta} W_\eta^{-1} G_\eta^T - \Delta t^4  G_{\eta} W_\eta^{-1} P_p T_p^{-1} P_p^T W_\eta^{-1} G_\eta^T \\   G_\lambda^T + \Delta t^2 G_p T_p^{-1} P_p^T W_\eta^{-1} G_\eta^T 
    \end{bmatrix},\\
    S_3 &= 
    \begin{bmatrix}
         0 \\
         -\Delta t G_{\lambda} - \Delta t^3 G_{\eta} W_\eta^{-1} P_p T_p^{-1} G_p^T \\
            \Delta t G_p T_p^{-1} G_p^T
    \end{bmatrix} . 
\end{align*}

\noindent The product $S \bm{y}$ is equivalent to $ S_1 \bm{p}_f + S_2 \bm{g} + S_3 \bm{\lambda}_p.$ %In this product, there are seven unique terms which involve products with matrix inverses:
%\begin{align*}
%    \bm{s_1} &:= W_f^{-1} P_f \bm{p}_f, \qquad
%    \bm{s_2} := W_\eta^{-1} G_\eta^T \bm{g}, \qquad
%    \bm{s_3} :=  T_p^{-1} P_p^T \bm{s_2}, \qquad
%    \bm{s_4} := W_f^{-1} G_u^T \bm{g}\\
%    \bm{s_5} &:= W_\eta^{-1} P_p \bm{s_3}, \qquad
%    \bm{s_6} := T_p^{-1} G_p^T \bm{\lambda}_p, \qquad
%    \bm{s_7} := W_\eta^{-1} P_p  \bm{s_6}.
%\end{align*}
%Then the product becomes:
%\begin{align*}
%    S \bm{y} &= \begin{bmatrix}
%        \Delta t P_f^T (\bm{s_1} +  \bm{s_4}) \\
%       \Delta t M_g \bm{g} + \Delta t^2 G_{u} (\bm{s_1} + \bm{s_4}) + \Delta t^2 G_{\eta} (\bm{s_2} - \Delta t^2 \bm{s_5} - \Delta t \bm{s_7}) - \Delta t G_{\lambda} \bm{\lambda}_p  \\
%   G_\lambda^T \bm{g} + \Delta t G_p  ( \Delta t \bm{s_3} + \bm{s_6}) 
 %   \end{bmatrix}.
%\end{align*}
Let $\bm{z}_j$, $j=1,2,3,4$, be the solutions to the following linear systems:
\begin{align}\label{eqn:zEqns}
\begin{split}
    W_f \bm{z_1} &= P_f \bm{p}_f + G_u^T \bm{g}, \qquad \qquad
   W_\eta \bm{z_2} =   G_\eta^T \bm{g} \\
     T_p \bm{z_3} &=   \Delta t P_p^T \bm{z_2} + G_p^T \bm{\lambda}_p, \qquad W_\eta \bm{z_4} =    P_p \bm{z_3} .
     \end{split}
\end{align}
When  $S \bm{y}$ is needed in the iterative solver for the Schur complement equation, solve the four linear systems in \eqref{eqn:zEqns} for a given $\bm{y} = [\bm{p}_f^T, \bm{g}^T, \bm{\lambda}_p^T]^T$, and output
%\begin{enumerate}
%    \item Solve the system $W_f \bm{z_1} = P_f \bm{p}_f + G_u^T \bm{g}$.
%    \item  Solve the system $W_\eta \bm{z_2} = G_\eta^T \bm{g}$.
%    \item Solve the system $T_p \bm{z_3} = \Delta t P_p^T \bm{z_2} + G_p^T \bm{\lambda}_p$.
%    \item Solve the system $W_\eta \bm{z_4} =   P_p \bm{z_3}$.
%\end{enumerate}
\begin{align*}
    S \bm{y} &= \begin{bmatrix}
        \Delta t P_f^T \bm{z_1} \\
        \Delta t M_g \bm{g} + \Delta t^2 G_{u} \bm{z_1} + \Delta t^2 G_{\eta} (\bm{z_2} - \Delta t \bm{z_4}) - \Delta t G_{\lambda} \bm{\lambda}_p \\   G_\lambda^T \bm{g} + \Delta t G_p  \bm{z_3} 
    \end{bmatrix}.
 %   &= \begin{bmatrix}
   %     \Delta t P_f^T \bm{z_1}\\
     % \Delta t \begin{bmatrix}
    %       0_{N_{g1},1} \\  M_{g2} \bm{g}_2
     %  \end{bmatrix}+ \Delta t^2 G_u \bm{z_1} + \Delta t^2 G_\eta  (\bm{z_2} - \Delta t \bm{z_4}) - \Delta t \begin{bmatrix}
       %     G_{1,\lambda}^T \bm{\lambda}_p\\ 0_{N_{g2},1}
       % \end{bmatrix}\\
       % G_{1,\lambda}\bm{g}_1 + \Delta t G_p \bm{z_3}
   % \end{bmatrix}
\end{align*}

\subsubsection{Preconditioner for $S$}\label{FPSI_FEM:sec:Preconditioner}
For a general survey of preconditioners for general saddle point problems, we refer to \cite{Benzi_2005}.
Observing the structure of the Schur complement matrix, expressed explicitly in Section \ref{FPSI_FEM:sec:MatVecs}, we see that $S$ may be rewritten as a sum such that $S = \sum_{k=0}^4 \Delta t^k A_k$, where for $Q:= G_\eta W_\eta^{-1} P_p$,
\begin{align*}
A_0 &:= 
    \begin{bmatrix}
        0_{N_{pf} \times N_{pf}} & 0_{N_{pf} \times {N_\gamma}} & 0_{ N_{pf} \times N_\lambda} \\ 
        0_{N_{\gamma} \times N_{pf}} & 0_{N_{\gamma} \times N_{\gamma } } & 0_{N_{\gamma} \times N_\lambda} \\ 0_{N_\lambda \times N_{pf}} & 
       G_{\lambda}^T & 0_{N_\lambda \times N_{\lambda}}
    \end{bmatrix}, \hspace{5mm}     
    A_3 := \begin{bmatrix}
          0_{N_{pf} \times N_{pf}} & 0_{N_{pf} \times {N_\gamma}} & 0_{ N_{pf} \times N_\lambda} \\ 
        0_{N_{\gamma} \times N_{pf}} & 0_{N_{\gamma} \times N_\gamma} & -G_{\eta} W_\eta^{-1} P_p T_p^{-1} G_p^T 
        \\ 0_{N_\lambda \times N_{pf}} & 0_{N_\lambda \times N_\gamma} & 0_{N_\lambda \times N_\lambda}
    \end{bmatrix} \\     
    A_1 &:= 
    \begin{bmatrix}
        P_f^T W_f^{-1} P_f & P_f^T W_f^{-1} G_u^T & 0_{N_{pf} \times N_\lambda} \\
         0_{N_{\gamma} \times  N_{pf}} & M_g & -G_{\lambda} \\ 0_{N_\lambda \times N_{pf}} & 0_{N_\lambda \times N_\gamma} & G_p T_p^{-1} G_p^T        
    \end{bmatrix}, \hspace{5mm}   
    A_4 := \begin{bmatrix}
          0_{N_{pf} \times N_{pf}} & 0_{N_{pf} \times {N_\gamma}} & 0_{ N_{pf} \times N_\lambda} \\ 
        0_{N_{\gamma} \times N_{pf}} & -Q T_p^{-1} Q^T & 0_{N_{\gamma } \times N_\lambda}
        \\ 0_{N_\lambda \times N_{pf}} & 0_{N_\lambda \times N_\gamma} & 0_{N_\lambda \times N_\lambda}
    \end{bmatrix} \\        
     A_2&:=
    \begin{bmatrix}
         0_{N_{pf} \times N_{pf}} & 0_{N_{pf} \times {N_\gamma}} & 0_{ N_{pf} \times N_\lambda} \\ 
         G_{u} W_f^{-1} P_f & G_{u} W_f^{-1} G_u^T + G_{\eta} W_\eta^{-1} G_\eta^T & 0_{ N_{\gamma} \times N_{\lambda}} \\
         0_{N_\lambda \times N_{pf}} & G_p T_p^{-1} P_{p}^T W_\eta^{-1}G_\eta^T & 0_{N_\lambda \times N_\lambda}
    \end{bmatrix}.
\end{align*}

As $\Delta t$ is small, the approximation $S \approx A_0 + \Delta t A_1$ can serve as the preconditioner, $M_{pre}$. 
If desired, the lower right $N_\lambda \times N_\lambda$ block 
of $A_1$, which involves the inverse of $T_p$, may be left out of the preconditioner as it is the smallest block in $M_{pre}$. Numerically we have observed that including this block in $M_{pre}$ does not improve the accuracy of the result, but it improves the convergence of the iterative solver. In Section \ref{FPSI_FEM:sec:Numerical} we show how the inclusion of this block affects the iteration count. The preconditioner has the form:
\begin{align}\label{FPSI_FEM:eqn:preconditioner}
    M_{pre} := A_0 + \Delta t A_1 = 
    \begin{bmatrix}
        \Delta t P_f^T W_f^{-1} P_f & \Delta t P_f^T W_f^{-1} G_u^T & 0_{N_{pf} \times N_\lambda} \\
         0_{N_{\gamma} \times N_{pf}} & \Delta t M_g & -\Delta t G_{\lambda} \\ 
         0_{N_\lambda \times N_{pf}} & G_\lambda^T & ( G_p T_p^{-1} G_p^T )
    \end{bmatrix}, 
\end{align}
where the lower right block is in parentheses to emphasize the fact that it is optional. As before, we do not explicitly construct $M_{pre}$; the solver only requires the product $\bm{a} = M_{pre}^{-1} \bm{x}$ for a given $\bm{x}$. Expand the definitions of $M_g, G_\lambda$ from \eqref{FPSI_FEM:eq:GMgBlockDef} and compute $\bm{a}$ by solving the system
\begin{align}\label{FPSI_FEM:eqn:Mpreinv}
     \begin{bmatrix}
        \Delta t P_f^T W_f^{-1} P_f & \Delta t P_f^T W_f^{-1} G_u^T & 0_{N_{pf} \times N_\lambda} \\
         0_{N_{g1} \times N_{pf}} & [0_{N_{ \gamma 1} \times N_{g1}} \ 0_{N_{g1} \times N_{g2}}] & -\Delta t G_{1,\lambda}^T \\ 
        0_{N_{g2} \times N_{pf}} & [0_{N_{g2} \times N_{g1}} \hspace{0.8mm} \Delta t   M_{g2}] & 0_{N_{g2} \times N_\lambda} \\ 0_{N_\lambda \times N_{pf}} & [G_{1,\lambda} \hspace{1mm} 0_{N_\lambda \times N_{g2}}] & ( G_p T_p^{-1} G_p^T )
    \end{bmatrix} \begin{bmatrix}
        \bm{a_1} \\ \bm{a_{21}} \\ \bm{a_{22}} \\ \bm{a_3}
    \end{bmatrix} = \begin{bmatrix}
        \bm{x_1} \\ \bm{x_{21}} \\ \bm{x_{22}} \\ \bm{x_3}
    \end{bmatrix},
\end{align}
where $\bm{a_1}, \bm{x_1} \in \mathbb{R}^{N_{pf}}$, $\bm{a_{21}}, \bm{x_{21}} \in \mathbb{R}^{N_{g1}}$, $\bm{a_{22}}, \bm{x_{22}} \in \mathbb{R}^{N_{g2}}$, and $\bm{a_3}, \bm{x_3} \in \mathbb{R}^{N_\lambda}$. 
Since the discrete functions $g_{1,h}, \lambda_{p,h}$ are each defined in a subspace of $H^{1/2}(\gamma)$, $\Lambda^{h_\gamma}_{g1}$ can be chosen such that $\Lambda^{h_\gamma}_{g1} = \Lambda^{h_\gamma}_\lambda$ with identical basis functions. This means that the matrix $G_{1,\lambda}$, composed of the inner product of basis functions for $\lambda_{p,h}$ with basis functions for $g_{1,h}$, can always be taken to be an invertible mass matrix over the discrete finite element space $\Lambda_{g1}^{h_\gamma} = \Lambda_{\lambda}^{h_\gamma}$. 

From the matrix system \eqref{FPSI_FEM:eqn:Mpreinv}, we may directly derive the following:
\begin{equation}\label{FPSI_FEM:eqn:solvea3}
    \bm{a_3} = -\frac{1}{\Delta t} G_{1,\lambda}^{-T} \bm{x_{21}}, \quad \bm{a_{22}} = \frac{1}{\Delta t} M_{g2}^{-1} \bm{x_{22}}, \quad
	\bm{a_{21}} = \begin{cases}
	    G_{1,\lambda}^{-1} \bm{x_3}, \text{ or} \\
        G_{1,\lambda}^{-1}\left(\bm{x_3} - G_p T_p^{-1}G_p^T \bm{a_3}\right).
	\end{cases} 
\end{equation}
The two definitions for $\bm{a_{21}}$ depend on the inclusion of the block $G_p T_p^{-1} G_p^T$ in $M_{pre}$. If it is included, the product with $T_p^{-1}$ may be efficiently computed using its Cholesky factorization, as discussed toward the start of Section \ref{FPSI_FEM:subsec:SchurSolving}.
Lastly, we find $\bm{a_1}$ by
\begin{align}\label{FPSI_FEM:eqn:solvea1}
     P_f^T W_f^{-1} P_f \bm{a_1} = \frac{1}{\Delta t}\bm{x_1} -  P_f^T W_f^{-1} G_u^T \begin{bmatrix}
        \bm{a_{21}} \\ \bm{a_{22}}
    \end{bmatrix}.
\end{align}
The matrix on the left hand side may be computed and factored in a similar way to that of constructing $T_p$. %Defining $Q = R_f^{-T} P_{C,f}^T P_f$, we see that $P_f^T W_f^{-1} P_f =  Q^T Q$, with $R_f, P_{C,f}$ the Cholesky factor and permutation matrix for $W_f$. One could factor $Q^T Q$ once on the first time step as well and use its factors in the computation of \eqref{FPSI_FEM:eqn:solvea1}. 
If wanting to avoid this build, the matrix may be expanded into a larger saddle point system. In summary, if given a vector $\bm{x} = [
        \bm{x_1}^T, \bm{x_{21}}^T, \bm{x_{22}}^T,  \bm{x_3}^T]^T$, the matrix-vector product $\bm{a} = M_{pre}^{-1} \bm{x}$ is given by the solution of equations \eqref{FPSI_FEM:eqn:solvea3} -- \eqref{FPSI_FEM:eqn:solvea1}.

	\section{Numerical Results}\label{FPSI_FEM:sec:Numerical}
In this section, we verify the performance of our proposed algorithm for the FPSI system by testing the spatial  convergence on a manufactured solution and then implementing a hydrological example. In particular, the effect of the proposed preconditioner $M_{pre}$ \eqref{FPSI_FEM:eqn:preconditioner} on the iterative solver for the Schur complement equation is examined. %After these studies of the full order model, we turn to both reproductive and predictive tests for the ROM-ROM coupling.

As each LM $g_1, g_2, \lambda_p$ spatially occupies the same region, $\gamma$, we take the mesh for each of the FEM spaces $\Lambda_{g1}^{h_\gamma}$, $\Lambda_{g2}^{h_\gamma}$, and $\Lambda_{\lambda}^{h_\gamma}$ to be the restriction to the interface of the mesh on $\Omega_p$.

\subsection{FEM-FEM coupling} 
\label{FPSI_FEM:sec:ManSol} 
We begin our studies with a manufactured solution from \cite{Caucao_2022} to check convergence of the method in time and space. On $\Omega_f = [0,1] \times [0,1]$ and $\Omega_p = [0,1] \times [-1,0]$, the fluid velocity $\bm{u} = [u_1, u_2]^T$, fluid pressure $p_f$, structural displacement $\bm{\eta} = [\eta_1, \eta_2]^T$, and pore pressure $p_p$ are defined as:
\begin{align}\label{FPSI_FEM:eq:manufacturedSoln}
\begin{split}
u_1  &= \pi \cos(\pi t) \left(
    -3x + \cos(y)\right), \hspace{5mm}
u_2 = \pi \cos(\pi t) \left(
    y+1\right), \hspace{5mm}
p_f = e^t \sin(\pi x) \cos\left(\frac{\pi y}{2}\right) + 2\pi \cos(\pi t)\\
\eta_1 &= \sin(\pi t) \left(
    -3x + \cos(y)\right), \hspace{8mm}
\eta_2 =  \sin(\pi t) \left(
    y+1\right), \hspace{9mm} p_p = e^t \sin(\pi x) \cos\left(\frac{\pi y}{2}\right).
    \end{split}
\end{align}
Forcing functions and Neumann boundary conditions are derived from the system equations \eqref{FPSI_FEM:StokesMom}-\eqref{FPSI_FEM:Structure2}. For the displacement, we use Dirichlet conditions on each boundary, while for the fluid velocity and pore pressure, Neumann conditions are applied on the boundaries $x=0$ and $x=1$ with Dirichlet elsewhere. All parameters are set to one, and $P_2$ elements are implemented for the LMs $g_1, g_2, \lambda_p$. Taylor-Hood elements are used for fluid velocity and pressure, and $P_2$ elements are used for both the displacement and pore pressure. Preconditioned BiCGstab(l) with a tolerance of $10^{-8}$ is used to solve the Schur complement equation $SM_{pre}^{-1} \bm{x}^{n+1} = \bm{b}^{n+1}$, where $\bm{x}^{n+1} = M_{pre} \bm{y}^{n+1}$, as defined in \eqref{FPSI_FEM:SchurCompEq} and \eqref{FPSI_FEM:eqn:preconditioner}.

In Tables \ref{FPSI_FEM:table:SpaceConvMixedBCs} and \ref{FPSI_FEM:table:SpaceConvMixedBCs2}, we examine convergence in space, taking 10 time steps with two different final times. The tables list both $L^2$ and $H^1$ errors along with the convergence rates in parentheses; unless otherwise stated, all errors in this section are computed at the final time $T$. We observe $O(\Delta x^3)$ convergence in $L^2$ norms for the displacement, pore pressure, and fluid velocity, and $O(\Delta x^2)$ convergence for $H^1$ norms as well as the $L^2$ norm for the fluid pressure, which agree with the expected rates. 
\begin{comment}{\color{red} Are these better than (3.18)? The estimate says linear convergence in $h_2$ and $h_{\gamma}$.  If you use $P_2$ for LMs and $p_p$, the estimate (3.18) is improved to overall quadratic convergence, which supports the results in the tables if you can get quadratic convergence for $H_1$ errors of $p_p$ (I expect this is true).   Present $H^1$ errors of $p_p$ in the tables since the estimate is for the $H^1$ errors. I think it is fine to use $P^2$
for $p_p$ since the $s_0$ and $\kappa$ terms in the Darcy equation work as stabilization terms, so no inf-sup stable pair needed.} \\
\end{comment}
The smaller time step (Table \ref{FPSI_FEM:table:SpaceConvMixedBCs}) yields slightly better rates for the structural displacement and fluid pressure, whereas the pore pressure converges at a marginally higher rate for the larger time step (Table \ref{FPSI_FEM:table:SpaceConvMixedBCs2}).
\begin{table}[h!]
\centering
\begin{adjustbox}{max width=\textwidth}
\begin{tabular}{c|c|c|c|c|c|c|c}
 $\Delta x$     & $||\bm{\eta} - \bm{\eta}_h||_0$ &   $||\bm{\eta} - \bm{\eta}_h||_1$ & $||p_p - p_{p,h}||_0$ & $||p_p - p_{p,h}||_1$ & $||\bm{u} - \bm{u}_h||_0$ &  $ ||\bm{u} - \bm{u}_h||_1$  &   $||p_f - p_{f,h}||_0$  \\
\hline
1/2 & 1.16e-08 & 1.50e-07 & 1.58e-02 & 2.32e-01 & 1.15e-03 & 1.50e-02  & 6.39e-02 \\
1/4 & 1.47e-09 (2.98) & 3.80e-08 (1.98) & 2.06e-03 (2.94) & 6.05e-02 (1.94) & 1.47e-04 (2.98)  & 3.80e-03 (1.98) & 1.77e-02 (1.85) \\
1/8 & 1.84e-10 (2.99) & 9.54e-09 (1.99) & 2.61e-04 (2.99)& 1.53e-02 (1.99)& 1.84e-05 (2.99) & 9.54e-04 (1.99) & 4.01e-03 (2.14)  \\
1/16 & 2.30e-11 (3.00) & 2.39e-09 (2.00) & 3.26e-05 (3.00)& 3.82e-03 (2.00) & 2.31e-06 (2.99) & 2.40e-04 (1.99) & 9.74e-04 (2.04) \\
1/32 & 2.88e-12 (3.00) & 5.97e-10 (2.00) & 4.08e-06  (3.00)& 9.56e-04 (2.00) & 2.91e-07 (2.99) & 6.08e-05 (1.98) & 2.42e-04 (2.01) \\
1/64 & 3.66e-13 (2.97) & 1.53e-10 (1.96) & 6.29e-07 (2.70) & 2.59e-04 (1.89) & 3.93e-08 (2.89) & 1.60e-05 (1.92) & 6.10e-05 (1.99)
\end{tabular}
\end{adjustbox}
\caption{Convergence in space for $\Delta t = 10^{-6}$; $T=10^{-5}$}
 \label{FPSI_FEM:table:SpaceConvMixedBCs}
\end{table}

\begin{table}[h!]
\centering
\begin{adjustbox}{max width=\textwidth}
\begin{tabular}{c|c|c|c|c|c|c|c}
$\Delta x$     & $||\bm{\eta} - \bm{\eta}_h||_0$ &   $||\bm{\eta} - \bm{\eta}_h||_1$ & $||p_p - p_{p,h}||_0$ & $||p_p - p_{p,h}||_1$ & $||\bm{u} - \bm{u}_h||_0$ &  $ ||\bm{u} - \bm{u}_h||_1$  &   $||p_f - p_{f,h}||_0$  \\
\hline
1/2 & 1.16e-07 & 1.50e-06 & 1.58e-02 & 2.31e-01 & 1.15e-03 & 1.50e-02 & 6.36e-02 \\
1/4 & 1.47e-08 (2.98) & 3.80e-07 (1.98) & 2.05e-03 (2.94)  & 6.00e-02 (1.94)& 1.47e-04 (2.97) & 3.83e-03 (1.97) & 1.74e-02 (1.87) \\
1/8 & 1.84e-09 (2.99) & 9.55e-08 (1.99) & 2.58e-04 (2.99) & 1.51e-02 (1.99) & 1.90e-05 (2.95) & 1.00e-03 (1.93) & 3.96e-03 (2.14)  \\
1/16 & 2.31e-10 (3.00) & 2.40e-08 (1.99) & 3.24e-05 (3.00) & 3.80e-03 (1.99) & 2.41e-06 (2.98) & 2.57e-04 (1.97) & 9.74e-04 (2.02) \\
1/32 & 2.92e-11 (2.98) & 6.08e-09 (1.98) & 4.07e-06  (2.99) & 9.54e-04 (1.99) & 2.95e-07 (3.03) & 6.21e-05 (2.05) & 2.48e-04 (1.97) \\
1/64 & 4.32e-12 (2.76) & 1.62e-09 (1.91) & 5.10e-07 (2.99) & 2.39e-04 (2.00) & 4.13e-08 (2.84) & 1.52e-05 (2.03) & 8.17e-05 (1.60)
\end{tabular}
\end{adjustbox}
\caption{Convergence in space for $\Delta t = 10^{-5}$; $T=10^{-4}$}
 \label{FPSI_FEM:table:SpaceConvMixedBCs2}
\end{table}

% didn't do analysis for temporal convergence, so can leave out
\begin{comment}
Next, we fix $\Delta x = 1/32$ and are able to observe first order convergence in time in Table \ref{FPSI_FEM:table:TimeConvMixedBCs}, as expected for our discretization.

\begin{table}[]
    \centering
    \begin{adjustbox}{max width=\textwidth}
\begin{tabular}{c|c|c|c|c|c|c}
$\Delta t$     & $||\bm{\eta} - \bm{\eta}_h||_0$ &   $||\bm{\eta} - \bm{\eta}_h||_1$ & $||p_p - p_{p,h}||_0$ &  $||\bm{u} - \bm{u}_h||_0$ &  $ ||\bm{u} - \bm{u}_h||_1$  &   $||p_f - p_{f,h}||_0$  \\
 \hline
 1/20& 8.45e-03 & 4.94e-02 & 1.57e-02 & 4.49e-02 & 1.22e-01 & 4.17e-01 \\
 1/40& 4.10e-03 (1.04) & 2.61e-02 (0.92) & 8.09e-03 (0.96) & 2.35e-02 (0.93) & 6.43e-02 (0.92) & 2.12e-01 (0.98) \\
 1/80& 2.00e-03 (1.03) & 1.35e-02 (0.95) & 4.10e-03 (0.98) & 1.21e-02 (0.96) & 3.31e-02 (0.96) & 1.07e-01 (0.99)\\
 1/160& 9.88e-04 (1.02) & 6.87e-03 (0.97) & 2.09e-03 (0.98) & 6.12e-03 (0.98) & 1.68e-02 (0.98) & 5.37e-02 (0.99) \\
 1/320& 4.89e-04 (1.02) & 3.46e-03 (0.99) & 1.13e-03 (0.89) & 3.09e-03 (0.99) & 8.47e-03 (0.99) & 2.70e-02 (0.99)\\
 1/640& 2.42e-04 (1.02) & 1.73e-03 (1.00) & 7.49e-04 (0.59) & 1.56e-03 (0.99) & 4.26e-03 (0.99) & 1.36e-02 (0.99)    
    \end{tabular}
    \end{adjustbox}
    \caption{Convergence in time for $\Delta x = 1/32$; $T = 0.1$}
    \label{FPSI_FEM:table:TimeConvMixedBCs}
\end{table}
\end{comment}
To demonstrate the effect of the preconditioner on both the condition number of the Schur complement matrix, $\kappa(S),$ and the convergence of BiCGstab(l), Table \ref{FPSI_FEM:table:ItersConds} compares the results of solving the Schur complement equation \eqref{FPSI_FEM:SchurCompEq} without a preconditioner to the results obtained using the preconditioner $M_{pre}$ derived in Section \ref{FPSI_FEM:sec:Preconditioner}. As discussed, the lower right hand block $ G_p T_p^{-1} G_p^T $ of the preconditioner is optional; results are shown for the case where this block is included in $M_{pre}$ (referred to as ``LB'') and for the case where it is excluded (``no LB'').  The number of iterations reported is the average number of iterations of BiCGstab(l) per time step over the 10 time steps. The implementation of BiCGstab(l) used contains two inner iterations per outer iteration, so a fractional number of average iterations is possible at each time step. For the N/A number of iterations, the algorithm reached the maximum number of iterations (100) without converging. 

\begin{table}[h!]
\centering
\begin{adjustbox}{max width=\textwidth}
\begin{tabular}{|c||c|c||c|c||c|c|}
\hline
 $\Delta x$  & $\kappa(S)$&  iters, S& $\kappa(S M_{pre}^{-1})$, no LB&  iters, $SM_{pre}^{-1}$, no LB & $\kappa(S M_{pre}^{-1})$, LB & iters, $SM_{pre}^{-1}$, LB \\
\hline
1/2& 5.9e5& 18.4& 4.9e2 & 1 & 1.0& 0.5 \\
1/4& 5.5e5& 44.8& 1.9e3& 1 & 1.0& 0.5 \\
1/8& 5.5e5  & 99.0& 7.4e3& 1 & 1.0& 0.5 \\
1/16& 5.5e5 & N/A & 2.5e4 & 1 &  1.1& 0.5 \\
1/32& 5.5e5 & N/A & 6.4e4 & 1 &1.2& 0.5 \\
\hline
\end{tabular}
\end{adjustbox}
\caption[Effect of preconditioner on BiCGstab(l): iteration count and $\kappa(S)$]{Condition numbers and average number of iterations using BiCGstab(l) for the Schur complement equation, with $\Delta t = 10^{-5}$, $T=10^{-4}$. LB refers to the optional lower right hand block in the preconditioner $M_{pre}$ \eqref{FPSI_FEM:eqn:preconditioner}. }
 \label{FPSI_FEM:table:ItersConds}
\end{table}

\begin{comment} %table from dissertation, P1 elements for LM and p_p
    \begin{table}[h!]
\centering
\begin{tabular}{|c||c|c||c|c||c|c|}
\hline
 $\Delta x$  & $\kappa(S)$&  iters, S& $\kappa(S M_{pre}^{-1})$, no LB&  iters, $SM_{pre}^{-1}$, no LB & $\kappa(S M_{pre}^{-1})$, LB & iters, $SM_{pre}^{-1}$, LB \\
\hline
1/2& 3.8e5& 9.6& 8.6e1& 1 & 1.0& 0.5 \\
1/4& 3.9e5& 20.8& 3.4e2& 1 & 1.0& 0.5 \\
1/8& 4.0e5& 52.0& 1.3e3& 1 & 1.0& 0.5 \\
1/16& 4.0e5& N/A & 5.0e3& 1 &  1.0& 0.5 \\
1/32& 4.0e5& N/A & 1.7e4& 1 &1.1& 0.5 \\
\hline
\end{tabular}
\caption[Effect of preconditioner on BiCGstab(l): iteration count and $\kappa(S)$]{Condition numbers and average number of iterations using BiCGstab(l) for the Schur complement equation, with $\Delta t = 10^{-5}$, $T=10^{-4}$. LB refers to the optional lower right hand block in the preconditioner $M_{pre}$ \eqref{FPSI_FEM:eqn:preconditioner}. }
 \label{FPSI_FEM:table:ItersConds}
\end{table}
\end{comment}

% rel residual: || b - Sy|| / || b|| for each iteration
As Table \ref{FPSI_FEM:table:ItersConds} shows, use of a preconditioner drastically improves convergence of the iterative solver for the Schur complement equation. Even without the use of the lower right hand block, use of $M_{pre}$ decreases the condition number and allows for convergence of the solver on smaller mesh sizes. However, including the lower block in $M_{pre}$ actually reduces the condition number of the system matrix to one for this manufactured solution.

\begin{figure}[h]
    \centering
    \includegraphics[width=0.5\textwidth]{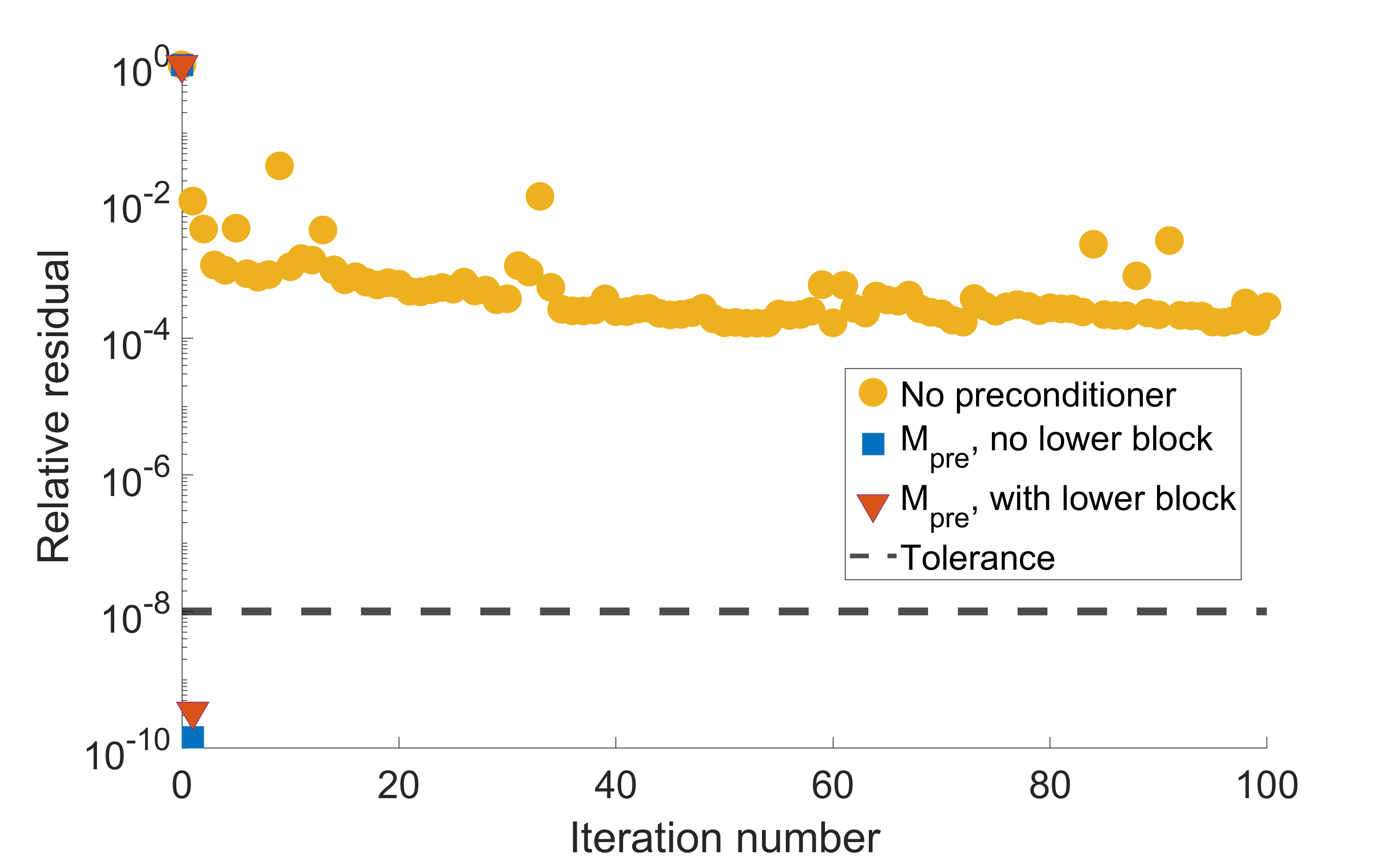}
    \caption[Residual vs. iteration number for manufactured solution]{Manufactured solution: relative residual versus iteration number for BiCGstab(l). }
    \label{FPSI_FEM:fig:ManSolution_Pre_bicg}
\end{figure}

In the case where $\Delta x = 1/32$, we see that without the preconditioner, BiCGstab(l) reaches the maximum number of iterations without converging. Although the Schur complement equation did not converge to a solution, the $L^2$ errors themselves are reasonable, with $||\bm{\eta} - \bm{\eta}_h||_0 $ = 1.32e-9, $ || p_p - p_{p,h}||_0$= 7.22e-4, $||\bm{u}-\bm{u}_h||_0  $=1.19e-4, and $ ||p_f-p_{f,h}||_0$ = 8.10. % Reducing the tolerance from $10^{-8}$ to $10^{-4}$ does not improve convergence. Setting the tolerance to $10^{-2}$ allows the Schur complement equation to converge in about 2 iterations, but the error for the fluid pressure is doubled, and the other variables experience a loss of accuracy by about two orders of magnitude. %3.83e-8, 1.73e-2, 2.88e-3, and 5.24  -- this is for P1 for LM / pp
Contrasted with the errors in the corresponding row of Table \ref{FPSI_FEM:table:SpaceConvMixedBCs2}, we see that the solution of the Schur complement equation particularly impacts the solution of the Stokes variables.

 To visually examine the results shown in Table \ref{FPSI_FEM:table:ItersConds}, we plot the relative residual of the Schur complement equation versus the number of iterations for BiCGstab(l) in Figure \ref{FPSI_FEM:fig:ManSolution_Pre_bicg} for $\Delta x = 1/32, T =10^{-4},$ and $\Delta t = 10^{-5}$. 
As discussed, without the preconditioner, BiCGstab(l) does not converge at the original tolerance level of $10^{-8}$.
If we utilize the preconditioner without the lower right block, we see from both Figure \ref{FPSI_FEM:fig:ManSolution_Pre_bicg} and Table \ref{FPSI_FEM:table:ItersConds} that BiCGstab(l) converges in one iteration on average. Including the lower right block does not change the accuracy of the solution, but it does reduce the number of iterations to 0.5 on average.  We conclude that the use of a preconditioner for the Schur complement equation is necessary for convergence and accuracy of the algorithm, but for a manufactured solution, the lower right block of $M_{pre}$ can be left out as it only increases the computational cost without providing any significant gains in accuracy or convergence.
% Without the preconditioner, the algorithm runs in approximately 40-55 seconds, whereas inclusion of the preconditioner reduces the runtime to about 20-28 seconds.

\subsection{Hydrological Example}\label{FPSI_FEM:sec:Numerical_Hydro}

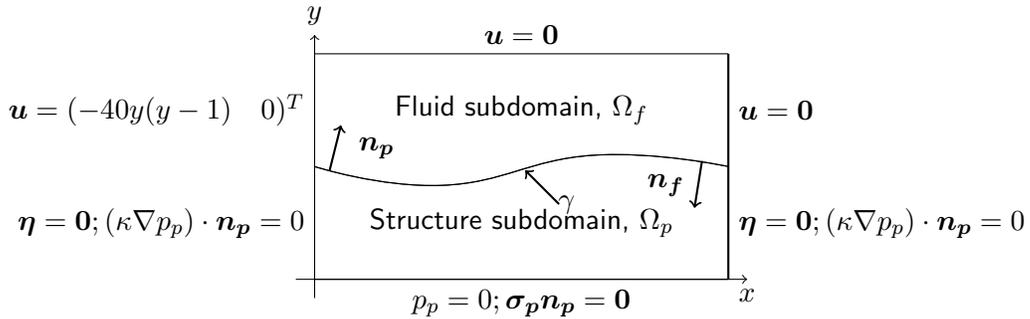
\begin{figure}[h!]
\begin{center}
\begin{tikzpicture}[scale=0.5,
    every path/.style = {},
    every node/.append style = {font=\sffamily}  ]
  \begin{scope}
   \draw[->] (-0.5,0) -- (11.5,0) node[below] {$x$};
    \draw[->] (0,-0.5) -- (0,6.5) node[above] {$y$};
\draw[name] (0,6)--(11,6);
\draw[name] (0,0)--(11,0);
\draw[name] plot [smooth,tension=0.8] coordinates{(0,3) (3.4,2.5) (7.4,3.3) (11,3)};
%\tikzfillbetween[of=A and B]{left color=blue!50, bottom color=red, right color=white, opacity=0.9};
%\tikzfillbetween[of=B and C]{left color=red!50, right color=white, opacity=0.9};
\draw[name] plot [smooth,tension=0.8] coordinates{(0,3) (3.4,2.5) (7.4,3.3) (11,3)};
\draw[thick, ->] (10.3,3.15)--(10.1,1.9);
\node[left] at (10.1,2.5) {$ \bm{n_f} $};
\draw[thick, ->] (0.4,2.9)--(0.7,4.1);
\node[right] at (0.9,3.5) {$ \bm{n_p} $};
\draw[thick, -] (11,0)--(11,6);
\node[left] at (0,1.5) {$ \bm{\eta} = \bm{0}; (\kappa \nabla p_p) \cdot \bm{n_p} = 0 $};
\node[left] at (0,4.5) {$ \bm{u} = (-40y(y-1) \quad 0)^T$};
\node[right] at (11,1.5) {$ \bm{\eta} = \bm{0}; (\kappa \nabla p_p)\cdot \bm{n_p} = 0$};
\node[right] at (11,4.5) {$ \bm{u} = \bm{0} $};
\node[below] at (5.5,0) {$ p_p = 0; \bm{\sigma_p}\bm{n_p} = \bm{0} $};
\node[above] at (5.5,6) {$ \bm{u} = \bm{0} $};
\node[] at (5.5,4.5) {Fluid subdomain, $ \Omega_f $};
\node[] at (5.5,1.5) {Structure subdomain, $ \Omega_p $};
\draw[thick, ->] (6.5,2)--(5.6,2.9);
\node[] at (6.7,2) {$ \gamma $};
  \end{scope}
\end{tikzpicture}
\caption{Fluid-poroelastic domain for hydrological example}
\label{FPSI_FEM:fig:HydrologicalDomainSetting}
\end{center} \vspace{-0.4cm}
\end{figure}

Next, we consider an example from \cite{Li_2022_Hydro} which models the coupling of a surface and subsurface hydrological system. As shown in Figure \ref{FPSI_FEM:fig:HydrologicalDomainSetting}, with $\Omega_f = [0,2]\times[0,1]$ and $\Omega_p = [0,2]\times[-1,0]$, we enforce homogeneous Dirichlet boundary conditions on the top and right boundaries of $\Omega_f$, with $\bm{u} = (-40y(y-1) \quad 0)^T$ on $x=0$. Homogeneous Dirichlet conditions are enforced for $p_p$ at $y=-1$, and for $\bm{\eta}$ on the right and left boundaries. Zero Neumann conditions are used elsewhere, and initial conditions and body forces are also set to zero. Using $\Delta t = 0.06$, $\Delta x = \Delta y = \frac{1}{32}$, we set the final time $T = 3$. In this example, we use $P_1$ elements for the LMs and Taylor Hood elements for the two-field Biot model \cite{Cesmelioglu_2020}. 
% with tilde is eta, without is p. 

\begin{figure}
    \centering
     \begin{subfigure}[b]{0.48\textwidth}   
    \centering    \includegraphics[width=\textwidth]{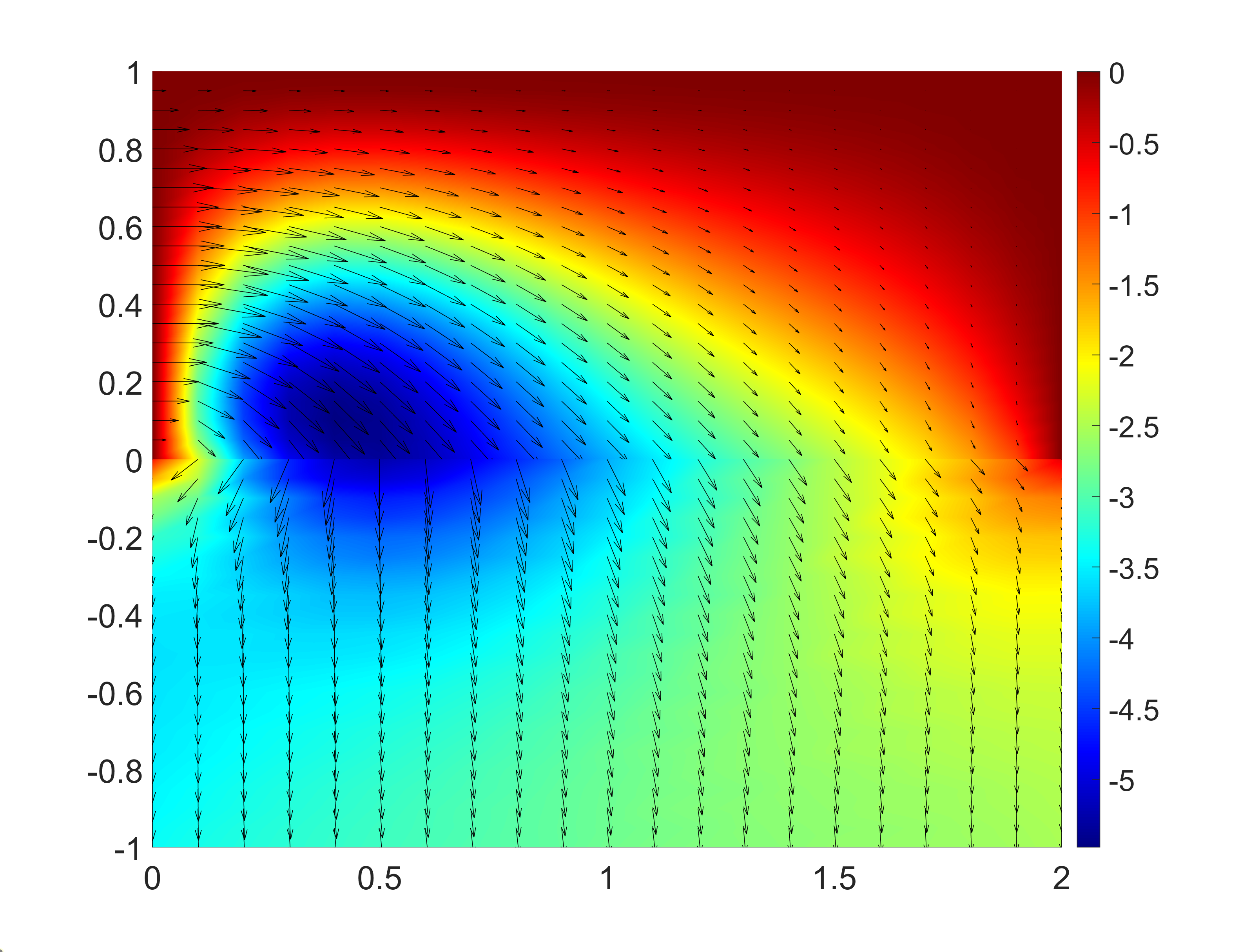}    \caption{Velocities: $\bm{u}$ and $\frac{\partial \bm{\eta}}{\partial t} - \kappa \nabla p_p$ (arrows); vertical components (color)} \label{FPSI_FEM:fig:Vel_Prob4}
    \end{subfigure}  
    \begin{subfigure}[b]{0.49\textwidth}   
    \centering    \includegraphics[width=\textwidth]{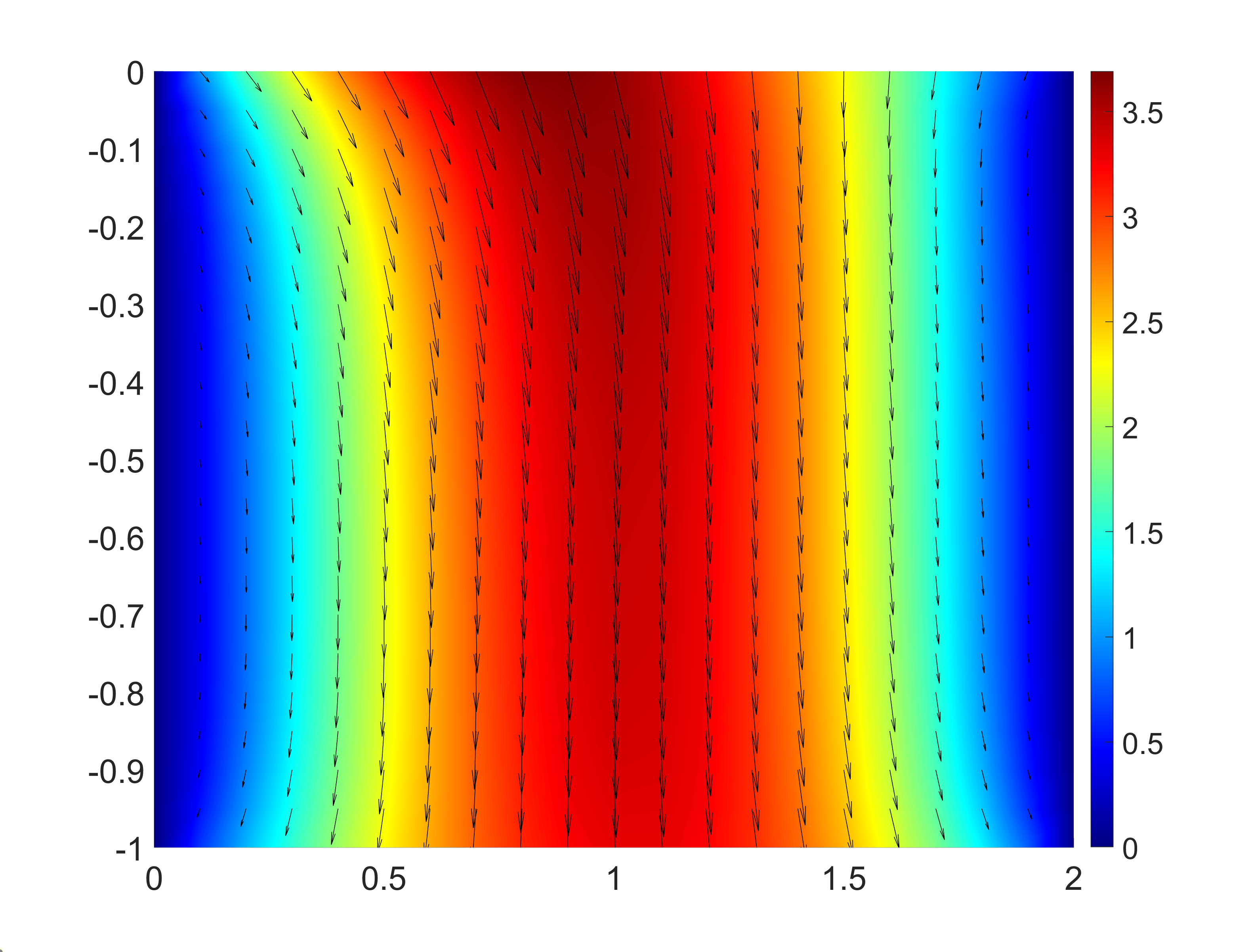}    \caption{Displacement: $\bm{\eta}$ (arrows), $|\bm{\eta}|$ (color)} \label{FPSI_FEM:fig:Disp_Prob4}
    \end{subfigure}  
    \begin{subfigure}[b]{0.48\textwidth}   
    \centering    \includegraphics[width=\textwidth]{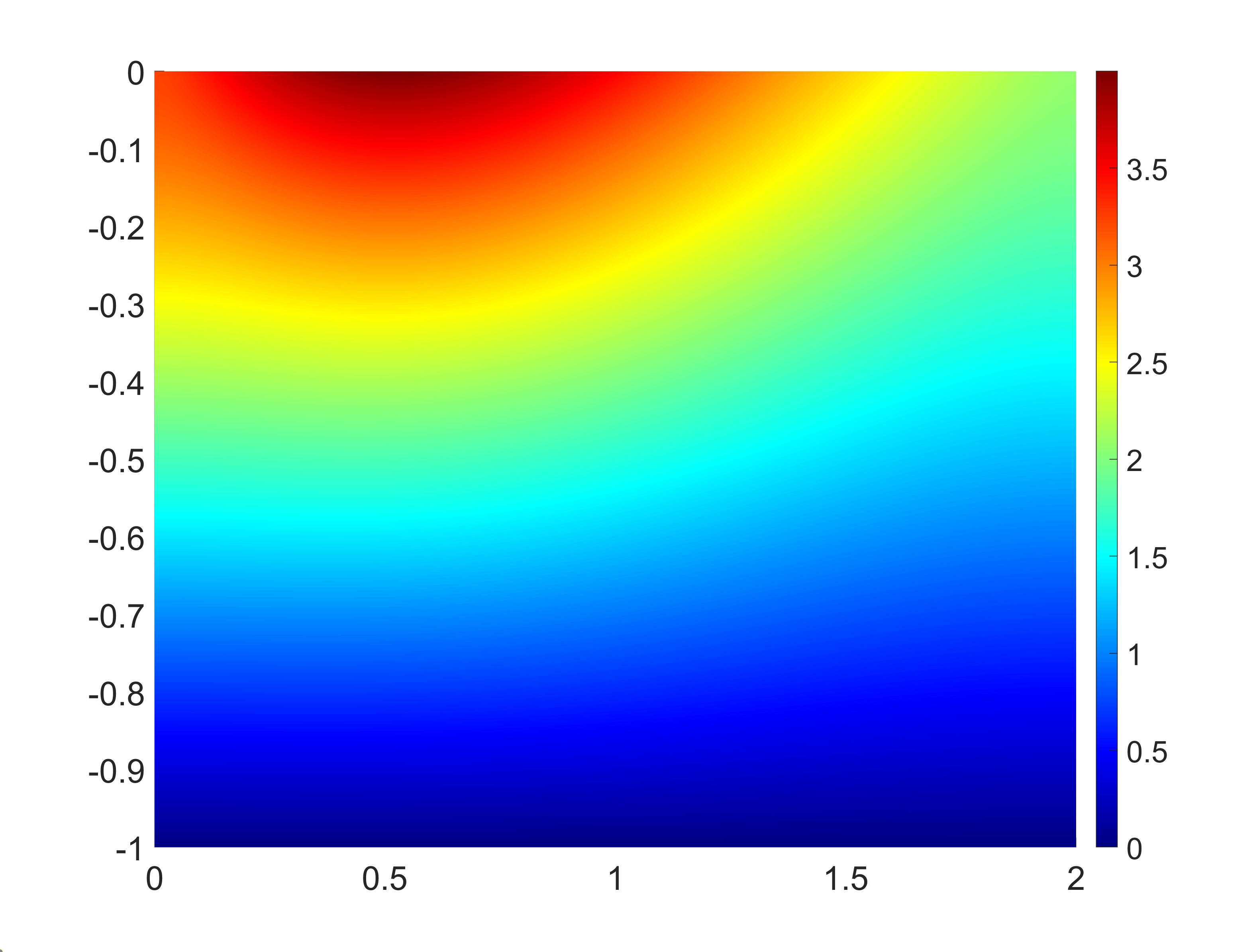}    \caption{Pore pressure $p_p$} \label{FPSI_FEM:fig:PoreP_Prob4}
    \end{subfigure}  
    \begin{subfigure}[b]{0.48\textwidth}   
    \centering    \includegraphics[width=\textwidth]{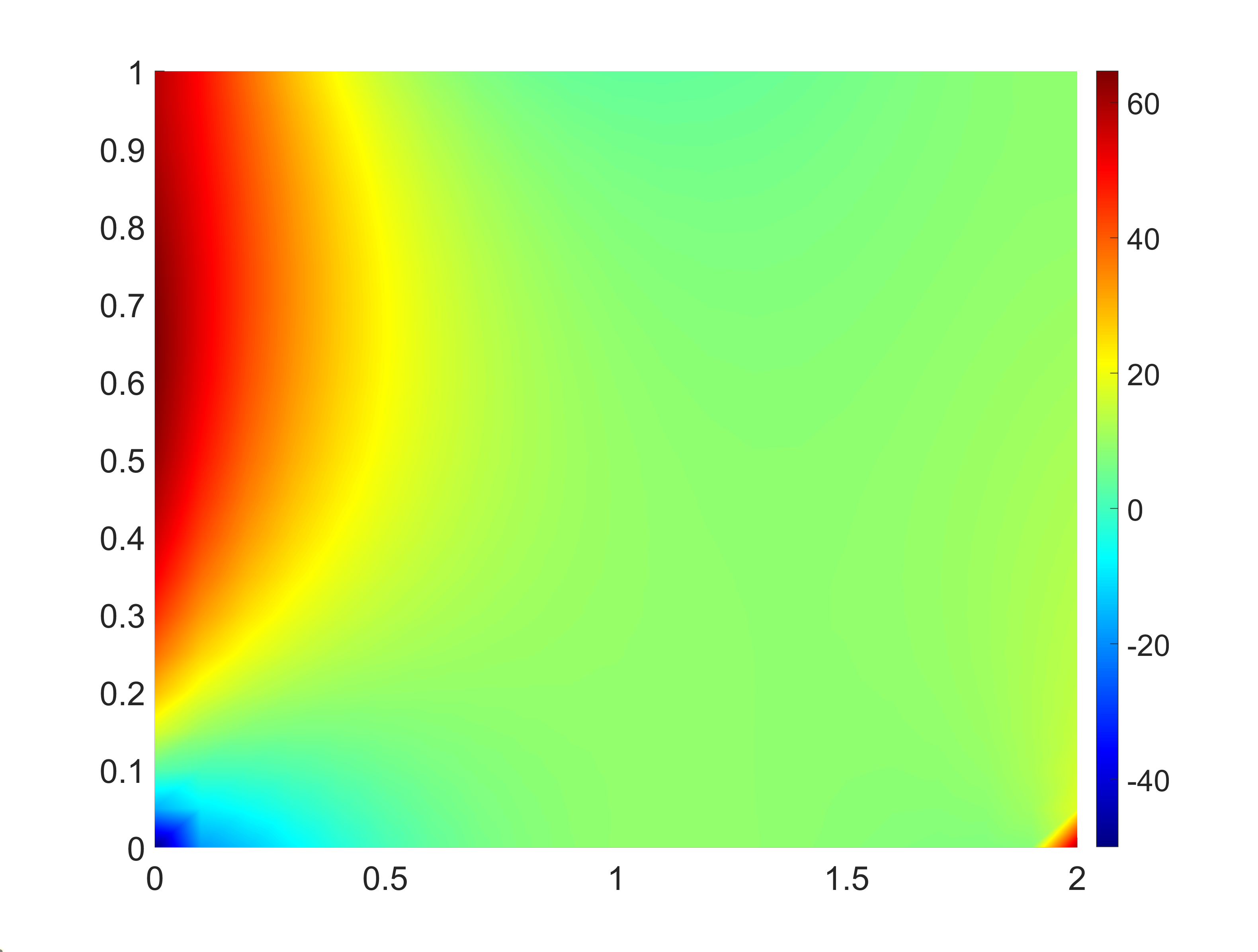}    \caption{Fluid pressure $p_f$} \label{FPSI_FEM:fig:FluidP_Prob4}
    \end{subfigure} 
    \caption[Physical plots of FF solution for case 1]{Hydrological example, case 1: all physical parameters set to 1.}
    \label{FPSI_FEM:fig:FEMResults_Prob4}
\end{figure}

We examine two cases which differ in the values of the parameters. For Case 1, all physical parameters are set to 1, which allows us to observe the qualitative behavior of the solution, shown at $T=3$ in Figure \ref{FPSI_FEM:fig:FEMResults_Prob4}. In Figure \ref{FPSI_FEM:fig:Vel_Prob4}, the arrows represent the fluid velocity in $\Omega_f$ and the vector $\frac{\partial \bm{\eta}}{\partial t} - \kappa \nabla p_p$ in $\Omega_p$ and the surface map represents the y-component of the velocity vectors. The arrows in Figure \ref{FPSI_FEM:fig:Disp_Prob4} show the magnitude $|\bm{\eta}|$ of the displacement and the arrows show the displacement itself. In particular, we observe that the velocities in $\Omega_f$ and $\Omega_p$ agree well on the interface $\gamma$.  The plots in Figure \ref{FPSI_FEM:fig:FEMResults_Prob4} may be compared with their corresponding images in Figure 1 of \cite{Li_2022_Hydro}.  We are solving the fully inertial case instead of quasi-static, so our results may differ slightly.

For Case 2, we set the parameters $\kappa = 10^{-4}, s_0 = 10^{-4}$, and $\lambda = 10^6$. Extreme values of these parameters increase the computational difficulty.
For example, small permeability ($\kappa$) and storativity constants ($s_0$) may cause the structure to behave as an incompressible material at early time steps, leading to pressure oscillations \cite{Phillips_2009}. A second source of instability, Poisson locking, may occur as $\lambda \rightarrow \infty$, observed in unreasonably small values of the displacement or oscillations in the stress. This is usually known to surface in the case of continuous linear or bilinear finite elements, as the FEM space becomes overconstrained by the requirement that $||\nabla \cdot \bm{\eta}||_1$ approaches zero as $\lambda $ approaches infinity \cite{Phillips_2009, Yi_2017}. As we use quadratic elements for the displacement, we do not expect to observe issues from Poisson locking in our results, but we recognize that the extreme values of $\kappa, s_0$, and $\lambda$ provide a scenario in which different forms of instabilities could potentially arise.

\begin{figure}
    \centering
     \begin{subfigure}[b]{0.48\textwidth}   
    \centering    \includegraphics[width=\textwidth]{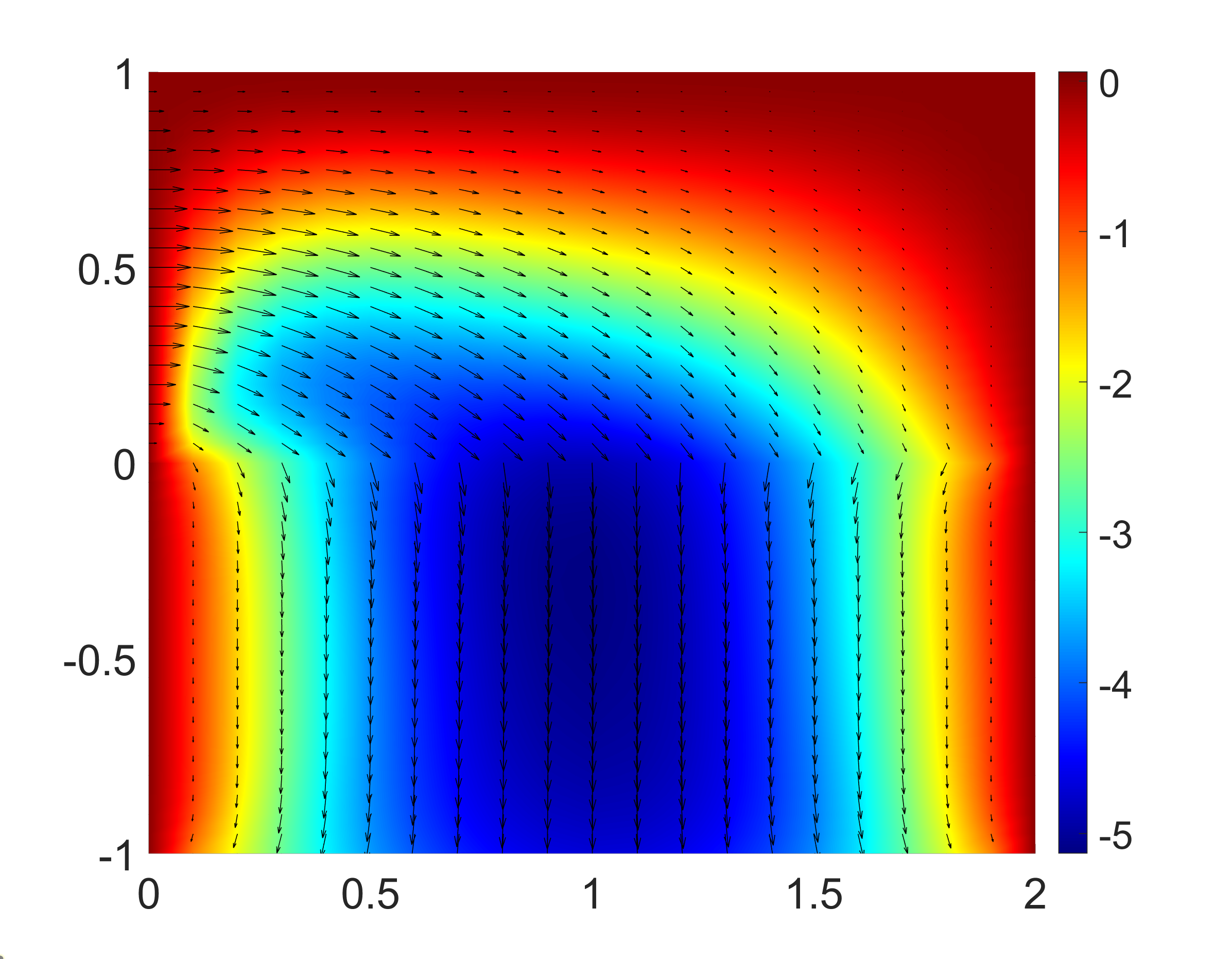}    \caption{Velocities: $\bm{u}$ and $\frac{\partial \bm{\eta}}{\partial t} - \kappa \nabla p_p$ (arrows); vertical components (color)} \label{FPSI_FEM:fig:Vel_Prob4_Case2}
    \end{subfigure}  
    \begin{subfigure}[b]{0.49\textwidth}   
    \centering    \includegraphics[width=\textwidth]{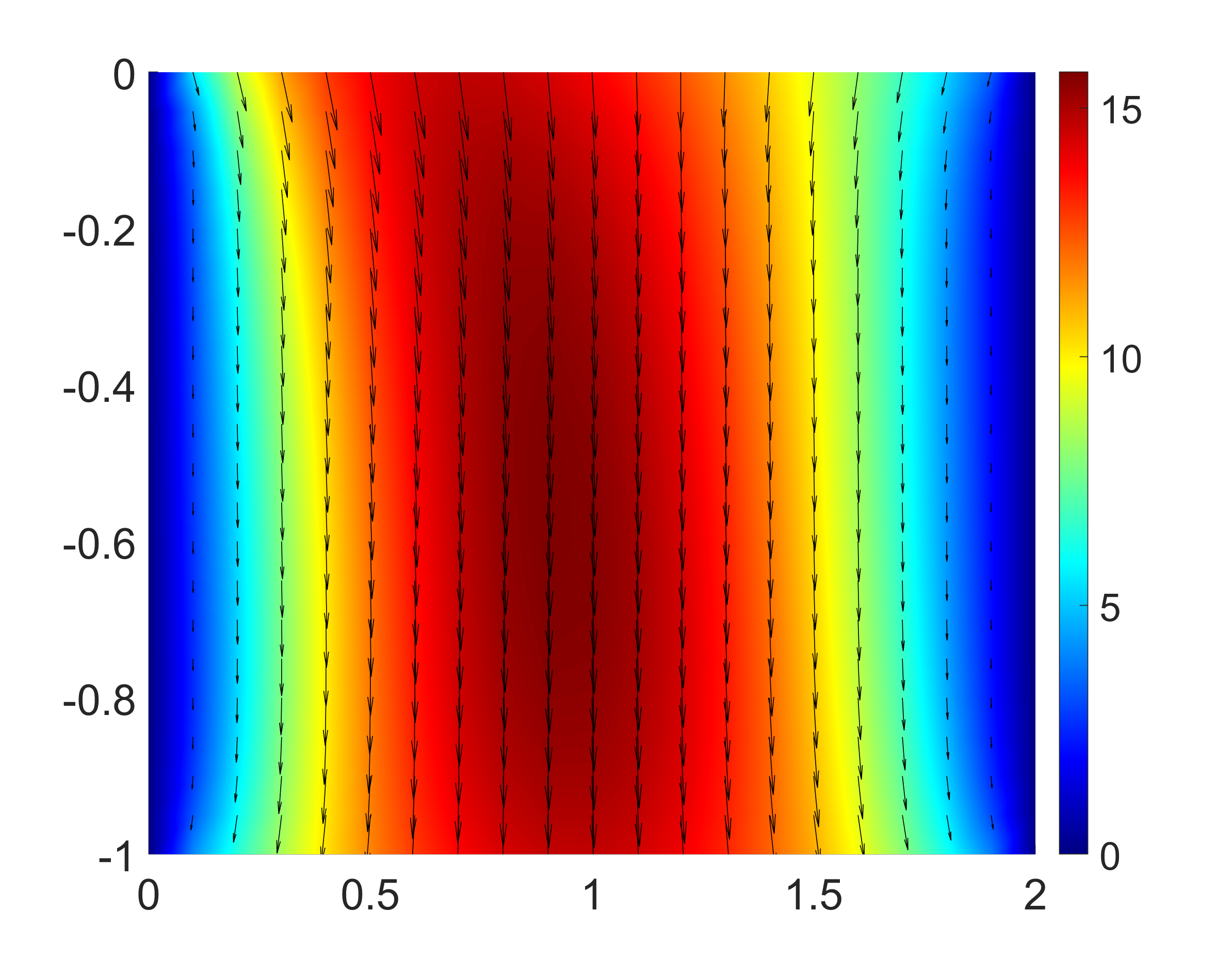}    \caption{Displacement: $\bm{\eta}$ (arrows), $|\bm{\eta}|$ (color)} \label{FPSI_FEM:fig:Disp_Prob4_Case2}
    \end{subfigure}  
    \begin{subfigure}[b]{0.48\textwidth}   
    \centering    \includegraphics[width=\textwidth]{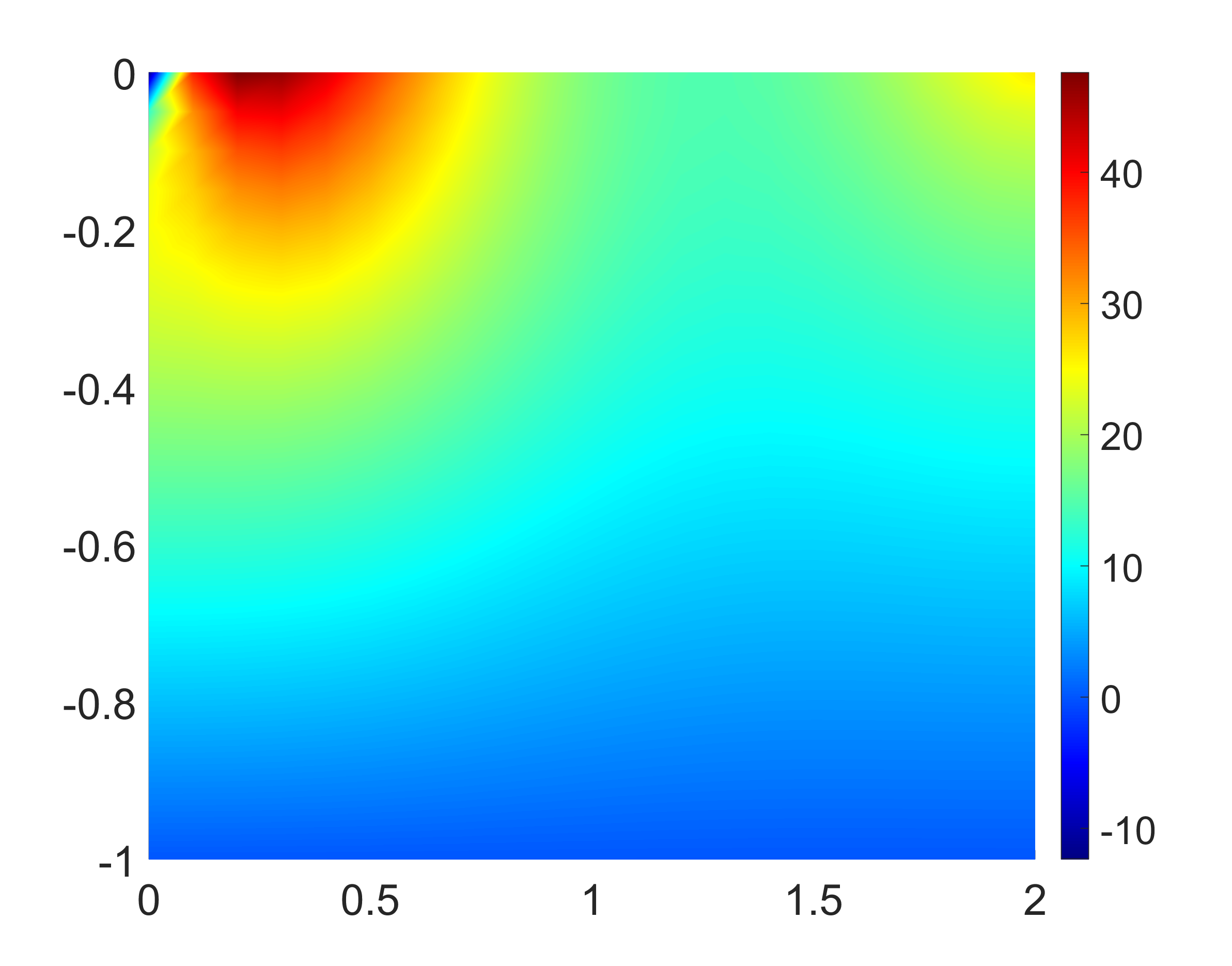}    \caption{Pore pressure $p_p$} \label{FPSI_FEM:fig:PoreP_Prob4_Case2}
    \end{subfigure}  
    \begin{subfigure}[b]{0.48\textwidth}   
    \centering    \includegraphics[width=\textwidth]{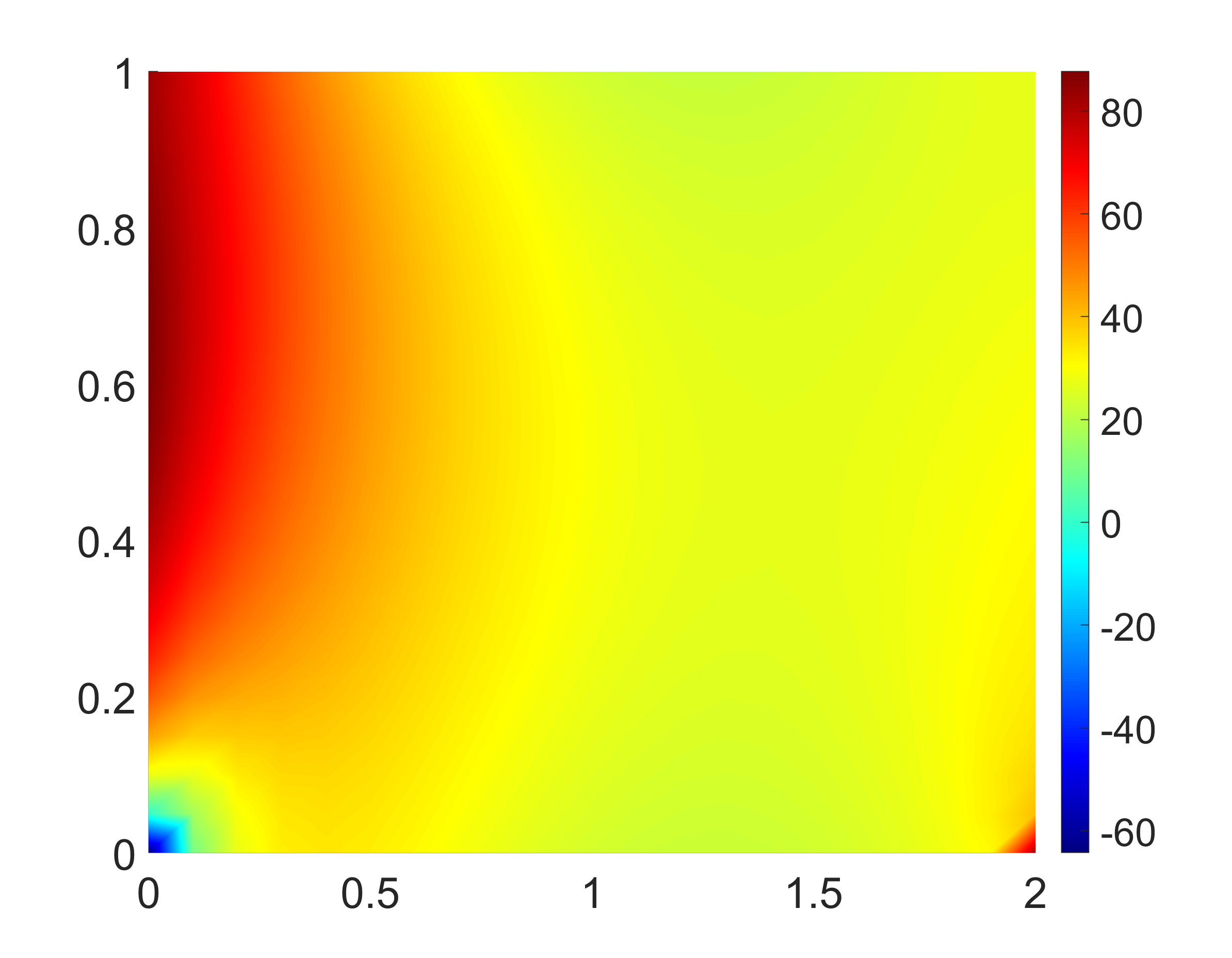}    \caption{Fluid pressure $p_f$} \label{FPSI_FEM:fig:FluidP_Prob4_Case2}
    \end{subfigure} 
    \caption[Physical plots of FF solution for case 2]{Hydrological example, case 2: $\kappa = 10^{-4}, s_0 = 10^{-4}, \lambda = 10^6$. Other parameters set to 1.}   \label{FPSI_FEM:fig:FEMResults_Prob4_Case2_FF}
\end{figure}

The physical solutions for Case 2 are plotted at the final time in Figure \ref{FPSI_FEM:fig:FEMResults_Prob4_Case2_FF}.
We do not observe spurious oscillations in the pore pressure or other effects of locking, and results align well with those obtained in Figure 2 of \cite{Li_2022_Hydro}. Qualitatively, the behavior is similar to that of Case 1, with the magnitude of the displacement and pressure fields increased.

\begin{figure}
    \centering
    \includegraphics[width=0.5\textwidth]{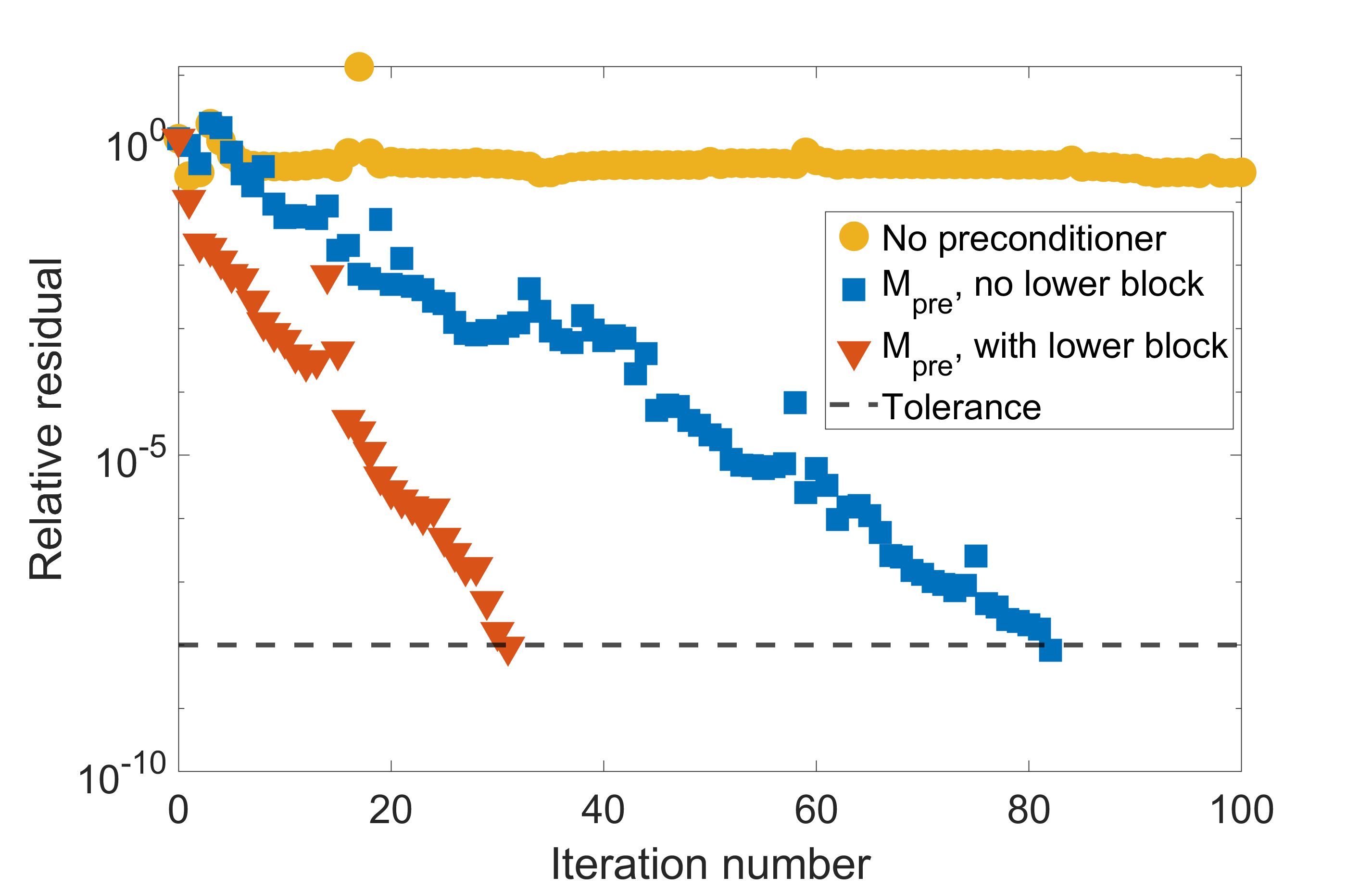}
    \caption[Residual vs. iteration number for case 2]{Hydrological example, case 2: relative residual versus iteration number for BiCGstab(l).}
    \label{FPSI_FEM:fig:P4C2_PreconditionerEffectonBiCGstab(l)}
\end{figure}

Lastly, we observe the behavior of the Schur complement equation with and without the preconditioner for this more computationally complex case. In Figure \ref{FPSI_FEM:fig:P4C2_PreconditionerEffectonBiCGstab(l)}, we see that without a preconditioner, the residual stagnates and the solver never converges. Adding a preconditioner, without the lower right block, allows the BiCGstab(l) to converge after about 80 iterations. However, inclusion of the lower right block in the preconditioner lowers the number of iterations from 80 to about 35. 

From these results, we see that the preconditioner is necessary to solve the Schur complement equation; for simpler problems, it is not necessary to include the lower right hand block in $M_{pre}$. However, as the complexity of the problem increases, the inclusion of that block does play a significant role in decreasing the number of iterations required for the Schur complement equation to converge to a solution.
%no pre: 520 seconds, pre w/o LB: 678  seconds, pre w/ LB: 334 sec

	\section{Conclusions}\label{FPSI_joint:sec:Conclusion}

We have presented a non-iterative, strongly coupled partitioned method for the solution of a fluid-poroelastic interaction system. The fully dynamic, two-field Biot model and dynamic Stokes equations are considered for the poroelastic subdomain and the fluid subdomain, respectively. Our formulation requires the use of three Lagrange multipliers representing interfacial quantities, and we performed convergence analysis for the finite element approximation of this problem.

Our partitioned method hinges on a Schur complement equation whose solution contains the Lagrange multipliers and fluid pressure at the unknown time step. Once this equation is solved, the fluid and poroelastic structure subdomains are fully decoupled and may be updated independently. If choosing to solve the Schur complement equation with an iterative linear solver, such as BiCG or related methods, we show an efficient technique for computing matrix-vector products with the Schur complement matrix $S$ and also derive a preconditioner for the system whose efficacy is demonstrated in several numerical examples. Expected rates of convergence are obtained in a manufactured solution, and the method's robustness is seen through a hydrological example with more physically and computationally challenging parameters.

%From our numerical tests thus far, our partitioned method has shown itself to be robust and it is a promising technique for future endeavors with different biomedical or geoscience applications. 

\small 
\bibliography{DissertationBib_FinalCombined}

@article{Benzi_2005,
  title={Numerical solution of saddle point problems},
  author={M. Benzi and G.H. Golub and J. Liesen},
  journal={Acta Numerica},
  volume={14},
  pages={1--137},
  year={2005},
  publisher={Cambridge University Press}
}

@article{Murad_2001,
  title={Micromechanical computational modeling of secondary consolidation and hereditary creep in soils},
  author={M.A. Murad and J.N. Guerreiro and A.F.D. Loula},
  journal={Computer Methods in Applied Mechanics and Engineering},
  volume={190},
  number={15-17},
  pages={1985--2016},
  year={2001},
  publisher={Elsevier}
}

@incollection{Detournay_1993,
  title={Fundamentals of poroelasticity},
  author={E. Detournay and A.H.-D. Cheng},
  booktitle={Analysis and Design Methods},
  pages={113--171},
  year={1993},
  publisher={Elsevier}
}

@article{Kim_2011,
  title={Stability and convergence of sequential methods for coupled flow and geomechanics: fixed--stress and fixed--strain splits},
  author={J. Kim and H.A. Tchelepi and R. Juanes},
  journal={Computer Methods in Applied Mechanics and Engineering},
  volume={200},
  number={13-16},
  pages={1591--1606},
  year={2011},
  publisher={Elsevier}
}

@article{Wheeler_2007,
  title={Iteratively coupled mixed and {G}alerkin finite element methods for poro--elasticity},
  author={M.F. Wheeler and X. Gai},
  journal={Numerical Methods for Partial Differential Equations: An International Journal},
  volume={23},
  number={4},
  pages={785--797},
  year={2007},
  publisher={Wiley Online Library}
}

@incollection{Showalter_2005,
  title={Poroelastic filtration coupled to {S}tokes flow},
  author={R.E. Showalter},
  booktitle={Control Theory of Partial Differential Equations},
  pages={243--256},
  year={2005},
  publisher={Chapman and Hall/CRC}
}

@article{Badia_2009_NSBiot,
  title={Coupling {B}iot and {N}avier--{S}tokes equations for modelling fluid--poroelastic media interaction},
  author={S. Badia and A. Quaini and A. Quarteroni},
  journal={Journal of Computational Physics},
  volume={228},
  number={21},
  pages={7986--8014},
  year={2009},
  publisher={Elsevier}
}

@article{Ambartsumyan_2018,
  title={A {L}agrange multiplier method for a {S}tokes--{B}iot fluid--poroelastic structure interaction model},
  author={I. Ambartsumyan and E. Khattatov and I. Yotov and P. Zunino},
  journal={Numerische Mathematik},
  volume={140},
  number={2},
  pages={513--553},
  year={2018},
  publisher={Springer}
}

@article{Bociu_2021,
    author = "L. Bociu and S. Canic and B. Muha and J.T. Webster",
    title = "Multilayered Poroelasticity Interacting with {S}tokes Flow",
    journal = "SIAM Journal on Mathematical Analysis",
    volume= "53",
    number="6",
    pages="6243--6279",
    year="2021",
    publisher="SIAM"}

@article{Caucao_2022,
author = "S. Caucao and T. Li and I. Yotov",
title = "A multipoint stress--flux mixed finite element method for the {S}tokes--{B}iot model",
journal = "Numerische Mathematik",
volume = "152",
year = "2022",
pages = "411--473"
}

@article{Avalos_2024,
  title="Weak and strong solutions for a fluid--poroelastic--structure interaction via a semigroup approach",
  author="G. Avalos and E. Gurvich and J.T. Webster",
  journal="Mathematical Methods in the Applied Sciences",
  year="2024",
  publisher="Wiley Online Library"
}

@article{Kunwar_2020,
  title={Second--order time discretization for a coupled quasi--{N}ewtonian fluid--poroelastic system},
  author={H. Kunwar and H. Lee and K. Seelman},
  journal={International Journal for Numerical Methods in Fluids},
  volume={92},
  number={7},
  pages={687--702},
  year={2020},
  publisher={Wiley Online Library}
}

@article{Cesmelioglu_2016,
  title={Optimization--based decoupling algorithms for a fluid--poroelastic system},
  author={A. Cesmelioglu and H. Lee and A. Quaini and K. Wang and S.-Y. Yi},
  journal={Topics in Numerical Partial Differential Equations and Scientific Computing},
  pages={79--98},
  year={2016},
  publisher={Springer}
}

@article{Cesmelioglu_2020,
  title={Numerical analysis of the coupling of free fluid with a poroelastic material},
  author={A. Cesmelioglu and P. Chidyagwai},
  journal={Numerical Methods for Partial Differential Equations},
  volume={36},
  number={3},
  pages={463--494},
  year={2020},
  publisher={Wiley Online Library}
}

@article{Wen_2020,
  title={A strongly conservative finite element method for the coupled {S}tokes--{B}iot model},
  author={J. Wen and Y. He},
  journal={Computers \& Mathematics with Applications},
  volume={80},
  number={5},
  pages={1421--1442},
  year={2020},
  publisher={Elsevier}
}

@article{Guidoboni_2009,
  title={Stable loosely--coupled--type algorithm for fluid--structure interaction in blood flow},
  author={G. Guidoboni and R. Glowinski and N. Cavallini and S. Canic},
  journal={Journal of Computational Physics},
  volume={228},
  number={18},
  pages={6916--6937},
  year={2009},
  publisher={Elsevier}
}

@article{Both_2017,
  title={Robust fixed stress splitting for {B}iot’s equations in heterogeneous media},
  author={J.W. Both and M. Borregales and J.M. Nordbotten and K. Kumar and F.A. Radu},
  journal={Applied Mathematics Letters},
  volume={68},
  pages={101--108},
  year={2017},
  publisher={Elsevier}
}

@article{Burman_2014NitscheRR,
  title={Explicit strategies for incompressible fluid--structure interaction problems: {N}itsche type mortaring versus {R}obin--{R}obin coupling},
  author={E. Burman and M.A. Fern{\'a}ndez},
  journal={International Journal for Numerical Methods in Engineering},
  volume={97},
  number={10},
  pages={739--758},
  year={2014},
  publisher={Wiley Online Library}
}

@article{Bukac_2015Nitsche,
  title={Partitioning strategies for the interaction of a fluid with a poroelastic material based on a {N}itsche’s coupling approach},
  author={M. Buka{\v{c}} and I. Yotov and R. Zakerzadeh and P. Zunino},
  journal={Computer Methods in Applied Mechanics and Engineering},
  volume={292},
  pages={138--170},
  year={2015},
  publisher={Elsevier}
}

@article{Bukac_2015OpSplit,
  title="An operator splitting approach for the interaction between a fluid and a multilayered poroelastic structure",
  author="M. Buka{\v{c}} and I. Yotov and P. Zunino",
  journal="Numerical Methods for Partial Differential Equations",
  volume="31",
  number="4",
  pages="1054--1100",
  year="2015",
  publisher="Wiley Online Library"
}

@article{Causin_2014,
  title={A poroelastic model for the perfusion of the lamina cribrosa in the optic nerve head},
  author={P. Causin and G. Guidoboni and A. Harris and D. Prada and R. Sacco and S. Terragni},
  journal={Mathematical Biosciences},
  volume={257},
  pages={33--41},
  year={2014},
  publisher={Elsevier}
}

@article{Banks_2017,
  title={Sensitivity analysis in poro--elastic and poro--visco--elastic models with respect to boundary data},
  author={H.T. Banks and K. Bekele-Maxwell and L. Bociu and M. Noorman and G. Guidoboni},
  journal={Quarterly of Applied Mathematics},
  volume={75},
  number={4},
  pages={697--735},
  year={2017},
  publisher={JSTOR}
}

@article{Calo_2008,
  title={Multiphysics model for blood flow and drug transport with application to patient--specific coronary artery flow},
  author={V.M. Calo and N.F. Brasher and Y. Bazilevs and T.J.R. Hughes},
  journal={Computational Mechanics},
  volume={43},
  pages={161--177},
  year={2008},
  publisher={Springer}
}

@article{Phillips_2009,
  title={Overcoming the problem of locking in linear elasticity and poroelasticity: an heuristic approach},
  author={P.J. Phillips and M.F. Wheeler},
  journal={Computational Geosciences},
  volume={13},
  pages={5--12},
  year={2009},
  publisher={Springer}
}

@article{Yi_2017,
  title={A study of two modes of locking in poroelasticity},
  author={S.-Y. Yi},
  journal={SIAM Journal on Numerical Analysis},
  volume={55},
  number={4},
  pages={1915--1936},
  year={2017},
  publisher={SIAM}
}

@article{Ambartsumyan_2019,
  title={A nonlinear {S}tokes--{B}iot model for the interaction of a non-{N}ewtonian fluid with poroelastic media},
  author={I. Ambartsumyan and V.J. Ervin and T. Nguyen and I. Yotov},
  journal={ESAIM: Mathematical Modelling and Numerical Analysis},
  volume={53},
  number={6},
  pages={1915--1955},
  year={2019},
  publisher={EDP Sciences}
}

@article{Cesmelioglu_2017,
  title={Analysis of the coupled {N}avier--{S}tokes/{B}iot problem},
  author={A. Cesmelioglu},
  journal={Journal of Mathematical Analysis and Applications},
  volume={456},
  number={2},
  pages={970--991},
  year={2017},
  publisher={Elsevier}
}

@article{Li_2024,
 title="An augmented fully mixed formulation for the quasistatic {N}avier--{S}tokes--{B}iot model",
  author="T. Li and S. Caucao and I. Yotov",
  journal="IMA Journal of Numerical Analysis",
  volume="44",
  number="2",
  pages="1153--1210",
  year="2024",
  publisher="Oxford University Press"
}

@article{Wilfrid_2020,
  title={Nonconforming finite element methods for a {S}tokes/{B}iot fluid--poroelastic structure interaction model},
  author={H.K. Wilfrid},
  journal={Results in Applied Mathematics},
  volume={7},
  pages={100127},
  year={2020},
  publisher={Elsevier}
}

@book{Ern_2004,
title = "Theory and {P}ractice of {F}inite {E}lements",
author = "A. Ern and J.-L. Guermond",
volume="159",
year = "2004",
publisher = "Springer-Verlag",
}

@article{Gunzburger_1992,
title = "Treating Inhomogeneous Essential Boundary Conditions in Finite Element Methods and the Calculation of Boundary Stresses",
author = "M. Gunzburger and S. Hou",
year = "1992",
volume = "29",
journal = "SIAM Journal of Numerical Analysis",
pages = "390--424"}

@misc{deCastro_2025_FPSIWP,
      title={Well--posedness of a novel {L}agrange multiplier formulation for fluid--poroelastic interaction}, 
      author={A. de Castro and H. Lee},
      year={2025},
      eprint={2512.08142},
      archivePrefix={arXiv},
      primaryClass={math.NA},
      url={https://arxiv.org/abs/2512.08142}, 
}

@book{Ciarlet_1978,
title = "The {F}inite {E}lement {M}ethod for {E}lliptic {P}roblems",
author = "P.G. Ciarlet",
year = "1978",
publisher = "North-Holland Publishing Company"}

@article{Li_2022,
  title="An Efficient {C}horin--{T}emam Projection Proper Orthogonal Decomposition Based Reduced--Order Model for Nonstationary {S}tokes Equations",
  author="X. Li and Y. Luo and M. Feng",
  journal="Journal of Scientific Computing",
  volume="93",
  pages="64",
  year="2022",
  publisher="Springer"
}

@article{Li_2022_Hydro,
  title={A mixed elasticity formulation for fluid--poroelastic structure interaction},
  author={T. Li and I. Yotov},
  journal={ESAIM: Mathematical Modelling and Numerical Analysis},
  volume={56},
  number={1},
  pages={1--40},
  year={2022},
  publisher={EDP Sciences}
}
    \bibliographystyle{plainurl}

\end{document}